\documentclass[a4paper,12pt]{amsart}

\usepackage{amssymb}
\usepackage{amsthm}
\usepackage{amscd}

\usepackage{latexsym}

\usepackage{amsmath}
\usepackage{epic,eepic}
\usepackage{graphicx}
\theoremstyle{plain}
\newtheorem{lemma}{Lemma}[subsection]
\newtheorem{prop}[lemma]{Proposition}
\newtheorem{thm}[lemma]{Theorem}
\newtheorem{cor}[lemma]{Corollary}
\newtheorem{aplemma}{Lemma~A.\hspace{-1.5mm}}
\newtheorem{approp}{Proposition~A.\hspace{-1.5mm}}
\newtheorem{apthm}{Theorem~A.\hspace{-1.5mm}}
\newtheorem{apcor}{Corollary~A.\hspace{-1.5mm}}
\newtheorem{intthm}{Theorem}

\usepackage[all]{xy}

\newtheorem{conj}[lemma]{Conjecture}

\theoremstyle{definition}

\newtheorem{rema}[lemma]{Remark}

\newtheorem{remb}{Remark}

\newtheorem{defi}[lemma]{Definition}
\newtheorem{exa}[lemma]{Example}
\newtheorem{aprem}{Remark~A.\hspace{-1.5mm}}
\newtheorem{apdefi}{Definition~A.\hspace{-1.5mm}}
\newcommand{\bde}{\begin{defi}}
\newcommand{\ede}{\end{defi}\vspace{1mm}}
\newcommand{\ble}{\begin{lemma}}
\newcommand{\ele}{\end{lemma}}
\newcommand{\bpr}{\begin{prop}}
\newcommand{\epr}{\end{prop}}
\newcommand{\bt}{\begin{thm}}
\newcommand{\et}{\end{thm}}
\newcommand{\bco}{\begin{cor}}
\newcommand{\eco}{\end{cor}}
\newcommand{\bre}{\begin{rema}}
\newcommand{\ere}{\end{rema}}
\newcommand{\brea}{\begin{rema}}
\newcommand{\erea}{\end{rema}\vspace{1mm}}
\newcommand{\breb}{\begin{remb}}
\newcommand{\ereb}{\end{remb}\vspace{1mm}}
\newcommand{\bex}{\begin{exa}}
\newcommand{\eex}{\end{exa}}
\newcommand{\bpf}{\begin{proof}}
\newcommand{\epf}{\end{proof}\vspace{1mm}}

\newcommand{\bade}{\begin{apdefi}}
\newcommand{\eade}{\end{apdefi}}
\newcommand{\bale}{\begin{aplemma}}
\newcommand{\eale}{\end{aplemma}}
\newcommand{\bapr}{\begin{approp}}
\newcommand{\eapr}{\end{approp}}
\newcommand{\bat}{\begin{apthm}}
\newcommand{\eat}{\end{apthm}}
\newcommand{\baco}{\begin{apcor}}
\newcommand{\eaco}{\end{apcor}}
\newcommand{\bare}{\begin{aprem}}
\newcommand{\eare}{\end{aprem}}

\newcommand{\be}{\begin{enumerate}}
\newcommand{\ee}{\end{enumerate}}
\newcommand{\bcd}{\[\begin{CD}}
\newcommand{\ecd}{\end{CD}\]}
\newcommand{\bit}{\begin{itemize}}
\newcommand{\eit}{\end{itemize}}
\newcommand{\bq}{\begin{quote}}
\newcommand{\eq}{\end{quote}}
\newcommand{\ba}{\begin{array}}
\newcommand{\ea}{\end{array}}
\newcommand{\mcA}{\mathcal{A}}
\newcommand{\mcB}{\mathcal{B}}

\newcommand{\mcD}{\mathcal{D}}
\newcommand{\mcE}{\mathcal{E}}
\newcommand{\mcF}{\mathcal{F}}
\newcommand{\mcG}{\mathcal{G}}
\newcommand{\mcH}{\mathcal{H}}
\newcommand{\mcI}{\mathcal{I}}

\newcommand{\mcK}{\mathcal{K}}
\newcommand{\mcL}{\mathcal{L}}
\newcommand{\mcM}{\mathcal{M}}
\newcommand{\mcN}{\mathcal{N}}
\newcommand{\mcO}{\mathcal{O}}

\newcommand{\mcS}{\mathcal{S}}
\newcommand{\mcT}{\mathcal{T}}
\newcommand{\mcU}{\mathcal{U}}
\newcommand{\mcV}{\mathcal{V}}

\newcommand{\mcZ}{\mathcal{Z}}

\newcommand{\mbB}{\mathbb{B}}
\newcommand{\mbC}{\mathbb{C}}

\newcommand{\mbF}{\mathbb{F}}
\newcommand{\mbG}{\mathbb{G}}
\newcommand{\mbH}{\mathbb{H}}

\newcommand{\mbN}{\mathbb{N}}

\newcommand{\mbP}{\mathbb{P}}
\newcommand{\mbQ}{\mathbb{Q}}
\newcommand{\mbR}{\mathbb{R}}
\newcommand{\mbS}{\mathbb{S}}
\newcommand{\mbT}{\mathbb{T}}

\newcommand{\mbV}{\mathbb{V}}
\newcommand{\mbW}{\mathbb{W}}

\newcommand{\mbZ}{\mathbb{Z}}

\newcommand{\mfC}{\mathfrak{C}}

\newcommand{\mfG}{\mathfrak{G}}

\newcommand{\mfM}{\mathfrak{M}}
\newcommand{\mfN}{\mathfrak{N}}
\newcommand{\mfO}{\mathfrak{O}}
\newcommand{\mfP}{\mathfrak{P}}

\newcommand{\mfS}{\mathfrak{S}}
\newcommand{\mfT}{\mathfrak{T}}

\newcommand{\mfX}{\mathfrak{X}}
\newcommand{\mfY}{\mathfrak{Y}}

\newcommand{\mfa}{\mathfrak{a}}
\newcommand{\mfb}{\mathfrak{b}}
\newcommand{\mfc}{\mathfrak{c}}

\newcommand{\mfg}{\mathfrak{g}}
\newcommand{\mfh}{\mathfrak{h}}
\newcommand{\mfi}{\mathfrak{i}}

\newcommand{\mfl}{\mathfrak{l}}

\newcommand{\mfn}{\mathfrak{n}}
\newcommand{\mfo}{\mathfrak{o}}
\newcommand{\mfp}{\mathfrak{p}}

\newcommand{\mfr}{\mathfrak{r}}
\newcommand{\mfs}{\mathfrak{s}}
\newcommand{\mft}{\mathfrak{t}}
\newcommand{\mfu}{\mathfrak{u}}
\newcommand{\mfv}{\mathfrak{v}}

\newcommand{\mpf}{\mathpzc{f}}
\newcommand{\mpg}{\mathpzc{g}}

\newcommand{\migi}{\rightarrow}
\newcommand{\longmigi}{\longrightarrow}

\newcommand{\isom}{\stackrel{\sim}{\migi}}

\newcommand{\migiincl}{\hookrightarrow}

\newcommand{\migisurj}{\twoheadrightarrow}

\newcommand{\mr}{\mathrm}
\newcommand{\hidden}[1]{\,}

\DeclareMathAlphabet{\mathpzc}{OT1}{pzc}{m}{it}
\newcommand{\vin}{\rotatebox{90}{$\in$}}

\begin{document}

\title[Moduli of Tango structures and dormant Miura opers]{Moduli of Tango structures \\
and dormant Miura opers}
\author{Yasuhiro Wakabayashi}
\date{}
\markboth{Yasuhiro Wakabayashi}{}
\maketitle
\footnotetext{Y. Wakabayashi: Department of Mathematics, Tokyo Institute of Technology, 2-12-1 Ookayama, Meguro-ku, Tokyo 152-8551, JAPAN;}
\footnotetext{e-mail: {\tt wkbysh@math.titech.ac.jp};}
\begin{abstract}
The purpose of the present paper is to develop the theory of (pre-)Tango structures and (dormant generic) Miura $\mathfrak{g}$-opers (for a semisimple Lie algebra $\mathfrak{g}$) defined on pointed stable curves in positive characteristic. A (pre-)Tango structure is a certain line bundle of an algebraic curve in positive characteristic, which gives some pathological (relative to zero characteristic) phenomena. In the present paper, we construct the moduli spaces of (pre-)Tango structures and (dormant generic)  Miura $\mathfrak{g}$-opers respectively and prove certain properties of them. One of the main results of the present paper states that there exists a bijective correspondence between the (pre-)Tango structures (of prescribed monodromies) and the dormant generic Miura $\mathfrak{s} \mathfrak{l}_2$-opers (of prescribed exponents). By using this correspondence, we achieve a detailed understanding of the moduli stack of (pre-)Tango structures. As an application,  we construct a family of algebraic surfaces in positive characteristic parametrized by a higher dimensional base space whose fibers are pairwise non-isomorphic and violate the Kodaira vanishing theorem.
\end{abstract}
\tableofcontents

\section*{Introduction}

The purpose of the present paper 
is to develop  the moduli  theory of (pre-)Tango structures and Miura $\mfg$-opers (for a semisimple Lie algebra $\mfg$)
 defined on (families of) pointed stable curves in positive characteristic.
One of the main results of the present paper states (cf. Theorem \ref{Td01})  that there exists 
 a bijective correspondence between the  (pre-)Tango structures (of prescribed monodromies) and 
the  dormant generic Miura $\mfs \mfl_2$-opers (of prescribed exponents).
By means of this correspondence,
we achieve a detailed understanding of the moduli stack of (pre-)Tango structures (cf. Theorem  \ref{T01}).
As an application, we construct  a family of algebraic surfaces in positive characteristic
  parametrized by a higher dimensional base space whose fibers are pairwise non-isomorphic and 
  violate the Kodaira vanishing theorem (cf. Corollary  \ref{P0ffg7}). 
In the rest of this Introduction, we shall provide more detailed discussions, including the content of the present paper.

\vspace{5mm}
\subsection*{0.1} \label{S01}

First, recall the notion of a Tango structure on an algebraic curve, which is one of the central objects   in  the present paper.
Let  $p$ be an odd prime, $k$ an algebraically closed field of characteristic $p$, and 
$X$  a proper smooth curve over $k$ of genus $g>1$.
Denote by $F_X : X \migi X$  the absolute (i.e., $p$-th power) Frobenius endomorphism of $X$.
In ~\cite{Tan}, H. Tango studied the injectivity  of the  map
\begin{align}
F_X^* : H^1 (X, \mcV) \migi H^1 (X, F_X^* (\mcV))
\end{align}
induced by  $F_X$ between  the first cohomology groups of a vector bundle $\mcV$ and its pull-back $F_X^* (\mcV)$.
In the case where  $\mcV = \mcO_X (-D)$ with  some  effective divisor $D$,   he described the kernel of the map $F^*_X$ in terms of  exact differentials, and characterized the injectivity by means of a certain numerical invariant which is now called the {\it Tango-invariant}.
The Tango-invariant is defined as 
\begin{align}
n (X) := \mr{max} \Big\{ \mr{deg}\left\lfloor  \frac{(df)}{p} \right\rfloor   \in \mbZ \ \Big|  \ f \in K (X) \setminus  K (X)^p\Big\}.
\end{align}
Here, 
 $K (X)$ denotes  the function field of $X$,   $K (X)^p$  denotes  the subfield of $K (X)$ consisting of  $p$-th powers (i.e., $K (X)^p:= \{ f^p  \ | \ f \in K (X) \}$), $(df)$ denotes the divisor $\sum_{x \in X} v_x (df) x$ (where $v_x$ denotes the valuation associated with $x$),  and $\lfloor \, -  \,  \rfloor$ denotes round down of coefficients.
One verifies  that the inequalities  $0 \leq n (X) \leq \frac{2g-2}{p}$ hold  (cf. ~\cite{Tan}, Lemma 12).
Moreover, the inequality  $n (X) > 0$ (resp., the equality $n (X) = \frac{2g-2}{p}$) implies that  there exists 
 an ample divisor $D$ on $X$ with $(df) \geq  p D$ (resp., $(df) = p D$) for   some $f \in K (X) \setminus K (X)^p$.
 We refer to the line bundle $\mcL := \mcO_X (D)$ for such a divisor $D$ as a {\it pre-Tango structure} (resp., a {\it Tango structure})  on $X$.
 (This  definition of a Tango structure  coincide with the definition   described  in Definition ~\ref{D13}.)
Then, according to  ~\cite{Tan}, Theorem 15, for each pre-Tango structure $\mcL$ on $X$,
the  map $F^* : H^1 (X, \mcL^\vee) \migi H^1 (X, F^*_{X/k} (\mcL^\vee))$
  (where $\mcL^\vee$ denotes the dual of $\mcL$)  is not injective.
Furthermore, the notion of a Tango structure has an impotence in the study of pathology of algebraic geometry in positive characteristic, as  discussed in \S\,\ref{GHTY} later.
The existence of a Tango structure  implies a strong restriction to the genus $g$ of the underlying curve, i.e., $p$ must divide $2g-2$.
At any rate, it will be natural to ask how many curves admitting a (pre-)Tango structure exist, or whether 
such curves
 are really exceptional or not.
These  questions lead us to study the moduli stack 
\begin{align} \label{E050}
\mfT \mfa \mfn_{g}
\end{align}
classifying  proper smooth  curves over $k$ of genus $g$ together with a Tango structure.
For instance, we want  to understand the image of the forgetting morphism  $\mfT \mfa \mfn_g \migi \mfM_g$, where $\mfM_g$ denotes the moduli stack classifying proper smooth curves over $k$ of genus $g$.

\vspace{5mm}
\subsection*{0.2} \label{Sff01}

In the present paper,  
we deal with  (pre-)Tango structures defined  in  a more general setting, i.e., 
 (pre-)Tango structures on  {\it families  of pointed stable curves}.
Let $(g,r)$ be a pair of nonnegative integers with $2g-2 +r >0$.
Write   $\overline{\mfM}_{g,r}$  for
   the moduli stack classifying $r$-pointed stable curves over $k$ of genus  $g$ and $\mfM_{g,r}$ for its dense open substack classifying {\it  smooth} curves.
  Suppose that we are given   an $r$-pointed genus $g$ stable curve $\mfX := (X, \{ \sigma_i \}_{i=1}^r)$  
  classified by a $k$-rational point of $\overline{\mfM}_{g,r}$, where $X$ denotes a proper nodal curve and  $\{ \sigma_i \}_{i=1}^r$ denotes an ordered set consisting of $r$ marked points in $X$.
Then, $X$ admits  a log structure  in a natural manner (cf. \S\,\ref{z041}) and we denote by $X^\mr{log}$ the resulting log scheme. 
A pre-Tango structure on $\mfX$ is defined (cf. Definition \ref{D113}) as a logarithmic  connection on the sheaf of logarithmic $1$-forms  $\Omega_{X^\mr{log}/k}$ with vanishing $p$-curvature whose horizontal sections are contained in the kernel of the Cartier operator.
If $r =0$ (i.e., $\{ \sigma_i \}_{i=1}^r= \emptyset$) and the underlying curve is smooth, then  this  definition  of a pre-Tango structure is equivalent to   the definition of a Tango structure  mentioned in the previous subsection  (cf. Proposition \ref{Pff110}).
Notice that our definition  does not require a type of  condition corresponding to the inequality $n (X) > 0$ as required  in the classical case.
But, we can  proceed to our discussion
  regardless of whether  such a condition should be imposed  or not.

Denote by $\mbF_p^{\times r}$ the product of $r$ copies of $\mbF_p := \mbZ/p\mbZ$ and let $\vec{\mu} \in \mbF_p^{\times r}$, where   $\vec{\mu} := \emptyset$ if $r =0$.
We shall write
\begin{align}
\overline{\mfT} \mfa \mfn_{g,r, \vec{\mu}}
\end{align}
(cf. (\ref{GGk})) for the moduli stack classifying $r$-pointed stable curves over $k$  of genus  $g$ together with a pre-Tango structure on it  of monodromies  $\vec{\mu}$.
In particular, we have $\mfT \mfa \mfn_g \cong  \overline{\mfT} \mfa \mfn_{g, 0, \emptyset} \times_{\overline{\mfM}_{g,0}} \mfM_{g,0}$.

\vspace{5mm}
\subsection*{0.3} \label{Sgtff01}

On the other hand, we recall the notion of a Miura $\mfg$-oper for a semismple Lie algebra $\mfg$.
A {\it Miura $\mfg$-oper}  is
a $\mfg$-oper equipped with an additional data, that is to say,
a   $\mbG$-torsor (where $\mbG$ denotes the identity component of the group of Lie algebra automorphisms of $\mfg$)  equipped with two  Borel reductions and  a flat connection satisfying some condition, including a certain transversality condition. 
For example, 
$\mfs \mfl_2$-opers
and Miura $\mfs \mfl_2$-opers
 may be identified with  projective and affine  connections respectively.
The Miura opers (over the field of complex numbers $\mbC$) plays an essential role  in integrable systems and representation theory of loop algebra, including
 the Drinfeld-Sokolov reduction and geometric Langlands correspondence (cf. ~\cite{DS}, ~\cite{FFR}, and ~\cite{Fr}).
The solutions of the Bethe Ansatz equations may be described
by means of Miura opers. (cf. ~\cite{Fr2}, ~\cite{Fr3}).

\vspace{5mm}
\subsection*{0.4} \label{Fgf01}

Let $k$, $\mfg$, and $\mbG$ be as above, and suppose that either one of the   two  conditions
  $(\mr{Char})_{p}$,  $(\text{Char})_{p}^{\mfs \mfl}$ described in \S\,\ref{SSQ1} is satisfied (in particular, $k$ has characteristic $p>0$).
Let $\mft$ be the Lie algebra of a split maximal torus of $\mbG$, and let
$\vec{\varepsilon} \in \mft (k)^{\times r}$, where $\vec{\varepsilon} := \emptyset$ if $r =0$.
Denote by
$\mfM \overline{\mfO} \mfp_{\mfg, g,r, \vec{\varepsilon}}$ (cf. (\ref{ee200}))
the moduli stack classifying $r$-pointed stable curves over $k$ of genus $g$ paired   with a generic Miura $\mfg$-oper on it of exponents $\vec{\varepsilon}$.
As shown in Proposition \ref{eqeq88}, it may be represented by a  Deligne-Mumford stack  over $k$.
Although  there is no generic Miura $\mfg$-oper on a proper smooth curve over $\mbC$ of genus $g>1$,
 this  moduli stack may not be empty.
One may find  the locus 
\begin{align}
\mfM \overline{\mfO} \mfp^{^\mr{Zzz...}}_{\mfg, g,r, \vec{\varepsilon}}
\end{align}
(cf. (\ref{eq35gg})) of $\mfM \overline{\mfO} \mfp_{\mfg, g,r, \vec{\varepsilon}}$  classifying generic  Miura $\mfg$-opers satisfying a nice condition  (regarding  $p$-curvature) called   {\it dormant generic  Miura $\mfg$-opers}.
It follows from Theorem  ~\ref{th013} (i) and (ii)  that
$\mfM \overline{\mfO} \mfp^{^\mr{Zzz...}}_{\mfg, g,r, \vec{\varepsilon}}$ is empty unless $\vec{\varepsilon} \in \mft (\mbF_p)^{\times r}$ or $\vec{\varepsilon} = \emptyset$, and  is finite over  $\overline{\mfM}_{g,r}$.
If, moreover, $\mfg = \mfs \mfl_2$, then one may obtain the following assertion, by which {\it a dormant generic Miura $\mfg$-oper  may be thought of as a generalization of a (pre-)Tango structure}.

\vspace{3mm}
\begin{intthm}[cf. Theorem  \ref{Td01}, Theorem  \ref{T01}] \label{ThD}
\leavevmode\\
 \ \ \ 
Let $\vec{\varepsilon} := (\varepsilon_i)_{i=1}^r \in \mbF_p^{\times r}$, where  $\vec{\varepsilon} := \emptyset$ if $r =0$.
\begin{itemize}
\item[(i)]
Both $\overline{\mfT} \mfa \mfn_{g,r, -\vec{\varepsilon}}$  and $\mfM \overline{\mfO} \mfp_{\mfs \mfl_2, g,r, [\vec{\varepsilon}\,]}^{^\mr{Zzz...}}$ (cf. Theorem \ref{Td01} and (\ref{ee480}) for the definitions  of $- \vec{\varepsilon}$ and $[\vec{\varepsilon}\,]$ respectively)
 may be represented  by  a (possibly empty)  smooth proper Deligne-Mumford stack  over $k$.
Also, there exists a canonical isomorphism
\begin{align}
\overline{\mfT} \mfa \mfn_{g,r, -\vec{\varepsilon}} \isom \mfM \overline{\mfO} \mfp_{\mfs \mfl_2, g,r, [\vec{\varepsilon}\,]}^{^\mr{Zzz...}}
\end{align}
 over $\overline{\mfM}_{g,r}$. 
\item[(ii)]
Write
\begin{align}
\mfM \mfO \mfp^{^\mr{Zzz...}}_{\mfs \mfl_2, g,r, [\vec{\varepsilon}\,]} := \mfM \overline{\mfO} \mfp^{^\mr{Zzz...}}_{\mfs \mfl_2, g,r,  [\vec{\varepsilon}\,]} \times_{\overline{\mfM}_{g,r}} \mfM_{g,r},
\end{align}
and let $\tau$ denote the natural bijection $\{0,1, \cdots, p-1 \} \isom \mbF_p$.
Then, $\mfM \mfO \mfp^{^\mr{Zzz...}}_{\mfs \mfl_2, g,r, [\vec{\varepsilon}\,]}$ is empty if $2g-2 + \frac{2g-2 + r +  \sum_{i=1}^r \tau^{-1}(\varepsilon_i)}{p} < 0$.
On the other hand, if  $2g-2 + \frac{2g-2 + r +  \sum_{i=1}^r \tau^{-1}(\varepsilon_i)}{p} \in \mbZ_{\geq  0}$, then any irreducible component $\mcN$ of $\mfM \overline{\mfO} \mfp_{\mfs \mfl_2, g,r, [\vec{\varepsilon}\,]}^{^\mr{Zzz...}}$ ($\cong \overline{\mfT} \mfa \mfn_{g,r, -\vec{\varepsilon}}$)  with $\mcN \times_{\overline{\mfM}_{g,r}} \mfM_{g,r} \neq \emptyset$ is  equidimensional of dimension $2g-2 + \frac{2g-2 +r + \sum_{i=1}^r \tau^{-1} (\varepsilon_i)}{p}$. 
\end{itemize}
 \end{intthm}
\vspace{5mm}

Notice that if $r =0$ and $g = \frac{(lp-1)(lp-2)}{2}$ for some integer $l$ with $lp \geq 4$, then $\mfM \mfO \mfp_{\mfs \mfl_2, g, 0, \emptyset}^{^\mr{Zzz...}}$ ($\cong \mfT \mfa \mfn_{g,0, \emptyset}$) is nonempty (cf. Remark \ref{RRR040} (i)).
As a corollary of the above theorem,  one may conclude the following assertion concerning the structure of the moduli stack $\mfT \mfa \mfn_g$.

\vspace{3mm}
\begin{intthm} \label{ThA}
\leavevmode\\
 \ \ \
If $\mfT \mfa \mfn_{g} \neq \emptyset$ (e.g.,  $g = \frac{(lp-1)(lp-2)}{2}$ for some integer $l$ with $lp \geq 4$), then $p | (g-1)$ and 
$\mfT \mfa \mfn_{g}$
 may be represented by
  an  equidimensional  smooth Deligne-Mumford stack over $k$ of dimension 
$2g-2 + \frac{2g-2}{p}$ which is finite over $\mfM_g$.

  \end{intthm}
\vspace{5mm}

In particular,
the locus of $\mfM_g$ classifying proper smooth curves admitting a Tango structure (i.e., the scheme-theoretic image of the projection $\mfT \mfa \mfn \mfg_g \migi \mfM_g$) forms, if it is nonempty,  an equidimensional closed substack of dimension 
$2g-2 + \frac{2g-2}{p}$.

\vspace{5mm}
\subsection*{0.5} \label{S0444}


In the last section of the present paper, we study 
the  pathology of algebraic geometry in positive characteristic, which  is of certain interest, since pathology reveals some completely different geometric phenomena from those in complex geometry.
It is well-known that the Kodaira vanishing theorem does not always hold if the characteristic of the base field  is positive.
 M. Raynaud has given (in ~\cite{Ray}) its counterexamples on smooth algebraic surfaces in  positive characteristic.
He constructed a smooth polarized surface $(X, \mcZ)$ with $H^1 (X, \mcZ^\vee) \neq 0$.
S. Mukai (cf. ~\cite{Muk}) generalized Raynaud's construction to obtain polarized smooth projective varieties $(X, \mcZ)$ of any dimension with $H^1 (X, \mcZ^\vee) \neq 0$.
The construction similar to Mukai's construction has been also studied by Y. Takeda (cf. ~\cite{Tak1}, ~\cite{Tak2}) and P. Russell (cf. ~\cite{Rus}).
The key ingredient of these constructions
   is the use of a Tango structure  on an algebraic curve.
By applying one of these construction and Theorem \ref{ThA} above, we obtain, as described in Theorem C below,  a family of algebraic varieties violating the  Kodaira vanishing theorem parametrized by a higher dimensional variety.
This result may be thought of as a refinement of the result given in ~\cite{Tak3}, Theorem 4.1 (or \cite{Tsu}.)

\vspace{3mm}
\begin{intthm}[cf. Corollary \ref{P0ffg7}] \label{ThBf} 
\leavevmode\\
 \ \ \ 
 Suppose that $p >3$, $p (p-1) | 2g-2$, and  $4 | p-3$.
Then,
there exists a flat family $\mfY \migi \mfT$ of proper smooth algebraic surfaces of general type parametrized by a Deligne-Mumford stack $\mfT$ over $k$ of dimension $\geq g-2 + \frac{2g-2}{p-1}$ all of whose  fibers are  pairwise non-isomorphic and have the automorphism group schemes being   non-reduced.
 \end{intthm}

\vspace{5mm}
\hspace{-4mm}{\bf Acknowledgement} \leavevmode\\
 \ \ \ 
 The author cannot  express enough his sincere and deep gratitude to Professors Shinichi Mochizuki and Kirti Joshi.
Without their philosophies and amazing insights, his study of mathematics would have remained ``{\it dormant}".
The author would like to thank
  all those who gave  the opportunity or   impart  the great joy of  studying 
  mathematics to him; the author wrote the present paper like a gratitude letter to them.
Also, the author would like to thank the referee for reading carefully his manuscript.
Finally,  special thanks go to the moduli stack  of dormant opers, who has guided him to  the beautiful world of mathematics.
The author was partially  supported by 
 the FMSP program at the Graduate School of Mathematical Sciences of the University of Tokyo, and the Grant-in-Aid for Scientific Research (KAKENHI No. 18K13385).

\vspace{10mm}
\section{Preliminaries} \vspace{3mm}

In this section, we recall some  definitions and notation
concerning  the theory of logarithmic connections on a torsor defined over  a log-curve. Basic references for the notion of a log scheme (or,  more generally, a log stack)  are ~\cite{KaKa}, ~\cite{ILL2}, and ~\cite{KaFu}.

Throughout the present paper, we shall fix a perfect field $k$ and  a pair of nonnegative integers $(g,r)$ with $2g-2 +r >0$. 
For each stack $\mfM$ over $k$, we shall denote by
$\mfS \mfc \mfh_{/\mfM}$ the category of $k$-schemes $S$ equipped with a $k$-morphism $S \migi \mfM$. 
For a log  stack (resp., a morphism of log stacks) indicated, say,  by $Y^\mr{log}$ (resp., $f^\mr{log} : Y^\mr{log} \migi Z^\mr{log}$), we shall write $Y$ (resp., $f : Y \migi Z$) for the underlying stack of $Y^\mr{log}$ (resp., the underlying morphism of stacks of $f^\mr{log}$).
If $Y^\mr{log}$ is as above and $Z^\mr{log}$ is a  log stack over $Y^\mr{log}$,
then we shall write
$\Omega_{Z^\mr{log}/Y^\mr{log}}$ for the sheaf of  logarithmic $1$-forms on $Z^\mr{log}$ over $Y^\mr{log}$ and write $\mcT_{Z^\mr{log}/Y^\mr{log}} := \Omega_{Z^\mr{log}/Y^\mr{log}}^\vee$, i.e., its dual. 
Also, we shall write $d$ for the universal (logarithmic) derivation $\mcO_Z \migi \Omega_{Z^\mr{log}/Y^\mr{log}}$.

\vspace{5mm}
\subsection{Log-curves} \label{z041}
\leavevmode\\ \vspace{-4mm}

Let $T^\mr{log}$ be an fs log scheme
over $k$.
A {\bf log-curve} over $T^\mr{log}$  (cf. ~\cite{ACGH},  Definition 4.5) is, by definition,  a log smooth integrable morphism $f^\mr{log} : U^\mr{log} \migi T^\mr{log}$ of fs log schemes
 such that the geometric fibers of the underlying morphism  of schemes $f : U \migi T$  are reduced and  connected $1$-dimensional schemes.
 In particular, both $\Omega_{U^\mr{log}/T^\mr{log}}$ and $\mcT_{U^\mr{log}/T^\mr{log}}$ are line bundles.

Denote by 
\begin{align}
\overline{\mfM}_{g,r}
\end{align}
 the moduli stack of $r$-pointed stable curves (cf. ~\cite{Kn2}, Definition 1.1) over $k$ of genus $g$, and by 
 \begin{align}
 \overline{f}_\mr{\mft \mfa \mfu} : \overline{\mfC}_{g,r} \migi \overline{\mfM}_{g,r}
 \end{align}
  the tautological curve over $\overline{\mfM}_{g,r}$, equipped with its $r$ marked points $\sigma_{\mft \mfa \mfu, 1} , \cdots, \sigma_{\mft \mfa \mfu, r} : \overline{\mfM}_{g,r} \migi \overline{\mfC}_{g,r}$.
Recall (cf. ~\cite{Kn2}, Corollary 2.6 and Theorem 2.7; ~\cite{DM}, \S\,5) that $\overline{\mfM}_{g,r}$ may be represented by a geometrically connected, proper, and smooth Deligne-Mumford stack over $k$ of dimension $3g-3+r$.
Denote by $\mcD_{g,r}$  the divisor of  $\overline{\mfM}_{g,r}$ at infinity.
In particular, its complement 
\begin{align}
\mfM_{g,r} := \overline{\mfM}_{g,r} \setminus \mcD_{g,r}
\end{align}
 in $\overline{\mfM}_{g,r}$ classifies smooth curves; we denote the tautological  smooth curve over $\mfM_{g,r}$ by
\begin{align}
f_{\mft \mfa \mfu} : \mfC_{g,r} \migi  \mfM_{g,r}
\end{align}
(i.e., $\mfC_{g,r} := \overline{\mfC}_{g,r} \times_{\overline{\mfM}_{g,r}} \mfM_{g,r}$).
 $\overline{\mfM}_{g,r}$ has a natural log structure given by $\mcD_{g,r}$ (cf. ~\cite{KaFu}, Theorem 4.5), 
where we shall denote the resulting log stack  by $\overline{\mfM}_{g,r}^{\mr{log}}$.
Also, by taking the divisor which is the union of the $\sigma_{\mft \mfa \mfu, i}$'s and the pull-back of  $\mcD_{g,r}$, we obtain a log structure on $\overline{\mfC}_{g,r}$; we denote  the resulting log stack by $\overline{\mfC}^{\mr{log}}_{g,r}$.
The structure morphism $\overline{f}_{\mft \mfa \mfu} : \overline{\mfC}_{g,r} \migi \overline{\mfM}_{g,r}$ extends naturally to a morphism 
$\overline{f}_{\mft \mfa \mfu}^\mr{log} : \overline{\mfC}^\mr{log}_{g,r} \migi \overline{\mfM}^\mr{log}_{g,r}$ of log stacks.

Next,  let $S$ be a scheme over $k$, or more generally, a stack  over $k$.
Also, let 
\begin{equation}
\label{X}
 \mfX : =(f :X \migi S, \{ \sigma_i : S \migi X\}_{i=1}^r)\end{equation}
 be an $r$-pointed stable curve over $S$ of genus $g$, where
 $f : X \migi S$ denotes 
  a proper nodal  curve  over $S$ of genus $g$ and $\{ \sigma_i \}_{i=1}^r$ denotes an ordered set of   $r$ marked points in $X$.
Then, $\mfX$ determines its classifying morphism $c_\mfX : S \migi \overline{\mfM}_{g,r}$ and an isomorphism $X \isom S \times_{c_\mfX, \overline{\mfM}_{g,r}, \overline{f}_{\mft \mfa \mfu}} \overline{\mfC}_{g,r}$ over $S$.
By pulling-back the  log structures of $\overline{\mfM}^{\mr{log}}_{g,r}$ and $\overline{\mfC}^{\mr{log}}_{g,r}$,
we obtain log structures on $S$ and $X$ respectively; we denote  the resulting log stacks by $S^{\mfX \text{-}\mr{log}}$ and $X^{\mfX \text{-}\mr{log}}$
respectively.
If there is no fear of causing confusion, we write $S^\mr{log}$ and $X^\mr{log}$ instead of $S^{\mfX \text{-}\mr{log}}$ and  $X^{\mfX \text{-}\mr{log}}$ respectively.
The structure morphism $f : X \migi S$ extends to a morphism $f^\mr{log} : X^{\mr{log}} \migi S^{\mr{log}}$ of log stacks, 
by which $X^\mr{log}$ determines a log-curve over $S^\mr{log}$ (cf. ~\cite{KaKa}, \S\,3;  ~\cite{KaFu}, Theorem 2.6).
For each $i = 1, \cdots, r$, 
there exists a canonical isomorphism  (i.e., the so-called {\it residue map})
\begin{align} \label{f029}
\mft \mfr \mfi \mfv_{\mfX, i} : \sigma^*_i (\Omega_{X^\mr{log}/S^\mr{log}})\isom \mcO_S,
\end{align}
(cf. ~\cite{Wak5}, \S\,1.6, (80)) which maps any local section of the form $\sigma_i^{*}(d\mr{log}(x)) \in \sigma_i^{*}(\Omega_{X^\mr{log}/S^{\mr{log}}})$  (for   a local function  $x$ defining 
the closed subscheme $\sigma_i : S \migi X$ of $X$) to  $1 \in \mcO_{S}$.

If, moreover, $k$ has characteristic $p >0$, then we shall  denote by  $\mfM_{g,r}^\mr{ord}$  the locus in $\mfM_{g,r}$ classifying  pointed proper smooth curves $\mfX := (X/S, \{ \sigma_i \}_{i=1}^r)$ such that the underlying curve $X/S$ is  {\it ordinary}  (i.e., the $p$-rank of its Jacobian is maximal).
It is well-known that $\mfM_{g,r}^\mr{ord}$ forms a dense open substack of $\mfM_{g,r}$.

\vspace{5mm}
\subsection{Logarithmic connections on a vector bundle} \label{z041e}
\leavevmode\\ \vspace{-4mm}

Let  $T^\mr{log}$ be an fs log scheme over $k$ and 
 $f^\mr{log} : U^\mr{log} \migi T^\mr{log}$    a log-curve over $T^\mr{log}$.
Also, let $\mcV$ be a vector bundle on $U$ (i.e., a locally free $\mcO_U$-module of finite rank).
By a {\bf $T^\mr{log}$-connection} on $\mcV$, we mean  an  
$f^{-1} (\mcO_T)$-linear morphism $\nabla_\mcV : \mcV \migi \Omega_{U^\mr{log}/ T^\mr{log}} \otimes \mcV$ satisfying the Leibniz rule:  $\nabla_\mcV (a \cdot  v) = d a \otimes v + a \cdot \nabla_\mcV (v)$, where $a$ and $v$ denote  any local sections of $\mcO_U$ and  $\mcV$ respectively.
A {\bf log flat bundle} on $U^\mr{log}/T^\mr{log}$  is a pair 
\begin{align}
\mcV^\flat := (\mcV, \nabla_\mcV)
\end{align}
 consisting of a vector bundle $\mcV$ on $U$ and a $T^\mr{log}$-connection $\nabla_\mcV$ on $\mcV$.
A {\bf log flat line  bundle} on $U^\mr{log}/T^\mr{log}$ is a log flat bundle $\mcV^\flat := (\mcV, \nabla_\mcV)$ 
such that 
 $\mcV$  is 
 of rank one.
We shall write
\begin{align} \label{ee203}
\mcO_U^\flat := (\mcO_U, d : \mcO_U \migi \Omega_{U^\mr{log}/T^\mr{log}}).
\end{align}

Next, let $\mcV^\flat := (\mcV, \nabla_\mcV)$ and ${\mcV'}^{\flat} := (\mcV', \nabla_{\mcV'})$ be log flat bundles on $U^\mr{log}/T^\mr{log}$.
An {\bf isomorphism of log flat bundles} from $\mcV^\flat$ to ${\mcV'}^\flat$  is an isomorphism of vector bundles $\mcV \isom \mcV'$ that is compatible with the respective $T^\mr{log}$-connections $\nabla_\mcV$ and $\nabla_{\mcV'}$.
If $\nabla_\mcV \otimes \nabla_{\mcV'}$ denotes the $T^\mr{log}$-connection  on the tensor product  $\mcV \otimes \mcV'$ induced  from $\nabla_\mcV$ and $\nabla_{\mcV'}$ (i.e., given by $v \otimes v' \mapsto \nabla_\mcV (v) \otimes v' + v \otimes \nabla_{\mcV'}(v')$),
then we shall write $\mcV^\flat \otimes {\mcV'}^\flat$ for the tensor product of $\mcV^\flat$ and ${\mcV'}^\flat$, i.e., the log flat bundle 
\begin{align} \label{EQ1}
\mcV^\flat \otimes {\mcV'}^\flat := (\mcV \otimes \mcV', \nabla_\mcV \otimes \nabla_{\mcV'})
\end{align}
on $U^\mr{log}/T^\mr{log}$.

\vspace{3mm}
\bde  \label{j010}\leavevmode\\
 \ \ \ 
Let $\mcV^\flat := (\mcV, \nabla_\mcV)$ and ${\mcV'}^\flat := (\mcV', \nabla_{\mcV'})$ be log flat bundles on $U^\mr{log}/T^\mr{log}$.
We shall say that $\mcV^\flat$ is {\bf $\mbG_m$-equivalent} to ${\mcV'}^\flat$ if 
there exists a log flat line bundle $\mcL^\flat := (\mcL, \nabla_\mcL)$ on $U^\mr{log}/T^\mr{log}$ such that
$\mcV^\flat \otimes \mcL^\flat$ is isomorphic to ${\mcV'}^\flat$.
\ede

\vspace{5mm}
\subsection{Logarithmic connections on a torsor} \label{09}
\leavevmode\\ \vspace{-4mm}

Let $U^\mr{log}/T^\mr{log}$ be as above,   $\mbG$   a connected smooth algebraic group over $k$, 
and  $\pi : \mcE \migi U$ be a right $\mbG$-torsor over $U$ in the  \'{e}tale 
 topology.
If $\mfh$ is a $k$-vector space equipped with a left $\mbG$-action, then we shall write $\mfh_{\mcE}$ for the vector bundle on $U$ associated with the relative affine space $\mcE \times^\mbG \mfh$ 
($:= (\mcE \times_k \mfh) /\mbG$) over $U$.

Let us equip $\mcE$ with a log structure pulled-back from $U^\mr{log}$ via $\pi : \mcE \migi U$;
we denote  the resulting log scheme  by   $\mcE^\mr{log}$.
The projection $\pi$ extends to  a morphism  $\mcE^\mr{log} \migi U^\mr{log}$, whose differential gives rise to the following short exact sequence:
\begin{equation}  \label{Ex0}
 0 \longmigi  \mfg_\mcE \longmigi  \widetilde{\mcT}_{\mcE^\mr{log}/T^\mr{log}} \stackrel{\mfa^\mr{log}_\mcE}{\longmigi}  \mcT_{U^\mr{log}/T^\mr{log}} \longmigi 0,
 \end{equation}
where 
$\widetilde{\mcT}_{\mcE^\mr{log}/T^\mr{log}}$ denotes the subsheaf $\pi_* (\mcT_{\mcE^\mr{log}/T^\mr{log}})^\mbG$ of $\mbG$-invariant sections of
$\pi_* (\mcT_{\mcE^\mr{log}/T^\mr{log}})$  (cf. ~\cite{Wak5}, \S\,1.2, (31)).
By 
a {\bf $T^\mr{log}$-connection}
 on  $\mcE$, we mean  an $\mcO_U$-linear morphism
$\nabla_\mcE : \mcT_{U^\mr{log}/T^\mr{log}} \migi \widetilde{\mcT}_{\mcE^\mr{log}/T^\mr{log}}$ such that
$\mfa^\mr{log}_\mcE \circ \nabla_\mcE = \mr{id}_{\mcT_{U^\mr{log}/T^\mr{log}}}$.
Also, by a  {\bf  log flat $\mbG$-torsor} 
 over  $U^\mr{log}/T^\mr{log}$, we mean  a pair $(\mcE, \nabla_\mcE)$ consisting of a right $\mbG$-torsor $\mcE$ over $U$ and a  $T^\mr{log}$-connection $\nabla_\mcE$ on $\mcE$.

\vspace{5mm}
\subsection{Monodromy of a logarithmic connection} \label{sc09}
\leavevmode\\ \vspace{-4mm}

Let  $\mfX := (X/S, \{ \sigma_i \}_{i=1}^r)$ be an $r$-pointed stable curve of genus $g$   over a $k$-scheme $S$.
Unless otherwise stated,  we  suppose, in this subsection,  that $r >0$.
Recall from  ~\cite{Wak5}, Definition 1.6.1, that,  to each  log flat $\mbG$-torsor $(\mcE, \nabla_\mcE)$ over $X^\mr{log}/S^\mr{log}$ ($= X^{\mfX \text{-} \mr{log}}/S^{\mfX \text{-} \mr{log}}$) and each $i \in \{ 1, \cdots, r \}$, 
one may associate an element 
\begin{align}
\mu_i^{(\mcE, \nabla_\mcE)} \in \Gamma (S, \sigma^*_i  (\mfg_{\mcE}))
\end{align}
 called the {\bf monodromy} of $(\mcE, \nabla_\mcE)$ at $\sigma_i$.

\vspace{3mm}
\bde  \label{D010}\leavevmode\\
 \ \ \ 
Let $\vec{\mu} := (\mu_i)_{i=1}^r$ be an element of $\prod_{i=1}^r \Gamma (S, \sigma^*_i  (\mfg_{\mcE}))$ 
 and 
$(\mcE, \nabla_\mcE)$  a log flat $\mbG$-torsor over $X^\mr{log}/S^\mr{log}$. 
Then, we shall say that  $(\mcE, \nabla_\mcE)$
 is {\bf of monodromies $\vec{\mu}$} if 
 $\mu_i^{(\mcE, \nabla_\mcE)} = \mu_i$ for any $i \in \{ 1, \cdots, r \}$.
If $r =0$, then we shall refer to any log flat $\mbG$-torsor 
as being {\bf of monodromy $\emptyset$}. 
\ede

\begin{rema}  \label{zz05g01} \leavevmode\\
 \vspace{-5mm}
 \begin{itemize}
 \item[(i)]
Let us consider the case where $\mbG = \mr{GL}_n$ (for some positive integer $n$).
Let $U^\mr{log}/T^\mr{log}$ be as in \S\,\ref{09}.
Recall  that giving a  $\mr{GL}_n$-torsor over  $U$ is
equivalent to giving a rank $n$ vector bundle on $U$.
This equivalence may be given by assigning 
$\mcE \mapsto (k^{\oplus n})_{\mcE}$.
Now, let $\mcE$ be a $\mr{GL}_n$-torsor over $U$ and $\mcV$ the rank $n$ vector bundle on $U$ corresponding to $\mcE$.
Then, there exists  a canonical  isomorphism 
 \begin{align} \label{eq1}
 \mcE nd_{\mcO_U} (\mcV)\isom (\mfg \mfl_n)_{\mcE}.
 \end{align} 
Moreover,  the notion of a $T^\mr{log}$-connection  on  $\mcE$  (defined above) 
 coincides, via this equivalence, with the notion of a $T^\mr{log}$-connection on  $\mcV$ (defined at the beginning of \S\,\ref{z041e}).
  We refer to ~\cite{Wak5}, \S\,4.2,  for a  detailed  discussion.

 \item[(ii)]
Suppose further that
$U^\mr{log}/T^\mr{log} = X^\mr{log}/S^\mr{log}$ for a pointed stable curve  $\mfX := (X/S, \{ \sigma_i \}_{i=1}^r)$.
Let $\nabla_\mcE$ be 
  an $S^\mr{log}$-connection 
    on $\mcE$, and 
denote by $\nabla_\mcV$ the $S^\mr{log}$-connection on $\mcV$ corresponding to $\nabla_\mcE$.
For each $i \in \{1, \cdots, r\}$,
we shall consider the composite
\begin{equation}
\mcV \xrightarrow{\nabla_\mcV} \Omega_{X^\mr{log}/S^\mr{log}}\otimes \mcV  \migi
\sigma_{i*}(\sigma_i^{*}(\Omega_{X^\mr{log}/S^\mr{log}})\otimes \sigma_i^{*}(\mcV))
\isom \sigma_{i*}(\sigma_i^{*}(\mcV)),
\end{equation}
where the second arrow arises from 
the adjunction relation ``$\sigma_i^{*}(-) \dashv \sigma_{i*}(-)$"
(i.e., ``the functor $\sigma^{*}_i(-)$ is left adjoint to the functor $\sigma_{i*}(-)$")
 and the third arrow arises   from  (\ref{f029}).
This composite corresponds (via the  relation ``$\sigma_i^{*}(-) \dashv \sigma_{i*}(-)$" again) to an  $\mcO_{ S}$-linear endomorphism  $\sigma_i^{*}(\mcV) \migi \sigma_i^{*}(\mcV)$.
Thus, we obtain
 an element  
  \begin{align} \label{eq2}
 \mu_i^{(\mcV, \nabla_\mcV)} \in  \Gamma (S, \mcE nd_{\mcO_{ S}} (\sigma_i^{*}(\mcV))) \ \left(= \Gamma (S, \sigma^*_i (\mcE nd_{\mcO_X} (\mcV))) \right),
  \end{align}
 which we shall refer to as the {\bf monodromy} of $(\mcV, \nabla_\mcV)$ (or just, of $\nabla_\mcV$) at $\sigma_i$.
This element $ \mu_i^{(\mcV, \nabla_\mcV)}$ coincides, via  (\ref{eq1}),  with the monodromy $\mu_i^{(\mcE, \nabla_\mcE)}$ of $(\mcE, \nabla_\mcE)$.
In the present paper, we shall not distinguish between  these notions of monodromy.
 \end{itemize}
 \end{rema}
\vspace{3mm}

\begin{rema}  \label{zz050} \leavevmode\\
 \ \ \  Let $\mcL^\flat : = (\mcL, \nabla_\mcL)$ be a log flat {\it line} bundle on $\mfX$.
Then, $\mcE nd_{\mcO_S} (\sigma_i^*(\mcL)) \cong \mcO_S$, and hence, $\mu_i^{\mcL^\flat}$ ($i = 1, \cdots, r$) may be thought of as an element of $\Gamma (S, \mcO_S)$.
In particular, it makes  sense to ask whether $\mu_i^{\mcL^\flat}$ lies in $k$ ($\subseteq \Gamma (S, \mcO_S)$) or not.
 \end{rema}

\vspace{5mm}
\subsection{Moduli of logarithmic  connections} \label{scfff091}
\leavevmode\\ \vspace{-4mm}

We shall write
\begin{align} \label{EQ2}
\overline{\mfC} \mfo_{g,r}
\end{align}
for the set-valued contravariant functor  on $\mfS \mfc \mfh_{/ \overline{\mfM}_{g,r}}$ which, to any
object $c_\mfX : S \migi \overline{\mfM}_{g,r}$ of $\mfS \mfc \mfh_{/ \overline{\mfM}_{g,r}}$ classifying a  pointed stable curve $\mfX$,  assigns the set of $S^\mr{log}$-connections on the line bundle $\Omega_{X^\mr{log}/S^\mr{log}}$. 
Also, we write
\begin{align}
\mfC \mfo_{g,r} := \overline{\mfC} \mfo_{g,r} \times_{\overline{\mfM}_{g,r}} \mfM_{g,r}.
\end{align}
It follows from a routine argument (or, an argument similar to the argument in the proof of Proposition \ref{P174} described later) that $\overline{\mfC} \mfo_{g,r}$ and $\mfC \mfo_{g,r}$ may be represented by  relative schemes of finite type over $\overline{\mfM}_{g,r}$ and $\mfM_{g,r}$ respectively.

Next,
 let $\vec{\mu} := (\mu_i)_{i=1}^r$ be an element of $k^{\times r}$ (= the product of $r$ copies of $k$), where $\vec{\mu} := \emptyset$ if $r =0$.
We shall write
\begin{align} \label{EQ3}
\overline{\mfC} \mfo_{g,r, \vec{\mu}}
\end{align}
for the closed substack of $\overline{\mfC} \mfo_{g,r}$ classifying connections 
 of monodromies $\vec{\mu}$.

\vspace{5mm}
\subsection{$p$-curvature} \label{sc091}
\leavevmode\\ \vspace{-4mm}

In the rest of this section, we discuss logarithmic connections in positive characteristic.
Suppose  that $\mr{char} (k)  = p >0$.
In the following, we shall recall (cf. ~\cite{Wak5}, \S\,3.2) the definition of the $p$-curvature of a connection.
Let $\mbG$ and $\mfX$  be as 
before
and let $(\mcE, \nabla_\mcE)$ be a log flat  $\mbG$-torsor over $X^\mr{log}/S^\mr{log}$.
Then, we obtain an $\mcO_X$-linear morphism $\mcT_{X^\mr{log}/S^\mr{log}}^{\otimes p} \migi \mfg_\mcE \ (\subseteq \widetilde{\mcT}_{\mcE^\mr{log}/S^\mr{log}})$ given by  $\partial^{\otimes p} \mapsto  (\nabla_\mcE (\partial))^{[p]} - \nabla_\mcE (\partial^{[p]})$, 
where $\partial$ is any local section of $\mcT_{X^\mr{log}/S^\mr{log}}$ and $\partial^{[p]}$ denotes the $p$-th symbolic power of $\partial$ (i.e., ``$\partial \mapsto \partial^{(p)}$" asserted in ~\cite{Og}, Proposition 1.2.1).
This morphism corresponds to an element
\begin{align}
\psi^{(\mcE, \nabla_\mcE)} \in  \Gamma (X, \Omega_{X^\mr{log}/S^\mr{log}}^{\otimes p}\otimes \mfg_\mcE),
\end{align}
which we shall refer to as the {\bf $p$-curvature} of $(\mcE, \nabla_\mcE)$  (cf. ~\cite{Wak5}, Definition 3.2.1).

If $\mbG = \mr{GL}_n$ and $(\mcV, \nabla_\mcV)$ denotes the log flat vector bundle 
 corresponds to $(\mcE, \nabla_\mcE)$, then $\psi^{(\mcE, \nabla_\mcE)}$ coincides, via (\ref{eq1}), with 
the classical definition of the $p$-curvature of $(\mcV, \nabla_\mcV)$  (cf., e.g.,  ~\cite{Wak6}, \S\,1.5). 

\vspace{5mm}
\subsection{Canonical connections arising from the Frobenius morphism} \label{sc091}
\leavevmode\\ \vspace{-4mm}

We shall recall the canonical connection arising from pull-back via Frobenius morphisms.
Let $Y$ be an $S$-scheme with structure morphism $f : Y \migi S$.
 Denote by $F_S : S \migi S$ (resp., $F_Y : Y \migi Y$) the absolute (i.e., $p$-th power) Frobenius endomorphism of $S$ (resp., $Y$).
The  {\bf Frobenius twist of $Y$ over $S$} is, by definition,  the base-change $Y^{(1)}_S$ ($:= Y \times_{S, F_S} S$) of $Y$ via
 $F_S : S \migi S$.
Denote by $f^{(1)} : Y^{(1)}_S \migi S$ the structure morphism of 
$Y^{(1)}_S$.
The {\bf relative Frobenius morphism of $Y$ over $S$} is  the unique morphism $F_{Y/S} : Y \migi Y^{(1)}_S$ over $S$ that makes the following diagram commute:
\begin{align}
\vcenter{\xymatrix{
Y \ar@/^10pt/[rrrrd]^{F_Y}\ar@/_10pt/[ddrr]_{f} \ar[rrd]_{ \ \ \ \ \ \ \  \  \   F_{Y/S}} & & &   &   \\
& & Y^{(1)}_S  \ar[rr]_{\mr{id}_Y \times F_S} \ar[d]^{f^{(1)}}  \ar@{}[rrd]|{\Box}  &  &  Y \ar[d]_{f} \\
&  & S \ar[rr]_{F_S} & &  S.
}}
\end{align}

 Now, let $\mfX := (f :X \migi S, \{ \sigma_i \}_{i=1}^r)$ be as before and $\mcU$   a vector bundle on $X^{(1)}_S$.
  Then,  one may construct (cf. \cite{Wak5}, \S\,3.3)  an $S^\mr{log}$-connection
\begin{equation} \label{eeqq1}
\nabla^\mr{can}_\mcU : F^*_{X/S}(\mcU) \migi \Omega_{X^\mr{log}/S^\mr{log}}\otimes F^*_{X/S}(\mcU) 
\end{equation}
on the pull-back $F^*_{X/S}(\mcU)$ of $\mcU$  which is uniquely determined by the condition that the sections of the subsheaf $F^{-1}_{X/S}(\mcU)$ ($\subseteq F^*_{X/S}(\mcU)$) are contained in $\mr{Ker}(\nabla^\mr{can}_\mcU)$.
We shall refer to $\nabla^\mr{can}_\mcU$ as the {\bf canonical $S^\mr{log}$-connection} on $F^*_{X/S}(\mcU)$.
One verifies immediately  that
\begin{equation} \label{k09}
\mr{Im} (\nabla^\mr{can}_\mcU) \subseteq \Omega_{X/S}\otimes F^*_{X/S}(\mcU) \ \left(\subseteq \Omega_{X^\mr{log}/S^\mr{log}}\otimes F^*_{X/S}(\mcU) \right)
\end{equation}
(i.e., $\nabla_\mcU^\mr{can}$ comes  from a {\it non-logarithmic} connection on $F^*_{X/S}(\mcU)$).
Moreover, we have
\begin{align} \label{eeqq2}
\psi^{(F^*_{X/S}(\mcU), \nabla^\mr{can}_\mcU)} = 0 
\end{align}
and (under the assumption that $r >0$)
\begin{align} \label{eeqq3}
\mu_i^{(F^*_{X/S}(\mcU), \nabla^\mr{can}_\mcU)} = 0
\end{align}
for any $i \in \{ 1, \cdots, r \}$.

\vspace{5mm}
\subsection{Moduli of connections on a line bundle} \label{sfc091}
\leavevmode\\ \vspace{-4mm}

We shall write
\begin{align} \label{eqeq44}
\widetilde{\mbF}_p := \{ 0, 1, \cdots, p-1 \} \ \left(\subseteq \mbZ \right), 
\end{align}
and  write $\tau$ for 
 the natural composite bijection
\begin{align} \label{eqeq45}
\tau  : \widetilde{\mbF}_p \migiincl \mbZ \migisurj \mbF_p \ (:= \mbZ /p \mbZ).
\end{align}

Let  $\vec{\mu} := (\mu_i)_{i=1}^r$ be an element of $k^{\times r}$  (where  $\vec{\mu} := \emptyset$ if $r =0$), $d$ an integer,  and $\mcL$ 
a line bundle on the tautological curve  $\overline{\mfC}_{g,r}$ (resp., $\mfC_{g,r}$) whose restriction to any   fiber of $\overline{f}_{\mft \mfa \mfu} : \overline{\mfC}_{g,r} \migi \overline{\mfM}_{g,r}$  (resp., $f_{\mft \mfa \mfu} : \mfC_{g,r} \migi \mfM_{g,r}$) has  degree $d$.   
We shall denote by 
\begin{align}  \label{eqeq33}
\overline{\mfC} \mfo^{\psi =0}_{\mcL, g,r, \vec{\mu}} \ \ (\text{resp.,} \ \mfC \mfo^{\psi =0}_{\mcL, g,r, \vec{\mu}})
\end{align}
the set-valued contravariant functor on $\mfS \mfc \mfh_{/\overline{\mfM}_{g,r}}$ (resp., $\mfS \mfc \mfh_{/\mfM_{g,r}}$)
which, to any object $S \migi \overline{\mfM}_{g,r}$  (resp., $S \migi \mfM_{g,r}$) (where we denote by $\mfX := (X/S, \{ \sigma_i \}_{i=1}^r)$ the pointed curve classified by this object), assigns the set of $S^{\mfX \text{-} \mr{log}}$-connections on $\mcL |_X$ with vanishing $p$-curvature. 
If $\mcL = \Omega_{\overline{\mfC}_{g,r}^\mr{log}/\overline{\mfM}_{g,r}^\mr{log}}$, then we shall write
\begin{align}
\overline{\mfC} \mfo^{\psi =0}_{g,r, \vec{\mu}} := \overline{\mfC} \mfo^{\psi =0}_{\mcL, g,r, \vec{\mu}} \ \ (\text{resp.,} \ \mfC \mfo^{\psi =0}_{g,r, \vec{\mu}} := \mfC \mfo^{\psi =0}_{\mcL, g,r, \vec{\mu}})
\end{align}
for simplicity.

\vspace{3mm}
\bpr\label{prfg001tt}
\leavevmode\\
 \ \ \ 
Let $\mcL$ be as above and $\vec{\mu} :=  (\mu_i)_{i=1}^r\in \mbF_p^{\times r}$ (where $\vec{\mu} := \emptyset$ if $r =0$).
\begin{itemize}
\item[(i)]
$\mfC \mfo^{\psi =0}_{\mcL, g,r, \vec{\mu}}$ is nonempty if and only if 
$p | (d + \sum_{i=1}^r \tau^{-1} (\mu_i))$.
\item[(ii)]
Suppose that $\mfC \mfo^{\psi =0}_{\mcL, g,r, \vec{\mu}}$ is nonempty.
Then, $\mfC \mfo^{\psi =0}_{\mcL, g,r, \vec{\mu}}$ may be represented by a Deligne-Mumford stack over $k$ which is finite and faithfully flat over $\mfM_{g,r}$ of degree $p^g$.
Moreover, the open substack $\mfC \mfo^{\psi =0}_{\mcL, g,r, \vec{\mu}} \times_{\mfM_{g,r}} \mfM_{g,r}^\mr{ord}$ is \'{e}tale over $\mfM_{g,r}^\mr{ord}$.
\end{itemize}
 \epr
\begin{proof}
Let $S$ be a $k$-scheme and $\mfX := (X/S, \{ \sigma_i \}_{i=1}^r)$ an $r$-pointed proper smooth curve over $S$ of genus $g$. 
Denote by $c_\mfX : S \migi \mfM_{g,r}$ the classifying morphism of $\mfX$.
For each $d' \in \mbQ$,  denote by $\mr{Pic}_{X/S}^{d'}$ (resp., $\mr{Pic}_{X^{(1)}_S/S}^{d'}$)  the relative Picard scheme of $X/S$ (resp., $X^{(1)}_S/S$) classifying the set of (equivalence classes, relative to the equivalence relation determined by tensoring with a line bundle pulled-back from the base $S$, of) degree $d'$ line bundles on $X$ (resp., on $X^{(1)}_S$).
(Here, we take $
\mr{Pic}_{X^{(1)}_S/S}^{d'} := \emptyset$ if $d' \in \mbQ \setminus \mbZ$.)
Denote by $c_{\mcL, \vec{\mu}} : S \migi \mr{Pic}^{d+\sum_{i=1}^r \tau^{-1} (\mu_i)}_{X/S}$
the classifying morphism of the line bundle $\mcL ( \sum_{i=1}^r \tau^{-1} (\mu_i) \cdot \sigma_{\mft \mfa \mfu, i}|_{\mfC_{g,r}})$ restricted to $X$ via $c_\mfX$, i.e., $\mcL |_X (\sum_{i=1}^r \tau^{-1}(\mu_i) \cdot \sigma_i)$.
Let us consider the morphism 
\begin{align}
\mcV er : \mr{Pic}^{\frac{1}{p} \cdot (d+\sum_{i=1}^r \tau^{-1} (\mu_i))}_{X^{(1)}_S/S} & \migi \mr{Pic}^{d+\sum_{i=1}^r \tau^{-1} (\mu_i)}_{X/S} \notag \\
\vin \hspace{28mm}& \hspace{10mm} \vin
 \\
[\mcN] \hspace{26mm}  & \mapsto \hspace{0mm} [F^*_{X/S}(\mcN)] \notag
\end{align}
determined by pull-back via
$F_{X/S}$.
In what follows, we shall prove the claim that 
$\mfC \mfo^{\psi =0}_{\mcL, g,r, \vec{\mu}} \times_{\mfM_{g,r}, c_\mfX} S$
is  isomorphic to the inverse image  $\mcV er^{-1} (c_{\mcL, \vec{\mu}})$ of $c_{\mcL, \vec{\mu}}$ via the morphism $\mcV er$.

First, let $\mcN$ be a line bundle on $X^{(1)}_S$ classified by 
$\mcV er^{-1} (c_{\mcL, \vec{\mu}})$, which admits, by definition, 
  an isomorphism  
\begin{align} \label{eqeq47}
F^*_{X/S} (\mcN) \isom  \mcL |_X ( \sum_{i=1}^r \tau^{-1} (\mu_i) \cdot \sigma_i).
\end{align}
Consider 
 the $S$-connection on $\mcL |_X ( \sum_{i=1}^r \tau^{-1} (\mu_i) \cdot \sigma_i)$
corresponding, via (\ref{eqeq47}),  to the canonical (logarithmic) $S$-connection $\nabla^\mr{can}_{\mcN}$;
it restricts to an 
 $S$-connection $\nabla_{\mcL |_X}$ on  $\mcL |_X$,
which    has  vanishing $p$-curvature and  monodromies $\vec{\mu}$.
The assignment $\mcN \mapsto \nabla_{\mcL |_X}$ is functorial with respect to $S$,
and hence, determines a morphism 
\begin{align} \label{EE1ff0}
\mcV er^{-1} (c_{\mcL, \vec{\mu}})  \migi \mfC \mfo^{\psi =0}_{\mcL, g,r, \vec{\mu}} \times_{\mfM_{g,r}, c_\mfX} S 
\end{align}
over $S$.

On the other hand, let us take an $S$-connection $\nabla_{\mcL |_X}$ on 
$\mcL |_X$ with vanishing $p$-curvature and of monodromies $\vec{\mu}$.
One may find a unique $S$-connection $\nabla_{\mcL |_X, + \vec{\mu}}$ on  $\mcL |_X (\sum_{i=1}^r \tau^{-1} (\mu_i) \cdot \sigma_i)$ whose restriction to $\mcL |_X$ coincides with $\nabla_{\mcL |_X}$.
Let us  regard  $F_{X/S*}(\mr{Ker} (\nabla_{\mcL|_X, + \vec{\mu}}))$ as an $\mcO_{X_S^{(1)}}$-module.
Then, 
 the natural  inclusion $F_{X/S*}(\mr{Ker} (\nabla_{\mcL|_X, + \vec{\mu}})) \migiincl F_{X/S*} (\mcL  |_X  (\sum_{i=1}^r \tau^{-1} (\mu_i) \cdot \sigma_i))$ is  $\mcO_{X_S^{(1)}}$-linear and  corresponds, via the adjunction relation ``$F_{X/S}^{*}(-) \dashv F_{*}(-)$",
 to a morphism 
 \begin{align} \label{eqeq49}
 F^*_{X/S} (F_{X/S*} (\mr{Ker} (\nabla_{\mcL |_X, + \vec{\mu}}))) \migi \mcL |_X.  (\sum_{i=1}^r \tau^{-1} (\mu_i) \cdot \sigma_i).
 \end{align} 
Since $\nabla_{\mcL |_X, + \vec{\mu}}$ has  monodromies $(0, 0, \cdots, 0) \in \mbF_p^{\times r}$,
 it turns out to be an isomorphism.
By this isomorphism,  $F_{X/S*} (\mr{Ker} (\nabla_{\mcL |_X, + \vec{\mu}}))$ specifies 
a morphism $S \migi  \mcV er (c_{\mcL, \vec{\mu}})$ as its classifying morphism.
The resulting assignment $\nabla_{\mcL |_X} \mapsto F_{X/S*} (\mr{Ker} (\nabla_{\mcL |_X, +\vec{\mu}}))$ determines  a morphism
\begin{align} \label{EE11}
\mfC \mfo^{\psi =0}_{\mcL, g,r,  \vec{\mu}} \times_{\mfM_{g,r}, c_\mfX} S  \migi \mcV er^{-1} (c_{\mcL, \vec{\mu}}) 
\end{align}
over $S$.
One verifies that  (\ref{EE1ff0}) and  (\ref{EE11}) are inverses to  each other,  and hence, 
$\mfC \mfo^{\psi =0}_{\mcL, g,r,  \vec{\mu}} \times_{\mfM_{g,r}, c_\mfX} S$ is isomorphic to $\mcV er^{-1} (c_{\mcL, \vec{\mu}})$.
This completes the proof of the claim.
Then,  assertions (i) and (ii) follow from the claim just proved and the well-known fact that  the morphism $\mcV er$ 
 is finite and faithfully flat of degree $p^p$ and moreover \'{e}tale if $c_{\mfX}$ lies in  $\mfM_{g,r}^\mr{ord}$. 
\end{proof}

\vspace{10mm}
\section{Cartan connections on a pointed stable curve} \vspace{3mm}

In this section, we shall discuss   logarithmic connections  on a certain  torsor (i.e., $\mcE^\dagger_{\mbT, U^\mr{log}/T^\mr{log}}$ defined in (\ref{f016})) called  {\it Cartan connections} (cf. Definition \ref{z012}).
As 
shown in \S\,\ref{Sec4}, such  connections correspond bijectively to generic Miura opers;
this result  may be thought of as  a generalization of  ~\cite{Fr}, Proposition 8.2.2.

Denote by $\mbG_m$ the multiplicative group over $k$.
In the following, for a scheme (or more generally, a stack) $Y$ and an $\mcO_Y$-module $\mcV$,  we shall write
$\mbV (\mcV)$ for  the relative affine space over $Y$ associated with $\mcV$.

\vspace{5mm}
\subsection{Algebraic groups and Lie algebras} \label{SSQ1}
\leavevmode\\ \vspace{-4mm}

Let   $\mbG$ be  a split connected semisimple algebraic group  of adjoint type over $k$. 
Assume   that
either one of the following  three  conditions (Char)$_{0}$,  $(\mr{Char})_{p}$,  and $(\text{Char})_{p}^{\mfs \mfl}$ is  satisfied:
\begin{itemize} 
\item[]
\begin{itemize}
\item[ $(\text{Char})_{0}$] : $\mr{char}(k)=0$;
\vspace{1.5mm}
\item[ $(\text{Char})_{p}$]  : $\mr{char}(k) =p> 2 \cdot h$, where $h$ denotes the Coxeter number of $\mbG$;
\vspace{1.5mm}
\item[ $(\text{Char})_{p}^{\mfs \mfl}$]  
 : $\mr{char}(k) =p>0$ and $\mbG = \mr{PGL}_n$
   for a positive integer $n$ with $n <p$.
\end{itemize}
\end{itemize}
In particular, 
$(\text{Char})_{p}$ and $(\text{Char})_{p}^{\mfs \mfl}$  imply respectively 
that  $p$ does not divide the order of  the Weyl group of $\mbG$  (and  is very good for $\mbG$).

Let us fix  a split maximal torus $\mbT_\mbG$
  of $\mbG$, 
  a Borel subgroup $\mbB_\mbG$
     of $\mbG$ containing $\mbT_\mbG$.
If there is no fear of causing confusion, then we shall write $\mbT := \mbT_\mbG$, $\mbB := \mbB_\mbG$ 
for simplicity.
Denote by 
$\Phi^+$  the set of positive  roots in $\mbB$ with respect to $\mbT$ and by   $\Phi^-$ the set of negative roots. 
Also, denote by $\Gamma$ ($\subseteq \Phi^+$) the set of  simple positive  roots.
Denote by $\mfg$, $\mfb$,  and $\mft$   the  Lie algebras of $\mbG$, $\mbB$,   and $\mbT$ respectively (hence $\mft \subseteq \mfb \subseteq \mfg$).
For each $\alpha \in \Phi^+ \cup \Phi^-$, we write
\begin{equation} \label{22223}
\mfg^\alpha := \big\{ x \in \mfg \ \big| \ \text{$\mr{ad}(t)(x) =  \alpha  (t) \cdot x$ for all  $t\in \mbT$}  \big\}.\end{equation}
$\mfg$ admits a canonical 
decomposition 
\begin{align} \label{22224}
\mfg =  \mft \oplus \left(\bigoplus_{\alpha \in \Phi^+} \mfg^\alpha \right) \oplus \left(\bigoplus_{\beta \in \Phi^-} \mfg^\beta \right)
\end{align}
(which restricts  to a decomposition $\mfb =   \mft \oplus (\bigoplus_{\alpha \in \Phi^+} \mfg^\alpha)$ of $\mfb$).
By means of this decomposition, we obtain  a unique decreasing filtration $\{ \mfg^j \}_{j \in \mbZ}$ on $\mfg$ such that 
\begin{itemize}
\item[$\bullet$]
$\mfg^0 = \mfb$, $\mfg^0 /\mfg^1 = \mft$, $\mfg^{-1}/\mfg^0 = \bigoplus_{\alpha \in \Gamma} \mfg^{-\alpha}$;
\item[$\bullet$]
$[\mfg^{j_1}, \mfg^{j_2}] \subseteq \mfg^{j_1 +j_2}$ for $j_1$, $j_2 \in \mbZ$.
\end{itemize}
We shall write
$\check{\rho}$ for the one-parameter subgroup $\mbG_m \migi \mbT$ of $\mbT$ defined to be  the sum $\sum_{\alpha \in \Gamma}  \check{\omega}_\alpha$, where each $\check{\omega}_\alpha$ ($\alpha \in \Gamma$) denotes  the fundamental coweight  of $\alpha$.
 By passing to differentiation, we consider $\check{\rho}$ as an element of $\mft$.

Denote by $\mbW$ the Weyl group of $(\mbG, \mbT)$, i.e., $\mbW := N_\mbG (\mbT)/\mbT$, where $N_\mbG (\mbT)$ denotes the normalizer of $\mbT$ in $\mbG$.
We shall identify $\mfg$ and  $\mft$ with $\mr{Spec}(\mbS_k (\mfg^\vee))$ and $\mr{Spec}(\mbS_k (\mft^\vee))$ respectively, where  for a $k$-vector space $\mfa$ we denote by  $\mbS_k (\mfa)$ the symmetric algebra on $\mfa$ over $k$.
Consider the GIT quotient $\mfg \ooalign{$/$ \cr $\,/$}\hspace{-0.5mm}\mbG$  (resp., $\mft \ooalign{$/$ \cr $\,/$}\hspace{-0.5mm}\mbW$) of 
   $\mfg$ (resp., $\mft$) by 
the adjoint action of $\mbG$ (resp.,  $\mbW$), i.e., the spectrum of the ring of polynomial invariants $\mbS_k (\mfg^\vee)^\mbG$ (resp., $\mbS_k (\mft^\vee)^\mbW$) on  $\mfg$ (resp., $\mft$).
Let us write 
 \begin{equation} \mfc_\mbG :=  \mft \ooalign{$/$ \cr $\,/$}\hspace{-0.5mm}\mbW.
 \label{c}
 \end{equation}
If there is no fear of causing confusion, then we shall write
$\mfc := \mfc_\mbG$ for simplicity.
 A Chevalley's theorem asserts (cf. ~\cite{Ngo}, Theorem 1.1.1; ~\cite{KW}, Chap.\,VI, Theorem 8.2) that
the natural morphism $\mbS_k (\mfg^\vee)^\mbG \migi \mbS_k (\mft^\vee)^\mbW$ is
an isomorphism.
Thus, one may define a morphism
\begin{equation} \chi :  \mfg \migi \mfc_\mbG
\label{chi}
\end{equation}
of $k$-schemes
to be the composite of the natural quotient $\mfg \migisurj \mfg \ooalign{$/$ \cr $\,/$}\hspace{-0.5mm}\mbG$ and the inverse of the resulting  isomorphism $\mfc_\mbG\isom \mfg \ooalign{$/$ \cr $\,/$}\hspace{-0.5mm}\mbG$.
  $\chi$  factors through the quotient  $\mfg \migi [\mfg /\mbG]$. 
 (Here,  $[\mfg /\mbG]$ is the quotient stack representing the contravariant  functor which, to any $k$-scheme $T$, assigns the groupoid of pairs $(\mcF, R)$ consisting of a right $\mbG$-torsor $\mcF$ over $T$ and $R \in \Gamma (T, \mfg_\mcF)$.)
We shall equip $\mfc_\mbG$ with  the  $\mbG_m$-action  that  comes from the homotheties on $\mfg$ (i.e., the natural grading on $\mbS_k (\mfg^\vee)$).

\vspace{5mm}
\subsection{Cartan  connections}
\leavevmode\\ \vspace{-4mm}

Let $T^\mr{log}$ be an fs log scheme over $k$ and $f^\mr{log}  : U^\mr{log} \migi T^\mr{log}$ a log-curve over $T^\mr{log}$.
We shall denote by  $\mcE^\dagger_{\mbT, U^\mr{log}/T^\mr{log}}$  the $\mbT$-torsor over $U$ defined to be
\begin{align} \label{f016}
\mcE^\dagger_{\mbT, U^\mr{log}/T^\mr{log}} := (\Omega_{U^\mr{log}/T^\mr{log}})^\times \times^{\mbG_m, \check{\rho}} \mbT,
\end{align}
where $(\Omega_{U^\mr{log}/T^\mr{log}})^\times$ denotes the $\mbG_m$-torsor over $U$ corresponding to the line bundle $\Omega_{U^\mr{log}/T^\mr{log}}$.
If there is no fear of causing confusion, then we shall write $\mcE^\dagger_\mbT := \mcE^\dagger_{\mbT, U^\mr{log}/T^\mr{log}}$ for simplicity.

\vspace{3mm}
\bde  \label{z012}\leavevmode\\
 \ \ \ A {\bf $\mfg$-Cartan connection} on $U^\mr{log}/T^\mr{log}$ is a $T^\mr{log}$-connection on $\mcE^\dagger_{\mbT}$. 
 If $U^\mr{log}/T^\mr{log} = X^{\mfX \text{-} \mr{log}}/S^{\mfX \text{-} \mr{log}}$
 for some pointed stable curve $\mfX := (X/S, \{ \sigma_i \}_{i=1}^r)$, then we shall refer to any $\mfg$-Cartan connection on $X^{\mfX \text{-} \mr{log}}/S^{\mfX \text{-} \mr{log}}$ as a {\bf $\mfg$-Cartan connection on $\mfX$}.
\ede
\vspace{3mm}

Next, let $\mfX := (f : X \migi S, \{ \sigma_i \}_{i=1}^r)$ be an $r$-pointed stable curve of genus $g$  over a $k$-scheme $S$, and 
suppose that $r >0$.
 Let us fix $i \in \{1, \cdots, r \}$.
By the definition of $\mcE^\dagger_\mbT$ ($:= \mcE^\dagger_{\mbT, X^\mr{log}/S^\mr{log}}$),
we have  a sequence of  isomorphisms
\begin{align} \label{eqeq39}
\sigma_i^* (\mcE^\dagger_\mbT) \isom (\sigma^*_i (\Omega^\times_{X^\mr{log}/S^\mr{log}})) \times^{\mbG_m, \check{\rho}} \mbT\isom (S \times_k \mbG_m) \times^{\mbG_m, \check{\rho}} \mbT\isom  S \times_k \mbT,
\end{align}
where the second arrow 
arises from 
 (\ref{f029}).
It follows that  the $\mbT$-torsor $\sigma_i^* (\mcE^\dagger_\mbT)$ is trivial
  and 
 we have a sequence  of isomorphisms:
\begin{align} \label{f004}
\Gamma (S, \sigma^*_i (\mft_{{\mcE}^{\dagger}_\mbT})) \isom \Gamma (S, \mft_{\sigma_i^*(\mcE^\dagger_\mbT)}) \isom \Gamma (S, \mft_{S \times_k \mbT}) \isom \mft (S).
\end{align}
The monodromy of any $S^\mr{log}$-connection on $\mcE^\dagger_\mbT$ (at each marked point $\sigma_i$) may be thought, via (\ref{f004}),  of as an element of   $\mft (S)$.
Hence, for each $\vec{\mu} \in \mft (S)^{\times r}$,  it makes sense to speak of a $\mfg$-Cartan connection on $\mfX$ {\it of monodromies $\vec{\mu}$}.

Let $\vec{\mu} 
 \in \mft (k)^{\times r}$ (where if $r =0$, then we take $\vec{\mu} := \emptyset$).
We shall write
\begin{align} \label{f011}
\mfC  \overline{\mfC} \mfo_{\mfg,  g,r}
  \ \ \ (\text{resp.,} \ \mfC \overline{\mfC} \mfo_{\mfg, g,r, \vec{\mu}} )
\end{align}
for the set-valued contravariant functor on $\mfS \mfc \mfh_{/ \overline{\mfM}_{g,r}}$ which, to each 
object $S \migi \overline{\mfM}_{g,r}$ of $\mfS \mfc \mfh_{/ \overline{\mfM}_{g,r}}$ classifying  a pointed stable curve $\mfX$,  assigns the set of $\mfg$-Cartan connections on $\mfX$ 
 (resp., the set of $\mfg$-Cartan connections on $\mfX$ of monodromies $\vec{\mu}$).
Also, we shall write
\begin{align}
\mfC  \mfC \mfo_{\mfg,  g,r} := \mfC  \overline{\mfC} \mfo_{\mfg,  g,r} \times_{\overline{\mfM}_{g,r}} \mfM_{g,r}, \ \ \  \mfC  \mfC \mfo_{\mfg, g,r, \vec{\mu}} := \mfC  \overline{\mfC} \mfo_{\mfg, g,r, \vec{\mu}} \times_{\overline{\mfM}_{g,r}} \mfM_{g,r}.
\end{align}
Then, the following proposition holds.
(Notice that we will apply, in advance, the result of 
Proposition \ref{prfg001} described later
 in order to prove assertion (iii).)

\vspace{3mm}
\bpr \label{P174}\leavevmode\\
\vspace{-5mm}
\begin{itemize}
\item[(i)]
Both  $\mfC \overline{\mfC} \mfo_{\mfg, g,r}$ and $\mfC \overline{\mfC} \mfo_{\mfg, g,r,  \vec{\mu}}$ may be represented by (possibly empty) relative affine schemes over $\overline{\mfM}_{g,r}$   of finite type.
In particular,  these moduli functors may be represented by Deligne-Mumford stacks over $k$.
\item[(ii)]
Assume  that $r >0$.
Then,  $\mfC \overline{\mfC} \mfo_{\mfg, g,r}$
  may be represented by a relative affine space over $\overline{\mfM}_{g,r}$ modeled on $\mbV (\overline{f}_{\mft \mfa \mfu *} (\Omega_{\overline{\mfC}_{g,r}^\mr{log}/\overline{\mfM}^\mr{log}_{g,r}}) \otimes_k \mft)$.
In particular, the projection $\mfC \overline{\mfC} \mfo_{\mfg, g,r} \migi \overline{\mfM}_{g,r}$ 
 is smooth of relative dimension $(g-1 +r) \cdot \mr{rk} (\mfg)$, where $\mr{rk} (\mfg)$ denotes the rank of $\mfg$.
\item[(iii)]
Assume  that the $\mr{char} (k) = p >0$,  $r =0$, and  $p | (2g-2)$.
Then, $\mfC \mfC \mfo_{\mfg, g,0}$ may be represented by a relative affine space over $\mfM_{g,0}$ modeled on $\mbV (f_{\mft \mfa \mfu *} (\Omega_{\mfC_{g, 0}^\mr{log}/\mfM^\mr{log}_{g,0}}) \otimes_k \mft)$.
In particular, the projection $\mfC \mfC \mfo_{\mfg, g,0} \migi \mfM_{g,0}$ is smooth of relative dimension 
$g \cdot \mr{rk} (\mfg)$.
\end{itemize}
 \epr
\begin{proof}
To begin with, let us fix an object $c_\mfX : S \migi \overline{\mfM}_{g,r}$ of $\mfS \mfc \mfh_{/\overline{\mfM}_{g,r}}$, which classifies a pointed stable curve
 $\mfX := (f: X \migi S, \{ \sigma_i \}_{i=1}^r)$.
We shall write 
\begin{align}
\mfC \overline{\mfC} \mfo_{\mfg, \mfX} := \mfC \overline{\mfC} \mfo_{\mfg, g,r}\times_{\overline{\mfM}_{g,r}, c_\mfX} S.
\end{align}

First, we shall  prove 
  assertion (i).
Since  $\mfC \overline{\mfC} \mfo_{\mfg, g,r,  \vec{\mu}}$ is a closed substack of $\mfC \overline{\mfC} \mfo_{\mfg, g,r}$,  it suffices  to 
prove that $\mfC \overline{\mfC} \mfo_{\mfg, \mfX}$ may be represented by a  (possibly empty) relative affine scheme over $S$ of finite type.
Consider the $\mcO_X$-linear  morphism
\begin{align} \label{eqeq222}
\mcH om_{\mcO_X} (\mcT_{X^\mr{log}/S^\mr{log}}, \widetilde{\mcT}_{\mcE^{\dagger\mr{log}}_\mbT/S^\mr{log}}) \migi  \mcE nd_{\mcO_X} (\mcT_{X^\mr{log}/S^\mr{log}}) \ \left(\cong \mcO_X \right)
\end{align}
obtained  by composition  with $\mfa^\mr{log}_{\mcE^\dagger_\mbT} : \widetilde{\mcT}_{\mcE^{\dagger\mr{log}}_\mbT/S^\mr{log}} \migi \mcT_{X^\mr{log}/S^\mr{log}}$.
This morphism gives rise to   a morphism 
\begin{align} \label{eq060}
\mbV (f_*(\mcH om_{\mcO_X} (\mcT_{X^\mr{log}/S^\mr{log}}, \widetilde{\mcT}_{\mcE^{\dagger\mr{log}}_\mbT/S^\mr{log}}))) \migi \mbV (f_* (\mcO_X)) \ \left(= \mbV (\mcO_S) \right)
\end{align}
of $S$-schemes.
By the definition of an $S^\mr{log}$-connection,
 $\mfC \overline{\mfC} \mfo_{\mfg,  \mfX}$ 
 is isomorphic to   the  inverse image, via (\ref{eq060}), of the section $S \migi \mbV (\mcO_S)$ corresponding to $1 \in \Gamma (S, \mcO_S)$.
In particular,  
  $\mfC \overline{\mfC} \mfo_{\mfg,  \mfX}$ may be represented by a closed subscheme of 
  $\mbV (f_*(\mcH om_{\mcO_X} (\mcT_{X^\mr{log}/S^\mr{log}}, \widetilde{\mcT}_{\mcE^{\dagger\mr{log}}_\mbT/S^\mr{log}})))$ and 
 admits a free transitive action of 
\begin{align}
\mbV (f_{*} (\mcH om_{\mcO_X} (\mcT_{X^\mr{log}/S^\mr{log}}, \mr{Ker} (\mfa^\mr{log}_{\mcE^\dagger_\mbT}))))
& \cong \mbV ( f_{*} (\mcH om_{\mcO_X} (\mcT_{X^\mr{log}/S^\mr{log}}, \mft_{\mcE^\dagger_\mbT}))) \\
& \cong \mbV ( f_* (\Omega_{X^\mr{log}/S^\mr{log}}) \otimes_k \mft). \notag
\end{align}
This completes the proof 
of assertion (i).

Next, we shall  prove  assertion (ii).
By the above  discussion, it suffices to prove that $\mfC \overline{\mfC} \mfo_{\mfg, \mfX}$  admits,  locally on $S$, a section $S \migi \mfC \overline{\mfC} \mfo_{\mfg, \mfX}$.
One may assume, without loss of generality, that $S = \mr{Spec} (R)$ for some $k$-algebra $R$.
 The  inverse image of the section $1 \in \Gamma (X, \mcO_X)$ via (\ref{eqeq222}) has, locally on $X$, a section.
That is to say,
there exists a collection $\{ (U_i, \nabla_i) \}_{i \in I}$ indexed by a set $I$, where  $\{ U_i \}_{i \in I}$ is  an open covering  of $X$ and each $\nabla_i$ 
 denotes an element of
$\Gamma (U_i, \mcH om_{\mcO_X} (\mcT_{X^\mr{log}/S^\mr{log}}, \widetilde{\mcT}_{\mcE^{\dagger\mr{log}}_\mbT/S^\mr{log}}))$ whose image via (\ref{eqeq222}) coincides with $1 \in \Gamma (U_i, \mcO_X)$.
Hence, $\nabla_{I \times I} := \{ \nabla_i |_{U_i \cap U_j} - \nabla_j |_{U_i \cap U_j} \}_{(i, j) \in I \times I}$ specifies a \v{C}ech $1$-cocycle  of $\{ U_i \}_{i \in I}$ with coefficients in $\mcH om_{\mcO_X} (\mcT_{X^\mr{log}/S^\mr{log}}, \mr{Ker} (\mfa^\mr{log}_{\mcE^\dagger_\mbT})) \ (\cong \Omega_{X^\mr{log}/S^\mr{log}} \otimes_k \mft )$.
Denote by $\overline{\nabla}_{I \times I}$
the element  of $H^1 (X, \Omega_{X^\mr{log}/S^\mr{log}} \otimes_k \mft)$ represented by  $\nabla_{I \times I}$.
By Serre duality, we obtain the following  sequence of isomorphisms of $R$-modules
\begin{align}
H^1 (X, \Omega_{X^\mr{log}/S^\mr{log}} \otimes_k \mft) & \isom H^1 (X, \Omega_{X^\mr{log}/S^\mr{log}})^{\oplus \mr{rk} (\mfg)} \\
& \isom    (H^0(X, \Omega_{X^\mr{log}/S^\mr{log}}^\vee \otimes \omega_{X/S})^\vee)^{\oplus \mr{rk} (\mfg)} \notag \\
& \isom   (H^0(X, \mcO_X (- D_{\mfX}))^\vee)^{\oplus \mr{rk} (\mfg)} \notag \\
& = 0, \notag
\end{align}
where $D_{\mfX}$ denotes the 
effective relative  divisor on $X$ defined to be the union of
  the marked points  $\sigma_i$ ($i = 1, \cdots, r$) and $\omega_{X/S}$ denotes the dualizing sheaf of $X$ over $S$, which is isomorphic to $\Omega_{X^\mr{log}/S^\mr{log}} (-D_{\mfX})$.
Hence, we have $\overline{\nabla}_{I \times I} = 0$.
This implies that 
 after possibly replacing each $\nabla_i$ by
$\nabla_i + \delta_i$ (for some $\delta_i \in \Gamma (U_i, \mcH om_{\mcO_X} (\mcT_{X^\mr{log}/S^\mr{log}}, \mr{Ker} (\mfa^\mr{log}_{\mcE^\dagger_\mbT})))$),
the sections $\nabla_i$ may be glued together to obtain  a global section $\nabla$ of  
$\mcH om_{\mcO_X} (\mcT_{X^\mr{log}/S^\mr{log}}, \widetilde{\mcT}_{\mcE^{\dagger\mr{log}}_\mbT/S^\mr{log}})$ mapped to $1 \in \Gamma (X, \mcO_X)$
 (i.e., $\nabla$ forms a $\mfg$-Cartan connection on $\mfX$).
$\nabla$ specifies a section $S \migi \mfC \overline{\mfC} \mfo_{\mfg, \mfX}$, and 
this completes the proof of assertion (ii).

Finally, assertion (iii) follows from Proposition \ref{prfg001}  (ii), which implies that after base-changing via the finite and faithfully flat morphism $\mfC \mfC \mfo_{\mfg, g, 0}^{^\mr{Zzz...}} \migi \mfM_{g,0}$ (cf. (\ref{eqeq14}) for the definition of  $\mfC \mfC \mfo_{\mfg, g, 0}^{^\mr{Zzz...}}$),
the projection $\mfC \mfC \mfo_{\mfg, g, 0} \migi \mfM_{g,0}$  admits a global section 
(i.e., the section corresponding to the closed immersion $\mfC \mfC \mfo_{\mfg, g, 0}^{^\mr{Zzz...}} \migi \mfC \mfC \mfo_{\mfg, g, 0}$).
\end{proof}

\vspace{5mm}
\subsection{$\mr{GL}_n$-Cartain connections} \label{sc020}
\leavevmode\\ \vspace{-4mm}

In this subsection, 
we restrict ourselves  to the case where $\mbG = \mr{PGL}_n$.
(In particular, either one of the two conditions $(\mr{Char})_{0}$,  $(\text{Char})_{p}^{\mfs \mfl}$ is satisfied.)
 Let  $S$ be a $k$-scheme and $\mfX := (X/S, \{ \sigma_i \}_{i =1}^r)$ an $r$-pointed stable curve of genus $g$  over $S$.
For each line bundle  $\mcN$ on $X$ and each  integer $l$, we shall write
\begin{align} \label{f050}
\mcF_{\mcN, l}^\dagger := \mcT_{X^\mr{log}/S^\mr{log}}^{\otimes l} \otimes \mcN  
\end{align} 
(hence, $\mcF^{\dagger}_{\mcN, 0} = \mcN$).
Also,   write
\begin{align} \label{eq0F28}
\mcF_{\mcN}^{[n] \dagger} := \bigoplus_{l=0}^{n-1} \mcF_{\mcN, l}^\dagger.
\end{align}
In particular, the $\mcO_X$-module  $\mcF_{\mcN}^{[n] \dagger}$ admits, by definition,  a natural   grading.

\bde  \label{z012}\leavevmode\\
\vspace{-5mm}
\begin{itemize}
\item[(i)]
  A {\bf $\mr{GL}_n$-Cartan connection} on $\mfX$ is
  a log  flat bundle on $X^\mr{log}/S^\mr{log}$ 
of the form 
\begin{align}
\widehat{\mcF}^\clubsuit := (\mcF^{[n]\dagger}_{\mcN}, \bigoplus_{l=0}^{n-1}  \nabla_l) \ \left(= \bigoplus_{l=0}^{n-1} (\mcF^\dagger_{\mcN, l}, \nabla_l) \right), 
\end{align}
where $\mcN$ is a line bundle on $X$ and   each $\nabla_l$
 ($l =0, \cdots, n-1$) is an $S^\mr{log}$-connection on 
 the $l$-th component 
  $\mcF_{\mcN, l}^{\dagger}$ 
  of   $\mcF_{\mcN}^{[n]\dagger}$.

\item[(ii)]
 Let $\mcN^\flat := (\mcN, \nabla_\mcN)$ be a log flat  line bundle on $X^\mr{log}/S^\mr{log}$.
   A
 {\bf $(\mr{GL}_n, \mcN^\flat)$-Cartan connection} on $\mfX$ is 
  a $\mr{GL}_n$-Cartan connection $\widehat{\mcF}^\clubsuit := (\mcF^{[n]\dagger}_{\mcN}, \bigoplus_{l=0}^{n-1}  \nabla_l)$ on $\mfX$
   satisfying the equality
$(\mcF^\dagger_{\mcN, 0}, \nabla_0) = \mcN^\flat$.
(Hence, any $\mr{GL}_n$-Cartan connection is a $(\mr{GL}_n, \mcN^\flat)$-Cartan connection for some log flat line bundle $\mcN^\flat$.)
\end{itemize}
\ede
\vspace{3mm}

\begin{rema} \label{Rf010}
\leavevmode\\
 \ \ \ 
  Let  $\widehat{\mcF}^\clubsuit := (\mcF_{\mcN}^{[n]\dagger}, \bigoplus_{l=0}^{n-1} \nabla_l)$ be 
  a $(\mr{GL}_n, \mcN^\flat)$-Cartan connection on $\mfX$ for
   some  log flat line bundle $\mcN^\flat = (\mcN, \nabla_\mcN)$.
   If  $\mcL^\flat := (\mcL, \nabla_\mcL)$ is a log flat line bundle on $X^\mr{log}/S^\mr{log}$,
then 
$\mcF_{\mcN \otimes \mcL}^{[n]\dagger}$ may be  canonically identified with 
$\mcF_{\mcN}^{[n]\dagger} \otimes \mcL$ and, via this identification,  
 the tensor product 
\begin{align} \label{eq020}
\widehat{\mcF}^\clubsuit \otimes \mcL^\flat = (\mcF_{\mcN}^{[n]\dagger} \otimes \mcL, \bigoplus_{l=0}^{n-1} \nabla_{l} \otimes \nabla_\mcL )
\end{align}
forms
a $(\mr{GL}_n, \mcN^\flat \otimes \mcL^\flat)$-Cartan connection.
 \end{rema}
\vspace{3mm}

Denote by $\mbT_{\mr{GL}_n}$ the maximal torus of $\mr{GL}_n$ consisting of diagonal matrices.
We take   the maximal torus  $\mbT_{\mr{PGL}_n}$  of $\mr{PGL}_n$
 to be the image of diagonal matrices via the quotient $\mr{GL}_n \migisurj \mr{PGL}_n$.
In particular, we obtain a natural projection $\mbT_{\mr{GL}_n} \migisurj \mbT_{\mr{PGL}_n}$ and    an   isomorphism $\mbG_m^{\times n} \isom \mbT_{\mr{GL}_n}$ (where $\mbG_m^{\times n}$ denotes the product over $k$ of $r$ copies of $\mbG_m$) which, to any element $(a_1, \cdots , a_n) \in \mbG_m^{\times n}$, assigns
the diagonal matrix in $\mbT_{\mr{GL}_n}$ with entries $a_1, \cdots, a_n$.

Now, let us take   a $\mr{GL}_n$-Cartan connection
$\widehat{\mcF}^\clubsuit$ on $\mfX$, which 
determines (since it is a direct sum of $n$ log flat line bundles)  a log flat $\mbG_m^{\times n}$-torsor. 
One verifies that its underlying  $\mbG_m^{\times n}$-torsor determines  the $\mbT_{\mr{PGL}_n}$-torsor  $\mcE^\dagger_{\mbT_{\mr{PGL}_n}}$ by 
 change of structure group via  the composite  $\mbG_m^{\times n} \isom \mbT_{\mr{GL}_n} \migisurj \mbT_{\mr{PGL}_n}$.
 That is to say, 
 $\widehat{\mcF}^\clubsuit$  induces  
 an $\mfs \mfl_n \ (=\mfp \mfg \mfl_n)$-Cartan connection $\nabla^{\widehat{\mcF}^\flat}$ on  $\mfX$.

\vspace{3mm}
\bpr \label{Pf004} \leavevmode\\
 \vspace{-5mm}
\begin{itemize}
\item[(i)]
 The following maps of sets are bijective:
   \begin{align} \label{efegg0991}
& \begin{pmatrix}
\text{the set  of $S^\mr{log}$-connections} \\
\text{on $\Omega_{X^\mr{log}/S^\mr{log}}$} 
\end{pmatrix}^{\times (n-1)}  \\
& \hspace{38mm} \stackrel{\sim}{\longmigi} \begin{pmatrix}
\text{the set  of $(\mr{GL}_n, \mcO_X^\flat)$-Cartan} \\
\text{connections on $\mfX$} 
\end{pmatrix} \notag  \\
& \hspace{70mm}
\stackrel{\sim}{\longmigi} 
\begin{pmatrix}
\text{the set  of $\mfs \mfl_n$-Cartan} \\
\text{connections on $\mfX$}\\
\end{pmatrix}, \notag 
\end{align}
where the first and second  maps  are  given by assigning 
$(\nabla_l)_{l=1}^{n-1} \mapsto (\mcF_{\mcO_X}^{[n]\dagger}, \bigoplus_{l=0}^{n-1} (\bigotimes_{j=1}^{l} \nabla_j^\vee ))$
 (where $\bigotimes_{j=1}^0\nabla_j^\vee := d$)
 and 
 $\widehat{\mcF}^\clubsuit \mapsto \nabla^{\widehat{\mcF}^\flat}$ respectively.
Moreover, these bijections are 
functorial with respect to $S$.
\item[(ii)]
The composite of  the bijections  in 
 (\ref{efegg0991})
    induces an isomorphism
 \begin{align}\label{eq070}
 \overline{\mfC} \mfo_{g,r}^{\times (n-1)} \isom \mfC \overline{\mfC} \mfo_{\mfs \mfl_n, g,r}
 \end{align}
over $\overline{\mfM}_{g,r}$, where the left-hand side denotes the product  over $\overline{\mfM}_{g,r}$ of $n-1$ copies of $\overline{\mfC} \mfo_{g,r}$.
\end{itemize}
\epr
\begin{proof}
The assertions follow immediately  from the various definitions involved.
 \end{proof}

\vspace{3mm}

Next, denote by 
$\mft_{\mr{GL}_n}$ and  $\mft_{\mr{PGL}_n}$ the Lie algebras of 
$\mbT_{\mr{GL}_n}$ and  $\mbT_{\mr{PGL}_n}$ respectively.
The set $\mft_{\mr{GL}_n} (S)$ of $S$-rational points of $\mft_{\mr{GL}_n}$ may be  identified, via the isomorphism $\mbG_m^{\times n} \isom \mbT_{\mr{GL}_n}$ mentioned above,  with
$\Gamma (S, \mcO_S)^{\oplus  n}$.
Also, the set   $\mft_{\mr{PGL}_n} (S)$ may be identified with $\mr{Coker} (\varDelta)$, where 
$\varDelta$ denotes the diagonal embedding $\Gamma (S, \mcO_S) \migiincl \Gamma (S, \mcO_S)^{\oplus  n}$.
If
\begin{align} \label{eq114}
\pi :  \mft_{\mr{GL}_n} (S) \migisurj \mft_{\mr{PGL}_n} (S)
\end{align}
denotes the quotient arising from the natural quotient  $\mr{GL}_n \migisurj \mr{PGL}_n$,
then it may be identified, under the identifications just discussed,  with the natural quotient $\Gamma (S, \mcO_S)^{\oplus  n} \migisurj \mr{Coker} (\varDelta)$.
For each  $\vec{\mu} := ((\mu_{i1}, \mu_{i2}, \cdots, \mu_{in}))_{i =1}^r \in (\Gamma (S, \mcO_S)^{\oplus n})^{\times r}$ ($= \mft_{\mr{GL}_n} (S)^{\times r}$), we shall  write
\begin{align} \label{ee300}
\pi (\vec{\mu}) := (\pi (\mu_{i1}, \mu_{i2}, \cdots, \mu_{in}))_{i=1}^r \in  \mft_{\mr{PGL}_n} (S)^{\times r}.
\end{align}
If $\vec{\mu} := \emptyset$, then we write $\pi (\vec{\mu}) := \emptyset$.
 By the construction of  (\ref{eq070}), the following proposition holds.

\vspace{3mm}
\bpr \label{P174p2}\leavevmode\\
 \ \ \ 
Assume that $r >0$.
\begin{itemize}
\item[(i)]
  Let   $\vec{\mu} := ((\mu_{i1}, \mu_{i2}, \cdots, \mu_{in}))_{i =1}^r \in (\Gamma (S, \mcO_S)^{\oplus n})^{\times r}$.
 Then, the bijections in (\ref{efegg0991}) restrict to the following bijections:
   \begin{align} \label{ee290}
& \prod_{l=1}^{n-1} \begin{pmatrix}
\text{the set  of $S^\mr{log}$-connections} \\
\text{on $\Omega_{X^\mr{log}/S^\mr{log}}$ of } \\
\text{monodromies $(\mu_{il} - \mu_{i (l+1)})_{i=1}^r$} 
\end{pmatrix} \\
& \hspace{38mm} \stackrel{\sim}{\longmigi} \begin{pmatrix}
\text{the set  of $(\mr{GL}_n, \mcO_X^\flat)$-Cartan} \\
\text{connections on $\mfX$ of monodromies} \\
\text{$((0, \mu_{i2} - \mu_{i1}, \cdots, \mu_{i n}-\mu_{i1}))_{i=1}^r$} 
\end{pmatrix} \notag  \\
& \hspace{70mm}
\stackrel{\sim}{\longmigi} 
\begin{pmatrix}
\text{the set  of $\mfs \mfl_n$-Cartan} \\
\text{connections on $\mfX$}\\
\text{of monodromies $\pi (\vec{\mu})$}
\end{pmatrix}. \notag 
\end{align}
\item[(ii)]
 Let 
  $\vec{\mu} := ((\mu_{i1}, \mu_{i2}, \cdots, \mu_{in}))_{i =1}^r \in (k^{\oplus n})^{\times r}$.
  Then, 
 the    composite  of the bijections in   (\ref{ee290}) induces
   an isomorphism
   \begin{align} \label{eq071}
   \prod_{l=1}^{n-1} \overline{\mfC} \mfo_{g,r, (\mu_{il} -\mu_{i(l+1)})_{i=1}^{r}} \isom \mfC  \overline{\mfC} \mfo_{\mfs \mfl_n, g,r, \pi(\vec{\mu})}
   \end{align}
over $\overline{\mfM}_{g,r}$,  where  the left-hand side is  the product over  $\overline{\mfM}_{g,r}$.
\end{itemize}
 \epr


\vspace{5mm}
\subsection{}{\bf The relative  affine spaces $\ooalign{$\bigoplus$ \cr $\bigotimes_{\Box, g,r}^{}$}$ over $\overline{\mfM}_{g,r}$.}  \label{sde1d3}
\leavevmode\\ \vspace{-4mm}

Now, let us go back to the general case but,
in the rest of this section, we assume that
  either one of two conditions $(\text{Char})_{p}$,  $(\text{Char})_{p}^{\mfs \mfl}$ is satisfied.
  Moreover, assume that $\mbG$, $\mbB$, and $\mbT$ are all defined over $\mbF_p$.
  In particular, it makes sense to speak of the subset $\mft (\mbF_p)$ of $\mft (k)$ (cf. Proposition \ref{pr0011}).
We shall write
\begin{align} \label{eq090}
\ooalign{$\bigoplus$ \cr $\bigotimes_{\mfb, g,r}$} 
&:= \mbV (\overline{f}_{\mft \mfa \mfu*} (\Omega_{\overline{\mfC}^\mr{log}_{g,r}/\overline{\mfM}_{g,r}^\mr{log}}^{\otimes p} \otimes \mfb_{\mcE^\dagger})), 
\\
\ooalign{$\bigoplus$ \cr $\bigotimes_{\mft, g,r}$} &:= \mbV (\overline{f}^{}_{\mft \mfa \mfu*} (\Omega_{\overline{\mfC}_{g,r}^{ \mr{log}}/\overline{\mfM}_{g,r}^\mr{log}}^{\otimes p} \otimes_k \mft)) \ \left(= \mbV (\overline{f}^{}_{\mft \mfa \mfu*} (\Omega^{\otimes p}_{\overline{\mfC}_{g,r}^{\mr{log}}/\overline{\mfM}_{g,r}^\mr{log}})) \times_k \mft \right), \notag 
\end{align}
where  $\mcE^\dagger := (\Omega_{\overline{\mfC}^\mr{log}_{g,r}/\overline{\mfM}_{g,r}^\mr{log}})^{\times} \times^{\mbG_m, \check{\rho}} \mbT$ (cf. (\ref{f016})).
Also, we shall write 
\begin{equation}  \label{spindle}
   \ooalign{$\bigoplus$ \cr $\bigotimes_{\mfc, g,r}$}
   \end{equation}
for the relative  scheme over $\overline{\mfM}_{g,r}$  (cf. ~\cite{FGA}, Theorem 5.23) representing  the set-valued contravariant functor on $\mfS \mfc \mfh_{/\overline{\mfM}_{g,r}}$  which, to 
each object $S \migi \overline{\mfM}_{g,r}$ of $\mfS \mfc \mfh_{/\overline{\mfM}_{g,r}}$ classifying pointed stable curve  $\mfX := (X/S, \{ \sigma_i \}_{i=1}^r)$,
 assigns the set of  morphisms  
   $X  \migi (\Omega_{X^\mr{log}/S^\mr{log}}^{\otimes p})^\times \times^{\mbG_m}\mfc$ over $X$.
 The composite $\mfb \migiincl \mfg \stackrel{\chi}{\migisurj} \mfc$ and the projection $\mft \migisurj \mfc$ induce    morphisms
\begin{align}  \label{spindle2}
  \ooalign{$\oplus$ \cr $\otimes_{\mfb \migi \mfc, g,r}$} :    \ooalign{$\bigoplus$ \cr $\bigotimes_{\mfb, g,r}$}   \migi \ooalign{$\bigoplus$ \cr $\bigotimes_{\mfc, g,r}$}
   \ \  \text{and} \ \ 
 \ooalign{$\oplus$ \cr $\otimes_{\mft  \migi \mfc, g,r}$} : \ooalign{$\bigoplus$ \cr $\bigotimes_{\mft, g,r}$}   \migi \ooalign{$\bigoplus$ \cr $\bigotimes_{\mfc, g,r}$}
 \end{align}
respectively.
Also, we shall write
\begin{align} \label{EQEQ1}
[0]_{\mft, g,r} : \overline{\mfM}_{g,r} \migi   \ooalign{$\bigoplus$ \cr $\bigotimes_{\mft, g,r}$}  \ \text{and} \    [0]_{\mfc, g,r} \ \left(:=  \ooalign{$\oplus$ \cr $\otimes_{\mft  \migi \mfc, g,r}$} \circ [0]_{\mft, g,r}\right): \overline{\mfM}_{g,r} \migi \ooalign{$\bigoplus$ \cr $\bigotimes_{\mfc, g,r}$}
\end{align}
for the zero sections.
\vspace{5mm}
\subsection{Cartan connections with vanishing $p$-curvature}
  \label{se13}
\leavevmode\\ \vspace{-4mm}

We shall introduce the definitions of a dormant  Cartan connections and  a $p$-nilpotent 
Cartan connection (cf. Definition \ref{j01dd}).

The assignment from each $\mfg$-Cartan connection  to its $p$-curvature determines a morphism
\begin{align} \label{eq010222}
\Psi_{\mfg, g,r} :   \mfC \overline{\mfC} \mfo_{\mfg, g,r} \migi \ooalign{$\bigoplus$ \cr $\bigotimes_{\mft, g,r}$}
\end{align}
over $\overline{\mfM}_{g,r}$, and
we obtain the  composite
\begin{align} \label{eq0102}
\Psi_{\mfg, g,r}^\chi : \mfC \overline{\mfC} \mfo_{\mfg, g,r} \xrightarrow{\Psi_{\mfg, g,r}} \ooalign{$\bigoplus$ \cr $\bigotimes_{\mft, g,r}$} \xrightarrow{\ooalign{$\oplus$ \cr $\otimes_{\mft \migi \mfc, g,r}$}} \ooalign{$\bigoplus$ \cr $\bigotimes_{\mfc, g,r}$}.
\end{align}
Let us  write
\begin{align} \label{eqeq14}
\mfC \overline{\mfC} \mfo_{\mfg, g,r}^{^\mr{Zzz...}} \ \ (\text{resp.,} \ \mfC \overline{\mfC} \mfo_{\mfg, g,r}^{^{p\text{-} \mr{nilp}}}) 
\end{align}
for the  inverse image via  $\Psi_{\mfg, g,r}$ (resp.,  $\Psi_{\mfg, g,r}^\chi$) of the zero section $[0]_{\mft, g,r}$
  (resp., $[0]_{\mft, g,r}$).

\vspace{3mm}
\bde  \label{j01dd}\leavevmode\\
 \ \ \ 
We shall say that a $\mfg$-Cartan connection is {\bf dormant} (resp., {\bf $p$-nilpotent}) if it is classified by the closed substack 
$\mfC \overline{\mfC} \mfo_{\mfg, g,r}^{^\mr{Zzz...}}$ (resp., $\mfC \overline{\mfC} \mfo_{\mfg, g,r}^{^{p\text{-} \mr{nilp}}}$).
(In particular, a dormant $\mfg$-Cartan connection is a $\mfg$-Cartan connection with vanishing $p$-curvature.)
 \ede

\bpr \label{pr0011} \leavevmode\\
 \ \ \ 
Assume that $r >0$.
Let $S$ be a  connected (resp.,  reduced and connected) $k$-scheme, 
$\mfX := (X/S, \{ \sigma_i \}_{i=1}^r)$ an $r$-pointed stable curve over $S$ of genus $g$, and $\nabla$ a  dormant (resp., $p$-nilpotent) $\mfg$-Cartan connection  on $\mfX$.
Then, 
the monodromy $\mu_i^{(\mcE^\dagger_\mbT, \nabla)} \in \mft (S)$ of $\nabla$ at any marked point $\sigma_i$ ($i =1, \cdots, r$) lies in $\mft (\mbF_p)$.
\epr
\begin{proof}
Let us fix  $i \in \{1, \cdots,r\}$.
The restriction $\sigma^*_i (\psi^{(\mcE^\dagger_\mbT, \nabla)})$ of the $p$-curvature $\psi^{(\mcE^\dagger_\mbT, \nabla)} \in \Gamma (X, \Omega_{X^\mr{log}/S^\mr{log}}^{\otimes p} \otimes \mft_{\mcE^\dagger_\mbT}) \  (= \Gamma (X, \Omega_{X^\mr{log}/S^\mr{log}}^{\otimes p} \otimes_k \mft))$
  may be thought of as an element of
 $\mft (S)$
    via the following  composite isomorphism 
   \begin{align}
   \Gamma (S, \sigma^*_i (\Omega_{X^\mr{log}/S^\mr{log}}^{\otimes p} \otimes_k \mft)) & \isom \Gamma (S, \sigma^*_i (\Omega_{X^\mr{log}/S^\mr{log}})^{\otimes p} \otimes_k \mft)   \\
 & \isom \Gamma (S, \mcO_S \otimes_k \mft) \notag \\
 & \isom \mft (S), \notag
\end{align}
where the second isomorphism arises from  $\mft \mfr \mfi \mfv_{\mfX, i}$ (cf. (\ref{f029})).
By the definition of $p$-curvature, 
 the  following equality holds:
\begin{align} \label{f031}
\sigma^*_i (\psi^{(\mcE^\dagger_\mbT, \nabla)}) = (\mu_i^{(\mcE^\dagger_\mbT, \nabla)})^{[p]} -\mu_i^{(\mcE^\dagger_\mbT, \nabla)}, 
\end{align}
where $(\mu_i^{(\mcE^\dagger_\mbT, \nabla)})^{[p]}$ denotes the image of $\mu_i^{(\mcE^\dagger_\mbT, \nabla)} \ \left(\in \mft (S)\right)$ via  the $p$-th power operation $(-)^{[p]}$  on  the Lie algebra $\mft$ (cf. ~\cite{Wak5}, Remark 3.2.2).
Since $\mft$ (considered as a scheme)
 is defined over $\mbF_p$,
it may be identified with $\mft^{(1)}_k$ via the isomorphism $\mr{id}_\mft \times F_{\mr{Spec}(k)} : \mft^{(1)}_k \isom \mft$.
Under this identification, the $p$-th power operation on $\mft$ coincides with the relative Frobenius morphism  $F_{\mft / k} : \mft \migi \mft \ (= \mft_k^{(1)})$.
Hence, (\ref{f031}) induces the equality
\begin{align} \label{eqeq59}
\sigma^*_i (\psi^{(\mcE^\dagger_\mbT, \nabla)}) = F_{\mft/k} \circ  \mu_i^{(\mcE^\dagger_\mbT, \nabla)} - \mu_i^{(\mcE^\dagger_\mbT, \nabla)}.
\end{align}
Since $F_{\mft/k} \circ  \mu_i^{(\mcE^\dagger_\mbT, \nabla)} = \mu_i^{(\mcE^\dagger_\mbT, \nabla)}$ if and only if 
$\mu_i^{(\mcE^\dagger_\mbT, \nabla)} \in \mft (\mbF_p)$, 
the non-resp'd assertion  follows from (\ref{eqeq59}).

Next,  let us consider the resp'd assertion.
Since $\chi \circ \sigma_i^*(\psi^{(\mcE^\dagger_\mbT, \nabla)}) = \chi (0)$,
the element $\sigma_i^*(\psi^{(\mcE^\dagger_\mbT, \nabla)})$ of $\mft (S)$ lies in
$(\mft \times_{\mfc, \chi (0)} \mr{Spec} (k))(S)$.
But, the closed subscheme  $(\mft \times_{\mfc, \chi (0)} \mr{Spec} (k))_{\mr{red}}$ of $\mft$ is isomorphic to the point $0 \in \mft (k)$.
Hence, under 
the assumetion that $S$ is reduced and connected,
the equality  $\sigma_i^*(\psi^{(\mcE^\dagger_\mbT, \nabla)}) = 0$ holds.
By applying the same discussion as above, we have  $\mu_i^{(\mcE^\dagger_\mbT, \nabla)} \in \mft (\mbF_p)$.
This completes  the proof of the resp'd assertion.
\end{proof}

\vspace{3mm}
\bco \label{prff001}
\leavevmode\\
 \ \ \ 
 Let $\vec{\mu} \in \mft (k)^{\times r}$.
 Then, both $\mfC \overline{\mfC} \mfo_{\mfg, g,r, \vec{\mu}}^{^\mr{Zzz...}}$ and 
  $\mfC \overline{\mfC} \mfo_{\mfg, g,r, \vec{\mu}}^{^{p \text{-} \mr{nilp}}}$ are empty unless $\vec{\mu}$ lies in $\mft (\mbF_p)^{\times r}$.
  Moreover, $\mfC \overline{\mfC} \mfo_{\mfg, g,r}^{^\mr{Zzz...}}$ decomposes as
  \begin{align} \label{ee260}
  \mfC \overline{\mfC} \mfo_{\mfg, g,r}^{^\mr{Zzz...}} = \coprod_{\vec{\mu} \in \mft (\mbF_p)^{\times r}} \mfC \overline{\mfC} \mfo_{\mfg, g,r, \vec{\mu}}^{^\mr{Zzz...}}.
  \end{align}
 \eco

\vspace{5mm}
\subsection{}{\bf Structure of the moduli stacks $\mfC \mfC \mfo_{\mfg, g,r,  \vec{\mu}}^{^\mr{Zzz...}}$,  
  $\mfC \mfC \mfo_{\mfg, g,r, \vec{\mu}}^{^{p \text{-} \mr{nilp}}}$.}  \label{sfe13}
\leavevmode\\ \vspace{-4mm}

For each $\vec{\mu} \in \mft (k)^{\times r}$ (where $\vec{\mu} := \emptyset$ if $r =0$), we shall write
\begin{align} \label{ee240}
\mfC \mfC \mfo_{\mfg, g,r,  \vec{\mu}}^{^\mr{Zzz...}} := \mfC \overline{\mfC} \mfo_{\mfg, g,r,  \vec{\mu}}^{^\mr{Zzz...}} \times_{\overline{\mfM}_{g,r}} \mfM_{g,r}, \hspace{5mm}
\mfC \mfC \mfo_{\mfg, g,r, \vec{\mu}}^{^{p \text{-} \mr{nilp}}} := \mfC \overline{\mfC} \mfo_{\mfg, g,r, \vec{\mu}}^{^{p \text{-} \mr{nilp}}} \times_{\overline{\mfM}_{g,r}} \mfM_{g,r}.
\end{align}

\vspace{3mm}
\bpr\label{prfg001}
\leavevmode\\
 \ \ \ 
Let $\vec{\mu} \in \mft (\mbF_p)^{\times r}$ (where $\vec{\mu} := \emptyset$ if $r =0$).
\begin{itemize}
\item[(i)]
 $\mfC \mfC \mfo^{^\mr{Zzz...}}_{\mfg, g,r, \vec{\mu}}$  may be represented by either the empty stack or a Deligne-Mumford stack over $k$ which  is finite and faithfully flat over $\mfM_{g,r}$ of degree $p^{g \cdot \mr{rk} (\mfg)}$.
 If $\mfC \mfC \mfo^{^\mr{Zzz...}}_{\mfg, g,r,  \vec{\mu}} \neq \emptyset$, then  the open substack  $\mfC \mfC \mfo^{^\mr{Zzz...}}_{\mfg, g,r, \vec{\mu}} \times_{\mfM_{g,r}} \mfM_{g,r}^{\mr{ord}}$ of $\mfC \mfC \mfo^{^\mr{Zzz...}}_{\mfg, g,r, \vec{\mu}}$  is \'{e}tale over $\mfM_{g,r}^{\mr{ord}}$.
 \item[(ii)]
  If $r =0$ and $p | 2g-2$, then  $\mfC \mfC \mfo^{^\mr{Zzz...}}_{\mfg, g,0,  \emptyset}$ is nonempty.
 \end{itemize}
 \epr
\begin{proof}
Let us consider assertion (i).
For each integer $m$, we shall write $\mcL^m := \Omega_{\mfC_{g,r}^\mr{log}/\mfM_{g,r}^\mr{log}}^{\otimes m}$.
Observe  that 
the $\mbT$-torsor $\mcE^\dagger_\mbT$  for the tautological family of curves $\mfC_{g,r}/\mfM_{g,r}$ is 
isomorphic to the product of $\mbG_m$-torsors corresponding to 
$\mcL^{m_j}$'s
 (where $j =1, \cdots, \mr{rk} (\mfg)$ and each $m_j$ is an integer).
Hence,  there exists an isomorphism $\mfC \mfC \mfo_{\mfg, g,r} \isom \prod_{j =1}^{\mr{rk} (\mfg)}  \mfC \mfo_{\mcL^{m_j}, g,r}$.
This isomorphism restricts to an isomorphism 
\begin{align} \label{ee270}
\mfC \mfC \mfo^{^\mr{Zzz...}}_{\mfg, g,r, \vec{\mu}} \isom \prod_{j =1}^{\mr{rk} (\mfg)}  \mfC \mfo^{\psi =0}_{\mcL^{m_j}, g,r, \vec{n}_j},
\end{align}
 where each $\vec{n}_j$ is an element of $\mbF_p^{\times r}$ (cf. Remark \ref{zz050}).
Thus, the assertion follows from Proposition \ref{prfg001tt} (ii).
Assertion (ii) follows from (\ref{ee270}) and Proposition \ref{prfg001tt} (i) (together with the fact that $(2g-2) | \mr{deg} (\mcL^{m_j})$ for any $j =1, \cdots, \mr{rk} (\mfg)$).
\end{proof}

\vspace{3mm}
\bco \label{P230}\leavevmode\\
\ \ \ 
Let $\vec{\mu} \in \mft (\mbF_p)^{\times r}$ (where $\vec{\mu} := \emptyset$ if $r =0$).
Then,  $\mfC \mfC \mfo_{\mfg,  g,r, \vec{\mu}}^{^{p\text{-} \mr{nilp}}}$ may be represented by  a (possibly empty) Deligne-Mumford stack over $k$ which is finite over $\mfM_{g,r}$.
  \eco
\begin{proof}
Let  $\mfN$  be the closed substack of $\ooalign{$\bigoplus$ \cr $\bigotimes_{\mft, g,r}$}$
defined to be 
$\ooalign{$\oplus$ \cr $\otimes^{-1}_{\mft \migi \mfc, g,r}$} ([0]_{\mfc, g,r})$.
The morphism $\Psi_{\mfg, g,r}$ (cf.   (\ref{eq010222})) restricts to a morphism $\mfC \mfC \mfo_{\mfg,  g,r, \vec{\mu}}^{^{p\text{-} \mr{nilp}}} \migi \mfN$.
Here, observe that the reduced stack $\mfN_{\mr{red}}$ associated with $\mfN$ is isomorphic to the closed substack $[0]_{\mft, g,r}$.
It follows that  
\begin{align} \label{FGHJ}
\mfC \mfC \mfo_{\mfg, g,r, \vec{\mu}}^{^{p\text{-} \mr{nilp}}} \times_{\mfN} \mfN_{\mr{red}} \cong \mfC \mfC \mfo_{\mfg, g,r, \vec{\mu}}^{^\mr{Zzz...}}.
\end{align}
Since $\mfC \mfC \mfo_{\mfg, g,r, \vec{\mu}}^{^{p\text{-} \mr{nilp}}}$ is of finite type over $\mfM_{g,r}$, 
(\ref{FGHJ}) and Proposition \ref{prfg001} (i) imply that $\mfC \mfC \mfo_{\mfg, g,r,  \vec{\mu}}^{^{p\text{-} \mr{nilp}}}$ is finite over $\mfM_{g,r}$.
\end{proof}

Also, in the case where $\mbG = \mr{PGL}_n$ (hence $\mfg = \mfp \mfg \mfl_n = \mfs \mfl_n$), the following proposition holds.

\vspace{3mm}
\bpr \label{prff001hy}
\leavevmode\\
\vspace{-5mm}
\begin{itemize}
\item[(i)]
We have  the following commutative square diagram:
\begin{align} \label{ee250}
\begin{CD}
(\mfC \mfo_{g,r}^{\psi =0})^{\times (n-1)} @> \sim >> \mfC \mfC \mfo^{^\mr{Zzz...}}_{\mfs \mfl_n,g,r}
\\
@V \mr{open \,  imm.} VV @VV \mr{open \,  imm.}  V
\\
(\overline{\mfC} \mfo_{g,r}^{\psi =0})^{\times (n-1)} @>> \sim > \mfC \overline{\mfC} \mfo^{^\mr{Zzz...}}_{\mfs \mfl_n,g,r},
\end{CD}
\end{align}
where the upper and lower horizontal arrows  are isomorphisms  obtained by restricting  (\ref{eq070}).
\item[(ii)]
Assume further that $r >0$, and 
let $\vec{\mu}$, $\pi (\vec{\mu})$ be as in Proposition \ref{P174p2}.
  Then, by restricting square diagram (\ref{ee250}) above, we have the following commutative square diagram:
  \begin{align} \label{ee251}
  \begin{CD}
  \prod_{l=1}^{n-1} \mfC \mfo^{\psi =0}_{g,r, (\mu_{il} -\mu_{i(l+1)})_{i=1}^{r}} @> \sim >> \mfC  \mfC \mfo^{^\mr{Zzz...}}_{\mfs \mfl_n, g,r, \pi(\vec{\mu})}
 \\
 @V  \mr{open \,  imm.} VV @VV  \mr{open \,  imm.} V
 \\
  \prod_{l=1}^{n-1} \overline{\mfC} \mfo^{\psi =0}_{g,r, (\mu_{il} -\mu_{i(l+1)})_{i=1}^{r}} @>> \sim > \mfC  \overline{\mfC} \mfo^{^\mr{Zzz...}}_{\mfs \mfl_n, g,r, \pi(\vec{\mu})},
  \end{CD}
  \end{align}
where the upper and lower horizontal arrows are isomorphisms  and  the two stacks in the left-hand side are   the products over  $\overline{\mfM}_{g,r}$.
In particular, if $\mu_{i l} \in \mbF_p$ for any pair $(i, l)$ and 
  $p | (2g-2 + \sum_{i=1}^r \tau^{-1}(\mu_{il}-\mu_{i (l+1)}))$ for any  $l = 1, \cdots, n-1$, 
 then  
 both $\mfC  \overline{\mfC} \mfo^{^\mr{Zzz...}}_{\mfs \mfl_n, g,r, \pi(\vec{\mu})}$ and $\mfC  \mfC \mfo^{^\mr{Zzz...}}_{\mfs \mfl_n, g,r, \pi(\vec{\mu})}$
 are nonempty.
\end{itemize}
 \epr
\begin{proof}
The assertions follow from Proposition \ref{prfg001tt} (i) and  the construction of
  (\ref{eq071}).
 (Indeed, if two connections $\nabla_1$ and $\nabla_2$ have vanishing $p$-curvature, then the tensor product $\nabla_1 \otimes \nabla_2^\vee$ of $\nabla_1$ and the dual $\nabla_2^\vee$ of $\nabla_2$   has vanishing $p$-curvature.)
\end{proof}

\vspace{10mm}
\section{Miura $\mfg$-opers on pointed stable curves} \label{Sec4} \vspace{3mm}

In this section, we discuss the definition and some basic properties of {\it Miura $\mfg$-opers} on a family of pointed stable curves.
It will be shown   that generic Miura $\mfg$-opers correspond bijectively to $\mfg$-Cartan connections.
(This result may be thought of as a global version of ~\cite{Fr}, Proposition 8.2.2.)
In particular,  the moduli functor,  denoted by $\mfM \overline{\mfO} \mfp_{\mfg, g,r, (\vec{\varepsilon})}$ (cf. (\ref{eqeq100}) and (\ref{ee200})), classifying  pointed stable curves paired  with a generic Miura $\mfg$-oper (of prescribed exponents $\vec{\varepsilon}$) may be represented by a Deligne-Mumford stack (cf. Proposition \ref{eqeq88}).
In the case of positive characteristic,
we introduce two kinds of   Miura $\mfg$-opers,
called 
 {\it dormant generic Miura $\mfg$-opers} and {\it $p$-nilpotent generic Miura $\mfg$-opers}.
The bijective correspondence mentioned above restricts to a correspondence between $p$-nilpotent Miura $\mfg$-opers and $p$-nilpotent $\mfg$-Cartan connections (cf. Theorem  \ref{th013}).

\vspace{5mm}
\subsection{$\mfg$-opers} \label{y033}
\leavevmode\\
\vspace{-4mm}

Let us keep the notation and assumptions in \S\,\ref{SSQ1}.
First, we shall recall the definition of a $\mfg$-oper (on a log-curve).
Let  $T^\mr{log}$ be an fs log scheme over $k$, $U^\mr{log}$ a log-curve over $T^\mr{log}$, and 
$\pi_\mbB : \mcE_\mbB  \migi U$  a right $\mbB$-torsor over $U$.
Denote by
$\pi_\mbG : (\mcE_\mbB \times^\mbB \mbG =:) \ \mcE_\mbG \migi U$ the right $\mbG$-torsor over $U$ obtained by 
    change of structure group via the inclusion $\mbB \migiincl \mbG$.
For each $j \in \mbZ$,  $\mfg^j$ 
is   closed under  the adjoint action of $\mbB$ on $\mfg$,  and hence, induces
  an $\mcO_U$-submodule $\mfg^j_{\mcE_\mbB}$
   of $\mfg_{\mcE_\mbB}$ ($\subseteq \widetilde{\mcT}_{\mcE_\mbG^\mr{log}/T^\mr{log}}$). 
If we write 
\begin{align}
 \widetilde{\mcT}_{\mcE_\mbG^\mr{log}/T^\mr{log}}^j :=  \iota(\widetilde{\mcT}_{\mcE_\mbB^\mr{log}/T^\mr{log}}) + \mfg^{j}_{\mcE_\mbB}
 \end{align}
($j \in \mbZ$),   where $\iota$ denotes the injection $\widetilde{\mcT}_{\mcE_\mbB^\mr{log}/T^\mr{log}} \migiincl \widetilde{\mcT}_{\mcE_\mbG^\mr{log}/T^\mr{log}}$ induced by the natural inclusion $\mcE_\mbB \migiincl \mcE_\mbG$, then the collection
$\{ \widetilde{\mcT}_{\mcE_\mbG^\mr{log}/T^\mr{log}}^j\}_{j \in  \mbZ}$
forms a decreasing filtration on $\widetilde{\mcT}_{\mcE_\mbG^\mr{log}/T^\mr{log}}$.
Since each $\mfg^{-\alpha}$ ($\alpha \in \Gamma$) is closed under the action of $\mbB$ (defined to be the composite $\mbB \migisurj \mbT \xrightarrow{\mr{adj.\,rep.}} \mr{Aut} (\mfg^{-\alpha})$), we have a decomposition
\begin{align} \label{decom97}
\widetilde{\mcT}^{-1}_{\mcE_\mbG^\mr{log}/T^\mr{log}}/ \widetilde{\mcT}^{0}_{\mcE_\mbG^\mr{log}/T^\mr{log}} \isom \bigoplus_{\alpha \in \Gamma} \mfg^{- \alpha}_{\mcE_\mbB} 
\end{align}
arising from the decomposition $\mfg^{-1}/\mfg^0 = \bigoplus_{\alpha \in \Gamma} \mfg^\alpha$.
We recall from ~\cite{Wak5}, Definition 2.2.1, the definition of a $\mfg$-oper, as follows.

\vspace{3mm}
\bde \label{y035} \leavevmode\\
\vspace{-5mm}
\begin{itemize}
\item[(i)]
 Let 
 \begin{equation} \mcE^\spadesuit := (\pi_\mbB: \mcE_\mbB \migi U, \nabla_\mcE : \mcT_{U^\mr{log}/T^\mr{log}} \migi \widetilde{\mcT}_{\mcE_\mbG^\mr{log}/T^\mr{log}})\end{equation}
  be  a pair consisting of a right $\mbB$-torsor $\mcE_\mbB$ over $U$ and a $T^\mr{log}$-connection $\nabla_\mcE$  
 on the right $\mbG$-torsor $\pi_\mbG : \mcE_\mbG \migi U$  induced by $\mcE_\mbB$.
We shall say that the pair $\mcE^\spadesuit = (\mcE_\mbB, \nabla_\mcE)$ is 
a {\bf  $\mfg$-oper}
on $U^\mr{log}/T^\mr{log}$ if it satisfies the following two conditions:
\begin{itemize}
\item[$\bullet$]
 $\nabla_\mcE(\mcT_{U^\mr{log}/T^\mr{log}}) \subseteq  \widetilde{\mcT}_{\mcE_\mbG^\mr{log}/T^\mr{log}}^{-1}$;
\item[$\bullet$]
For any $\alpha \in \Gamma$, the composite
\begin{equation}
\label{isomoper}
\mcK \mcS_{\mcE^\spadesuit}^{\alpha} :
 \mcT_{U^\mr{log}/T^\mr{log}} \xrightarrow{\nabla_\mcE} \widetilde{\mcT}_{\mcE_\mbG^\mr{log}/T^\mr{log}}^{-1} \migisurj \widetilde{\mcT}_{\mcE^\mr{log}_\mbG/T^\mr{log}}^{-1} /\widetilde{\mcT}^0_{\mcE_\mbG^\mr{log}/T^\mr{log}} \migisurj   \mfg^{-\alpha}_{\mcE_\mbB},
  \end{equation}
is an isomorphism, where the third arrow denotes the natural projection relative  to the decomposition (\ref{decom97}).
\end{itemize}

\item[(ii)]
Let $\mcE^\spadesuit := (\mcE_\mbB, \nabla_\mcE)$, $\mcF^\spadesuit := (\mcF_\mbB, \nabla_\mcF)$ be $\mfg$-opers on $U^\mr{log}/T^\mr{log}$.
An {\bf isomorphism of $\mfg$-opers} 
 from $\mcE^\spadesuit$ to $\mcF^\spadesuit$ is an isomorphism $\mcE_\mbB \isom \mcF_\mbB$ of right $\mbB$-torsors  such that the induced isomorphism $\mcE_\mbG \isom \mcF_\mbG$ of right $\mbG$-torsors is compatible with the respective $T^\mr{log}$-connections $\nabla_\mcE$ and  $\nabla_\mcF$.
\end{itemize}
  \ede

\vspace{5mm}
\subsection{Miura $\mfg$-opers}
\leavevmode\\ \vspace{-4mm}

Next, we shall introduce the definition of a Miura $\mfg$-oper on a log-curve.

\vspace{3mm}
\bde \label{170} \leavevmode\\
 \ \ \ 
 \vspace{-5mm}
\begin{itemize}
\item[(i)]
 A {\bf Miura $\mfg$-oper} on $U^\mr{log}/T^\mr{log}$
   is a collection of data
 \begin{equation} \label{123}
 \widehat{\mcE}^\spadesuit := (\mcE_\mbB, \nabla_\mcE, \mcE'_{\mbB}, \eta_\mcE)
 \end{equation}
consisting of 
\begin{itemize}
\item[$\bullet$]
 a $\mfg$-oper $(\mcE_\mbB, \nabla_\mcE)$  on $U^\mr{log}/T^\mr{log}$;
 \item[$\bullet$]
  a $\mbB$-torsor $\mcE'_{\mbB}$ over $U$, 
   where   we shall write $\mcE'_\mbG := \mcE'_{\mbB} \times^{\mbB} \mbG$;
\item[$\bullet$]
  an isomorphism $\eta_\mcE : \mcE'_\mbG \isom \mcE_\mbG$ of $\mbG$-torsors
such that if $d \eta_\mcE$ denotes the $\mcO_U$-linear isomorphism
$(\widetilde{\mcT}_{{\mcE'}^{\mr{log}}_{\mbB}} \subseteq) \ \widetilde{\mcT}_{{\mcE'}_\mbG^{\mr{log}}} \isom \widetilde{\mcT}_{\mcE_\mbG^{\mr{log}}}$ obtained by differentiating $\eta_\mcE$, then we have 
\begin{align}
\nabla_\mcE (\mcT_{U^\mr{log}/T^\mr{log}}) \subseteq d \eta_\mcE (\widetilde{\mcT}_{{\mcE'}^{ \mr{log}}_{\mbB}}).
\end{align}
\end{itemize}
The $T^\mr{log}$-connection $\nabla_\mcE$ specifies a $T^\mr{log}$-connection
$\nabla_{\mcE'_\mbB} : \mcT_{U^\mr{log}/T^\mr{log}} \migi \widetilde{\mcT}_{{\mcE'}^\mr{log}_\mbB/T^\mr{log}}$ on $\mcE'_\mbB$ in such a way that
the composite 
\begin{align}
\mcT_{U^\mr{log}/T^\mr{log}} \xrightarrow{\nabla_{\mcE'_\mbB}} \widetilde{\mcT}_{{\mcE'}^\mr{log}_\mbB/T^\mr{log}} \migiincl \widetilde{\mcT}_{{\mcE'}^\mr{log}_\mbG/T^\mr{log}}  \stackrel{d \eta_\mcE}{\migi} \widetilde{\mcT}_{\mcE^\mr{log}_\mbG/T^\mr{log}}
\end{align}
coincides with $\nabla_\mcE$.
We shall refer to $(\mcE'_\mbB, \nabla_{\mcE'_\mbB})$ as the {\bf  log flat $\mbB$-torsor associated with $\widehat{\mcE}^\spadesuit$}.
Also, we shall refer to $\mcE^\spadesuit := (\mcE_\mbB, \nabla_\mcE)$ as the {\bf underlying $\mfg$-oper of $\widehat{\mcE}^\spadesuit$}.
If $U^\mr{log}/T^\mr{log} = X^{\mfX \text{-} \mr{log}}/S^{\mfX \text{-} \mr{log}}$ for some pointed stable curve $\mfX := (X/S, \{ \sigma_i \}_{i=1}^r)$,
then we shall refer to any  Miura $\mfg$-oper on $X^{\mfX \text{-} \mr{log}}/S^{\mfX \text{-} \mr{log}}$ as a {\bf Miura $\mfg$-oper on $\mfX$}.

 \item[(ii)]
 Let $\widehat{\mcE}^\spadesuit := (\mcE_\mbB, \nabla_\mcE, \mcE'_{\mbB}, \eta_\mcE)$, 
 $\widehat{\mcF}^\spadesuit:= (\mcF_\mbB,  \nabla_\mcF, \mcF'_{\mbB}, \eta_\mcF)$
  be Miura $\mfg$-opers on $U^\mr{log}/T^\mr{log}$.
 An {\bf isomorphism of Miura $\mfg$-opers} from $\widehat{\mcE}^\spadesuit$ to $\widehat{\mcF}^\spadesuit$ is a pair 
 \begin{equation}
 (\alpha_\mbB, \alpha'_{\mbB}),
 \end{equation}
 consisting of 
  \begin{itemize}
  \item[$\bullet$]
an isomorphism   $\alpha_\mbB : (\mcE_\mbB, \nabla_\mcE) \isom (\mcF_\mbB, \nabla_\mcF)$ of $\mfg$-opers (i.e., an isomorphism $\alpha_\mbB : \mcE_\mbB \isom \mcF_\mbB$ of $\mbB$-torsor respecting the structures of $T^\mr{log}$-connection);
  \item[$\bullet$]
   an isomorphism $\alpha'_{\mbB} : \mcE'_{\mbB} \isom \mcF'_{\mbB}$ of right $\mbB$-torsors  such that the induced isomorphism $\alpha'_{\mbG} : \mcE'_{\mbG} \isom \mcF'_{\mbG}$ of right $\mbG$-torsors 
 satisfies the equality
$\alpha_\mbG \circ \eta_\mcE =  \eta_\mcF \circ \alpha'_{\mbG}$. 
  \end{itemize}
  
   \end{itemize}
 \ede

\vspace{3mm}
\bpr \label{171}\leavevmode\\
 \ \ \ 
Any  Miura $\mfg$-oper on $U^\mr{log}/T^\mr{log}$ does not have nontrivial automorphisms.
 \epr
\begin{proof}
Each automorphism of a Miura $\mfg$-oper determines and is determined by  an automorphism of its underlying $\mfg$-oper.
Hence, 
the assertion follows directly  from
  ~\cite{Wak5}, Proposition 2.2.5.
\end{proof}

\vspace{5mm}
\subsection{Generic Miura $\mfg$-opers}
\leavevmode\\ \vspace{-4mm}

We shall define the notion of a {\it generic} Miura $\mfg$-oper.
Let $\widehat{\mcE}^\spadesuit := (\mcE_\mbB, \nabla_\mcE, \mcE'_\mbB, \eta_\mcE)$ be 
 a Miura $\mfg$-oper  on $U^\mr{log}/T^\mr{log}$.
Recall that the flag variety associated with $\mbG$ is the quotient $\mbG/ \mbB$, which
classifies all Borel subgroups of $\mbG$.
 (Indeed, the  point of $\mbG/\mbB$ represented by $h \in \mbG$ classifies the Borel subgroup $\mcA d_{h}(\mbB)$, where $\mcA d_{h}$ denotes the automorphism of $\mbG$ given by conjugation by $h$.)
The $\mbB$-orbits (with respect to  the left $\mbB$-action on $\mbG$)  in the flag variety $\mbG/\mbB$ are parametrized by the Weyl group $\mbW$.
That is to say, we have a decomposition $\mbG = \coprod_{w \in \mbW} \mbB  w  \mbB$, which is known as the {\it Bruhat decomposition}.
(More precisely, for each $w \in \mbW = N_\mbG (\mbT)/\mbT$ we choose a representative $\dot{w} \in N_\mbG (\mbT)$ of $w$ and consider the double coset $\mbB  \dot{w} \mbB$.
Since $\mbB  \dot{w}  \mbB$ is independent of the choice of the  representative $\dot{w}$ of $w$, we simply denote it by $\mbB  w  \mbB$.)
Let $w_0$ be the longest element of $\mbW$.
The orbit $\mbB  w_0 \mbB \subseteq \mbG/\mbB$
is a   dense open subscheme of $\mbG/\mbB$ (called the {\it big cell}).
The morphism
$\mbN \migi \mbG/\mbB$ given by assigning
 $h \mapsto h  \dot{w}_0  \mbB$ defines an isomorphism $\mbN \isom \mbB  w_0 \mbB$ of $k$-schemes.
By passing to this isomorphism, we shall regard $\mbN$  as a dense open subscheme of $\mbG /\mbB$ that is closed under the left $\mbB$-action on $\mbG/\mbB$.
By twisting $\mbG/\mbB$ by the right $\mbB$-torsor $\mcE_\mbB$, we obtain a proper $U$-scheme 
\begin{equation}
\mcE_\mbB \times^\mbB (\mbG/\mbB),
\end{equation}
 which contains 
 a dense open subscheme $\mcE_\mbB \times^\mbB \mbN$.
It follows from the definition of $\mcE_\mbB \times^\mbB (\mbG/\mbB)$ that
the image of $\mcE'_{\mbB}$ via the isomorphism $\eta_\mcE : \mcE'_\mbG \isom \mcE_\mbG$ determines its classifying morphism
\begin{equation} \label{eq23}
\sigma_{\widehat{\mcE}^\spadesuit} : U \migi \mcE_\mbB \times^\mbB (\mbG/\mbB).
\end{equation}

\vspace{3mm}
\bde \label{GGlk} \leavevmode\\
 \ \ \ 
We shall say that  a Miura  $\mfg$-oper $\widehat{\mcE}^\spadesuit$
 is {\bf generic}
if the image of the morphism 
$\sigma_{\widehat{\mcE}^\spadesuit}$ lies in  $\mcE_\mbB \times^\mbB \mbN$ ($\subseteq \mcE_\mbB \times^\mbB (\mbG/\mbB)$). 
 \ede

\vspace{5mm}
\subsection{Special Miura $\mfg$-opers} \label{SeSe1}
\leavevmode\\
\vspace{-4mm}

Next, we define the notion of a {\it special} Miura $\mfg$-oper.
Write
\begin{align} \label{f015}
\mcE^\dagger_{\mbB, U^\mr{log}/T^\mr{log}} := \mcE^\dagger_{\mbT, U^\mr{log}/T^\mr{log}} \times^\mbT \mbB, \hspace{5mm}
\mcE^\dagger_{\mbG, U^\mr{log}/T^\mr{log}} := \mcE^\dagger_{\mbT, U^\mr{log}/T^\mr{log}} \times^\mbT \mbG.
\end{align}
(cf. (\ref{f016}) for the definition of $\mcE^\dagger_{\mbT, U^\mr{log}/T^\mr{log}}$).
The subscheme $w_0 \mbB$ of $\mbG$
is closed under the left action by $\mbT$.
Hence, we have a 
$\mbB$-torsor 
\begin{align} \label{f014}
{\mcE'}^\dagger_{\mbB, U^\mr{log}/T^\mr{log}} := \mcE^\dagger_{\mbT, U^\mr{log}/T^\mr{log}} \times^\mbT w_0 \mbB,
\end{align}
which admits  a natural  isomorphism
 \begin{align} \label{f013}
 \eta^\dagger_\mcE : {\mcE'}^\dagger_{\mbB, U^\mr{log}/T^\mr{log}} \times^\mbB \mbG \isom \mcE^\dagger_{\mbG, U^\mr{log}/T^\mr{log}}
 \end{align}
of $\mbG$-torsors.
If there is no fear of causing confusion, then we shall write
\begin{align}
\mcE^\dagger_\mbB := \mcE^\dagger_{\mbB, U^\mr{log}/T^\mr{log}},   \ \ \ \mcE^\dagger_\mbG := \mcE^\dagger_{\mbG, U^\mr{log}/T^\mr{log}}, \ \ \  \ {\mcE'}^\dagger_\mbB := {\mcE'}^\dagger_{\mbB, U^\mr{log}/T^\mr{log}}
\end{align}
 for simplicity.

For each $\alpha \in \Gamma$,  let us fix a generator $x_\alpha$ of $\mfg^{\alpha}$.
Write $p_1 := \sum_{\alpha \in \Gamma} x_\alpha$.
Then, one may find 
a unique collection $(y_\alpha)_{\alpha \in \Gamma}$, where each $y_\alpha$ is a generator of $\mfg^{-\alpha}$, such that if we write
\begin{equation} \label{p_{-1}}
p_{-1}  := \sum_{\alpha \in \Gamma} y_\alpha \in \mfg_{-1},
\end{equation} 
then the set 
 $\{  p_{-1}, 2 \check{\rho}, p_1  \}$ forms an $\mfs \mfl_2$-triple (cf. \S\,\ref{SSQ1} for the definition of $\check{\rho}$).
For each $\alpha \in \Gamma$, we shall write
\begin{align} \label{f012}
\eta^\alpha : \mfg_{\mcE^{\dagger \mr{log}}_{\mbB}}^{-\alpha} \ \left(\isom  \mfg_{\mcE^{\dagger \mr{log}}_{\mbT}}^{-\alpha}\right) \isom \mcT_{U^\mr{log}/T^\mr{log}}
\end{align}
for the isomorphism determined uniquely by the following condition: 
for each local trivialization  $\tau :  \mcO_V \isom \mcT_{U^\mr{log}/T^\mr{log}}|_V$ of $\mcT_{U^\mr{log}/T^\mr{log}}$ (where $V$ denotes  an open subscheme of $U$), the composite isomorphism   
\begin{align}
 \mfg_{\mcE^{\dagger \mr{log}}_{\mbB}}^{-\alpha} |_V
&  \isom   \mfg_{\mcE^{\dagger \mr{log}}_{\mbT}}^{-\alpha}|_V 
 \isom   (k \cdot y_\alpha)_{(\Omega_{U^\mr{log}/T^\mr{log}}|_V)^{\times} \times^{\mbG_m, \check{\rho}}\mbT} \\ 
 & \isom   (k \cdot y_\alpha)_{\mcO_V^\times \times^{\mbG_m, \check{\rho}} \mbT}  \isom  k_{V \times_k \mbT} 
  \isom  \mcO_V 
   \stackrel{\tau}{\isom}  \mcT_{U^\mr{log}/T^\mr{log}} |_V \notag
\end{align}
coincides with $\eta^\alpha |_V$,  where the third isomorphism arises from the dual isomorphism $\Omega_{U^\mr{log}/T^\mr{log}} |_V \isom \mcO_V$ of $\tau$ and the fourth isomorphism is given by $a \cdot y_\alpha \mapsto a$ (for any $a \in k$).

\begin{rema} \label{RRRH}
\leavevmode\\
 \ \ \ 
In the case where  $\mbG = \mr{PGL}_n$ (hence $\mfg = \mfp \mfg \mfl_n = \mfs \mfl_n$),  we fix 
an   $\mfs \mfl_n$-triple  $\{  p_{-1}, 2 \check{\rho}, p_1  \}$  given by 
\begin{align}
& p_{-1} := \begin{pmatrix} 0 & 0 & \cdots  & 0 & 0 \\  1 & 0 & \cdots  & 0 & 0 \\  0 & 1 & \cdots  & 0 & 0 \\ \vdots  & \vdots  & \ddots  & \vdots  & \vdots  \\  0 & 0 & \cdots  & 1 & 0    \end{pmatrix}, \ \ 
2 \check{\rho} := \begin{pmatrix} n-1 & 0 & 0 & \cdots   & 0 \\  0 & n-3 & 0 &  \cdots   & 0 \\  0 & 0 & n-5 & \cdots   & 0 \\ \vdots  & \vdots  & \vdots  & \ddots  & \vdots  \\  0 & 0 & 0 &  \cdots   & - (n-1)
\end{pmatrix}, \\
 & \hspace{25mm} p_1 := \begin{pmatrix} 0 & n-1 & 0 & 0 & \cdots   & 0 \\  0 & 0 & 2(n-2) & 0 &  \cdots   & 0 \\  0 & 0 & 0  & 3(n-2) & \cdots   & 0 \\ \vdots  & \vdots  & \vdots  & \vdots &  \ddots  & \vdots  \\  0 & 0 & 0 & 0&  \cdots   & n-1   \\ 0 & 0 & 0 & 0&  \cdots   &0    \end{pmatrix}. \notag
\end{align} 
\end{rema}

\vspace{3mm}
\bde \leavevmode\\
 \ \ \ 
A {\bf $p_{-1}$-special Miura $\mfg$-oper} on $U^\mr{log}/T^\mr{log}$  is a Miura $\mfg$-oper on $U^\mr{log}/T^\mr{log}$ of the form 
\begin{align}
\widehat{\mcE}^{\spadesuit \diamondsuit} :=  (\mcE^\dagger_{\mbB, U^\mr{log}/T^\mr{log}}, \nabla_\mcE, {\mcE'}^\dagger_{\mbB, U^\mr{log}/T^\mr{log}}, \eta^\dagger_\mcE)
\end{align}
 (for some $T^\mr{log}$-connection $\nabla_\mcE$ on ${\mcE}^\dagger_{\mbG, U^\mr{log}/T^\mr{log}}$) satisfying the equality
$\eta^\alpha \circ \mcK \mcS^\alpha_{\mcE^{\spadesuit \diamondsuit}} = \mr{id}_{\mcT_{U^\mr{log}/T^\mr{log}}} $ for any  $\alpha \in \Gamma$, where $\mcE^{\spadesuit \diamondsuit} := (\mcE^\dagger_{\mbB, U^\mr{log}/T^\mr{log}}, \nabla_\mcE)$.
(In particular, any $p_{-1}$-special Miura $\mfg$-oper is a generic Miura $\mfg$-oper.) 
 \ede

\vspace{3mm}
\bpr \label{P171}\leavevmode\\
 \ \ \ 
For any generic Miura $\mfg$-oper $\widehat{\mcE}^\spadesuit := (\mcE_\mbB, \nabla_\mcE, \mcE'_\mbB, \eta_\mcE)$ on $U^\mr{log}/T^\mr{log}$,
there exists  a unique pair $(\widehat{\mcE}^{\spadesuit \diamondsuit}, \iota_{\widehat{\mcE}^\spadesuit})$ consisting of a $p_{-1}$-special Miura $\mfg$-oper  on $U^\mr{log}/T^\mr{log}$ and an isomorphism $\iota_{\widehat{\mcE}^\spadesuit} : \widehat{\mcE}^\spadesuit \isom \widehat{\mcE}^{\spadesuit \diamondsuit}$ of Miura $\mfg$-opers.
 \epr
\begin{proof}
Since Miura $\mfg$-opers may be constructed by means of descent with respect to \'{e}tale morphisms,
we are (by taking account of Proposition \ref{171}) free to replace $U$ with its \'{e}tale covering.
Thus, we may 
assume that there exists a global section  
$\partial  \in \Gamma (U, \Omega_{U^\mr{log}/T^\mr{log}})$ with $\mcO_U \cdot \partial = \mcT_{U^\mr{log}/T^\mr{log}}$.
Also, we may assume that
$\mcE_\mbB = U \times_k \mbB$ and  $\mcE'_\mbB = U \times_k w_0 \mbB$.

First, let us consider the existence assertion.
It follows from ~\cite{Wak5}, Lemma 2.2.4, that
there exists 
a unique $U$-rational point $h : U \migi \mbT$ of $\mbT$ such that
the underlying $\mfg$-oper $\mcE^\spadesuit$ of $\widehat{\mcE}^\spadesuit$ is {\it of  precanonical type relative to the triple $(U, x, \mfl_{h})$} (cf. ~\cite{Wak5}, Definition 2.2.3), where $\mfl_h$ denotes the left translation $U \times_k \mbG \isom U \times_k \mbG$ determined by $h : U \migi \mbT$ ($\subseteq \mbG$). 
Let us write
\begin{align}
\widehat{\mcE}^{\spadesuit \diamondsuit} := (U \times_k \mbB, \mfl_h^* (\nabla_\mcE), U \times_k w_0 \mbB, \mr{id}_{U \times_k \mbG}),
\end{align}
where $\mfl_h^* (\nabla_\mcE)$ denotes the $T^\mr{log}$-connection on $U \times_k \mbG$  obtained  from $\nabla_\mcE$ via  pull-back  by  $\mfl_h$.
Then,  the pair $(\widehat{\mcE}^{\spadesuit \diamondsuit}, \iota_{\widehat{\mcE}^\spadesuit})$ specifies a desired pair, where $\iota_{\widehat{\mcE}^\spadesuit} := (\mfl_h |_{U \times_k \mbB}, \mfl_h |_{U \times_k w_0\mbB})$.
This completes the existence assertion.

Next, let us prove the uniqueness assertion.
To this end, it suffices to prove the claim that if 
$\widehat{\mcE}_1^{\spadesuit \diamondsuit}$ and  $\widehat{\mcE}_2^{\spadesuit \diamondsuit}$
are $p_{-1}$-special Miura $\mfg$-opers on $U^\mr{log}/T^\mr{log}$  admitting an isomorphism $(\alpha_\mbB, \alpha'_\mbB) : \widehat{\mcE}_1^{\spadesuit \diamondsuit} \migi \widehat{\mcE}_2^{\spadesuit \diamondsuit}$, then $\widehat{\mcE}_1^{\spadesuit \diamondsuit}$ is identical to  $\widehat{\mcE}_2^{\spadesuit \diamondsuit}$.
Since  the automorphisms $\alpha_\mbB : U \times_k \mbB \isom U \times_k \mbB$
 and $\alpha'_\mbB : U \times_k w_0 \mbB \isom U \times_k w_0 \mbB$
  are compatible in the evident  sense,
they come from an automorphism of $U \times_k \mbT$.
Strictly speaking, there exists an automorphism $\alpha_\mbT$
of $U \times_k \mbT$ (i.e., a left-translation on $U \times_k \mbT$ by some $U \migi \mbT$) which induces both $\alpha_\mbB$ and $\alpha'_\mbB$.
But, because of  the equalities $\eta^\alpha \circ \mcK \mcS^\alpha_{\mcE_1^{\spadesuit \diamondsuit}} = \mr{id}_{\mcT_{U^\mr{log}/T^\mr{log}}}$ and $\eta^\alpha \circ \mcK \mcS^\alpha_{\mcE_2^{\spadesuit \diamondsuit}} = \mr{id}_{\mcT_{U^\mr{log}/T^\mr{log}}}$ (for every $\alpha \in \Gamma$), $\alpha_\mbT$ must be the identity morphism (cf. the proof of ~\cite{Wak5}, Lemma 2.2.4).
Hence, we have  $\widehat{\mcE}_1^{\spadesuit \diamondsuit} = \widehat{\mcE}_2^{\spadesuit \diamondsuit}$.
This completes the proof of the uniqueness assertion, and hence, Proposition \ref{P171}. 
\end{proof}

\vspace{3mm}
\bde \leavevmode\\
 \ \ \ 
For each generic Miura $\mfg$-oper $\widehat{\mcE}^\spadesuit$ on $U^\mr{log}/T^\mr{log}$,  we shall refer to the pair $(\widehat{\mcE}^{\spadesuit \diamondsuit}, \iota_{\widehat{\mcE}^\spadesuit})$ obtained by applying  Proposition \ref{P171} to $\widehat{\mcE}^\spadesuit$
as the {\bf $p_{-1}$-specialization} of $\widehat{\mcE}^\spadesuit$.
 \ede

\vspace{5mm}
\subsection{From  $\mfg$-Cartan  connections  to special Miura $\mfg$-opers}
\leavevmode\\ \vspace{-4mm}

In what follows, let us 
construct a bijective correspondence between the $\mfg$-Cartan  connections and the generic Miura $\mfg$-opers (cf. Proposition \ref{PP05}).
To begin with, we shall write 
\begin{align} \label{EE00}
\iota : \widetilde{\mcT}_{\mcE^{\dagger \mr{log}}_{\mbT}/T^\mr{log}} \oplus \bigoplus_{\alpha \in \Gamma}\mfg^{-\alpha}_{\mcE^\dagger_{\mbB}}\migiincl \widetilde{\mcT}_{\mcE^{\dagger \mr{log}}_{\mbG}/T^\mr{log}}
\end{align}
for the $\mcO_U$-linear morphism determined by the natural inclusions $\widetilde{\mcT}_{\mcE^{\dagger \mr{log}}_{\mbT}/T^\mr{log}} \migiincl \widetilde{\mcT}_{\mcE^{\dagger \mr{log}}_{\mbG}/T^\mr{log}}$ (obtained by differentiating   the inclusion $\mcE_{\mbT}^\dagger \migi \mcE^\dagger_{\mbG}$) and $\mfg^{-\alpha}_{\mcE^\dagger_{\mbB}} \migiincl \widetilde{\mcT}_{\mcE^{\dagger \mr{log}}_{\mbG}/T^\mr{log}}$ (for $\alpha \in \Gamma$).

Now, let $\nabla : \mcT_{U^\mr{log}/T^\mr{log}} \migi \widetilde{\mcT}_{\mcE^{\dagger \mr{log}}_{\mbT}/T^\mr{log}}$ be  a 
$\mfg$-Cartan connection on $U^\mr{log}/T^\mr{log}$.
The composite
\begin{align}
\nabla_\mcE := \iota \circ (\nabla \oplus  \bigoplus_{\alpha \in \Gamma} \eta^{\alpha}) : \mcT_{U^\mr{log}/T^\mr{log}} \migi \widetilde{\mcT}_{\mcE^{\dagger \mr{log}}_{\mbG}/T^\mr{log}}
\end{align}
 specifies  a $T^\mr{log}$-connection, and   the quadruple 
 \begin{align}
 \widehat{\mcE}^{\spadesuit \diamondsuit}_{\nabla} := (\mcE^\dagger_\mbB, \nabla_\mcE, {\mcE'}^\dagger_\mbB, \eta^\dagger_\mcE)
 \end{align}
  forms  a $p_{-1}$-special Miura $\mfg$-oper on $U^\mr{log}/T^\mr{log}$.

Conversely, let us take 
a generic Miura $\mfg$-oper $\widehat{\mcE}^\spadesuit$ on $U^\mr{log}/T^\mr{log}$.
Denote by  $(\widehat{\mcE}^{\spadesuit \diamondsuit}, \iota_{\widehat{\mcE}^\spadesuit})$ the $p_{-1}$-specialization of $\widehat{\mcE}^\spadesuit$.
Also, denote by $({\mcE'}^\dagger_\mbB, \nabla_{{\mcE'}^\dagger_\mbB})$ the log flat $\mbB$-torsor associated with $\widehat{\mcE}^{\spadesuit \diamondsuit}$.
Then,  the composite 
\begin{align}
\nabla_{\widehat{\mcE}^{\spadesuit}} : \mcT_{U^\mr{log}/T^\mr{log}} \xrightarrow{\nabla_{{\mcE'}^\dagger_\mbB}} \widetilde{\mcT}_{{\mcE'}^{\dagger \mr{log}}_{\mbB}/T^\mr{log}} \migisurj \widetilde{\mcT}_{\mcE_\mbT^{\dagger \mr{log}}/T^\mr{log}}
\end{align}
forms a  $\mfg$-Cartan connection on $U^\mr{log}/T^\mr{log}$,
where the second arrow denotes the surjection arising from  the quotient $\mbB \migisurj \mbT$. 
The $\mfg$-Cartan connection $\nabla_{\widehat{\mcE}^{\spadesuit}}$ depends only on the isomorphism class of $\widehat{\mcE}^\spadesuit$. 
Moreover, the following proposition holds.

\vspace{3mm}
\bpr \label{PP05}
 \leavevmode\\
 \ \ \ 
 The assignments $\nabla \mapsto \widehat{\mcE}^{\spadesuit \diamondsuit}_\nabla$ and $\widehat{\mcE}^{\spadesuit} \mapsto \nabla_{\widehat{\mcE}^{\spadesuit}}$  discussed above define the following bijection:
 \begin{equation} \label{equicat}
\begin{pmatrix}
\text{the set of $\mfg$-Cartan} \\
\text{connections  on $U^\mr{log}/T^\mr{log}$}
\end{pmatrix}
\isom
\begin{pmatrix}
\text{the set  of isomorphism classes of  } \\
\text{generic Miura $\mfg$-opers on  $U^\mr{log}/T^\mr{log}$}
\end{pmatrix}.
\end{equation}
Moreover, this bijection is 
functorial with respect to
  $S$.
 \epr
\begin{proof}
The assertion follows immediately from the definitions of the assignments  involved.
\end{proof}

\vspace{5mm}
\subsection{Miura $\mfg$-opers of prescribed exponents}
\leavevmode\\ \vspace{-4mm}

Let $S$ be a $k$-scheme and $\mfX := (X/S, \{ \sigma_i \}_{i=1}^r)$ an $r$-pointed stable curve over $S$ of genus $g$.
Denote by $\mbB^-$ ($\subseteq \mbG$) 
the opposite Borel subgroup of $\mbB$ relative to $\mbT$, which is, by definition, a unique Borel subgroup satisfying that 
 $\mbB \cap \mbB^- = \mbT$.
Now, suppose that $r >0$, and let us fix $i \in \{ 1, \cdots, r \}$. 
By 
 change of structure group, 
the composite isomorphism (\ref{eqeq39}) induces an isomorphism $\sigma^* (\mcE'^\dagger_\mbB) \isom S \times_k w_0 \mbB$, and hence, induces a sequence of  isomorphisms
\begin{align} \label{eqeq51}
\Gamma (S, \sigma^*_i (\mfb_{\mcE'^\dagger_\mbB})) \isom \Gamma (S, \mfb_{\sigma^* (\mcE'^\dagger_\mbB)}) \isom \Gamma (S, \mfb_{S \times_k w_0 \mbB}) \isom \mfb^{-} (S),
\end{align}
where $\mfb^{-}$ denotes the Lie algebra of $\mbB^{-}$.
In particular, we obtain a surjection
\begin{align} \label{eqeq50}
\Gamma (S, \sigma^*_i (\mfb_{\mcE'^\dagger_\mbB}))  \xrightarrow{(\ref{eqeq51})} \mfb^{-} (S) \migisurj \mft (S).
\end{align}

Now, let $\widehat{\mcE}^\spadesuit := (\mcE_\mbB, \nabla_\mcE, \mcE'_\mbB, \eta_\mcE)$ be a generic  Miura $\mfg$-oper on $\mfX$.
Denote by $(\widehat{\mcE}^{\spadesuit \diamondsuit}, \iota_{\widehat{\mcE}^\spadesuit})$
the $p_{-1}$-specialization of $\widehat{\mcE}^\spadesuit$.
 The  monodromy 
  of the flat $\mbB$-torsor associated with $\widehat{\mcE}^{\spadesuit \diamondsuit}$ at each marked point $\sigma_i$ ($i =1, \cdots, r$)
is sent, 
 via the composite (\ref{eqeq50}), 
 to an element
\begin{align} \label{f01}
 \varepsilon^{\widehat{\mcE}^\spadesuit}_i
 \in \mft (S).
\end{align}

\vspace{3mm}
\bde \leavevmode\\
\vspace{-5mm}
\begin{itemize}
\item[(i)]
For each  $i \in \{ 1, \cdots, r\}$,  we shall refer to $\varepsilon^{\widehat{\mcE}^\spadesuit}_i$ as the {\bf exponent} of $\widehat{\mcE}^\spadesuit$ at $\sigma_i$.
\item[(ii)]
Let $\vec{\varepsilon} := (\varepsilon_i)_{i=1}^r \in \mft (S)^{\times r}$, and let $\widehat{\mcE}^\spadesuit$ be  a generic Miura $\mfg$-oper  on $\mfX$.
We shall say that  $\widehat{\mcE}^\spadesuit$  is {\bf of exponents $\vec{\varepsilon}$} if $\varepsilon^{\widehat{\mcE}^\spadesuit}_i  = \varepsilon_i$ for any $i \in \{ 1, \cdots, r \}$.
If $r =0$, then we shall refer to any generic Miura $\mfg$-oper as being {\bf of exponent  $\emptyset$}.
\end{itemize}
\ede

\vspace{3mm}
\bpr \label{pr051}
 \leavevmode\\
 \ \ \ 
Let $\vec{\varepsilon}$ be an element of $\mft (S)^{\times r}$ (where  $\vec{\varepsilon} := \emptyset$ if $r =0$).
Then, the  bijection (\ref{equicat})  (of the case where $U^\mr{log}/T^\mr{log}$ is taken to be $X^{\mfX \text{-}\mr{log}}/S^{\mfX \text{-}\mr{log}}$) restricts to a functrial (with respect to 
$S$) bijection:
 \begin{equation} \label{equicat2}
\begin{pmatrix}
\text{the set of $\mfg$-Cartan connections} \\
\text{on $\mfX$ of monodromies $\vec{\varepsilon}$} 
\end{pmatrix}
\isom
\begin{pmatrix}
\text{the set  of isomorphism classes} \\
\text{of generic Miura $\mfg$-opers} \\
\text{on  $\mfX$ of exponents $\vec{\varepsilon}$}
\end{pmatrix}.
\end{equation}
 \epr
\vspace{3mm}

\begin{rema} \label{f010}
\leavevmode\\
 \ \ \ 
 Let $\widehat{\mcE}^\spadesuit$ be a generic Miura $\mfg$-oper on $\mfX$ of exponents $\vec{\varepsilon} : = (\varepsilon_i)_{i =1}^r \in \mft(S)^{\times r}$, where
each $\varepsilon_i$ ($i \in \{1, \cdots, r \}$) is supposed to be  regular.
Then, one verifies that the  radii (cf. ~\cite{Wak5}, Definition 2.9.2 (i)) of   the underlying $\mfg$-oper $\mcE^\spadesuit$ of $\widehat{\mcE}^\spadesuit$ is given by
\begin{align} \label{eqeq18}
\chi (\vec{\varepsilon}) := (\chi (\varepsilon_i))_{i=1}^r \in \mfc (S)^{\times r}.
\end{align}
  \end{rema}

\vspace{5mm}
\subsection{Moduli of generic Miura $\mfg$-opers} \label{sc333}
\leavevmode\\ \vspace{-4mm}

We shall write
\begin{align} \label{eqeq100}
\mfM \overline{\mfO} \mfp_{\mfg, g,r}
\end{align}
for the set-valued contravariant functor on $\mfS \mfc \mfh_{/\overline{\mfM}_{g,r}}$ which,
to any  object $S \migi \overline{\mfM}_{g,r}$ of $\mfS \mfc \mfh_{/\overline{\mfM}_{g,r}}$ classifying a pointed stable curve $\mfX$, assigns the set of isomorphism classes of generic Miura $\mfg$-opers on $\mfX$.
Also, for each $\vec{\varepsilon} \in \mft(k)^{\times r}$ (where $\vec{\varepsilon} := \emptyset$ if $r =0$), we shall write
\begin{align} \label{ee200}
\mfM \overline{\mfO} \mfp_{\mfg, g,r, \vec{\varepsilon}}
\end{align}
for the subfunctor of $\mfM \overline{\mfO} \mfp_{\mfg, g,r}$ classifying generic Miura $\mfg$-opers of exponents $\vec{\varepsilon}$.
Then, the following proposition holds.

\vspace{3mm}
\bpr \label{eqeq88} \leavevmode\\
 \ \ \ 
The functorial bijection (\ref{equicat}) determines   an  isomorphism
\begin{align} \label{eqd30}
\Xi_{\mfg, g,r, p_{-1}} : \mfC \overline{\mfC} \mfo_{\mfg, g,r} \isom \mfM \overline{\mfO} \mfp_{\mfg, g,r}
\end{align}
over $\overline{\mfM}_{g,r}$.
Moreover, if $\vec{\varepsilon} \in \mft (k)^{\times r}$, then  the isomorphisms  (\ref{eqd30}) restricts to an  isomorphism
\begin{align} \label{eq30}
\Xi_{\mfg, g,r, \vec{\varepsilon}, p_{-1}}  :\mfC  \overline{\mfC} \mfo_{\mfg, g,r, \vec{\varepsilon}} \isom \mfM \overline{\mfO} \mfp_{\mfg, g,r, \vec{\varepsilon}}. 
\end{align}
In particular, both $\mfM \overline{\mfO} \mfp_{\mfg, g,r} $ and $\mfM \overline{\mfO} \mfp_{\mfg, g,r, \vec{\varepsilon}} $ (for any $\vec{\varepsilon}$)  may be represented by 
 Deligne-Mumford stacks over $k$.
 \epr
\begin{proof}
The last assertion follows from Proposition \ref{P174} (i).
\end{proof}
\vspace{3mm}

\begin{rema} \label{f039}
\leavevmode\\
 \ \ \ 
If $\mr{char} (k) =0$, then we have
 \begin{align} \label{ee3001}
 \mfM \overline{\mfO} \mfp_{\mfs \mfl_2, g,0} \times_{\overline{\mfM}_{g,0}} \mfM_{g, 0} \  (\stackrel{\Xi_{\mfs \mfl_2, g,0, \emptyset, p_{-1}}}{\cong} \mfC \mfC \mfo_{\mfs \mfl_2, g, 0} \stackrel{(\ref{eq070})}{\cong} \mfC \mfo_{g, 0}) = \emptyset.
 \end{align} 
Indeed, it is well-known (cf. ~\cite{HTT}, Remark 1.2.10) that if $k$ is an algebraically closed field of characteristic zero and $X$ is a proper smooth curve over $k$  of genus $g >1$, then
 the degree of  any line bundle on $X$ admitting  a $k$-connection must be    $0$.
On the other hand, $\Omega_{X/k}$ has positive degree (more precisely,  $\mr{deg} (\Omega_{X/k}) = 2g -2 >0$).
It follows that there is no generic Miura $\mfs \mfl_2$-oper  on $X$.
That is to say, the stack  $\mfM \mfO \mfp_{\mfs \mfl_2, g,0}$
  turns out to be empty.
  \end{rema}

\vspace{5mm}
\subsection{Dormant and $p$-nilpotent  Miura opers} \label{sc334}
\leavevmode\\ \vspace{-4mm}

In this subsection, we consider Miura $\mfg$-opers {\it in positive characteristic}.
Assume  
that   either $(\text{Char})_{p}$ or $(\text{Char})_{p}^{\mfs \mfl}$ is satisfied.

Let $\widehat{\mcE}^\spadesuit$ be a generic Miura $\mfg$-oper on a pointed stable curve $\mfX$, 
$(\widehat{\mcE}^{\spadesuit \diamondsuit}, \iota_{\widehat{\mcE}^\spadesuit})$ 
 the specialization of $\widehat{\mcE}^\spadesuit$, and 
$(\mcE'_\mbB, \nabla_{\mcE'_\mbB})$ the log flat $\mbB$-torsor associated with $\widehat{\mcE}^{\spadesuit \diamondsuit}$.
Then, the assignment $\widehat{\mcE}^\spadesuit \mapsto \psi^{(\mcE'^\dagger_\mbB, \nabla_{\mcE'_\mbB})}$ determines a morphism
\begin{equation}  \label{dormant2}  \rotatebox{50}{{\large $\kappa$}}_{\mfb, g,r}  : \mfM \overline{\mfO} \mfp_{\mfg, g,r}  \migi  \ooalign{$\bigoplus$ \cr $\bigotimes_{\mfb, g,r}$} 
\end{equation}
over $\overline{\mfM}_{g,r}$.
We shall write
\begin{equation}  \label{dormant8} 
\rotatebox{50}{{\large $\kappa$}}^{\mr{H} \text{-} \mr{M}}_{\mfb, g,r} :  \mfM \overline{\mfO} \mfp_{\mfg, g,r} \migi  \ooalign{$\bigoplus$ \cr $\bigotimes_{\mfc,  g,r}$} 
  \end{equation}
for the composite $\ooalign{$\oplus$ \cr $\otimes_{\mfb \migi \mfc, g,r}$} \circ \rotatebox{50}{{\large $\kappa$}}_{\mfb, g,r}$ (cf. (\ref{spindle2})).

Denote by 
\begin{align}
\mfM \overline{\mfO} \mfp^{^\mr{Zzz...}}_{\mfg, g,r} \ \  (\text{resp.,} \ \mfM \overline{\mfO} \mfp^{^{p \text{-}\mr{nilp}}}_{\mfg,  g,r} )
\end{align}
the closed substack of $\mfM \overline{\mfO} \mfp_{\mfg, g,r}$ defined to be  the inverse image of the zero section  $[0]_{\mfg, g,r}$
  (resp., $[0]_{\mfc, g,r}$)
 (cf. (\ref{EQEQ1})) via the morphism  $ \rotatebox{50}{{\large $\kappa$}}_{\mfb, g,r} $ (resp., $\rotatebox{50}{{\large $\kappa$}}^{\mr{H} \text{-} \mr{M}}_{\mfb, g,r}$).
Also, for each $\vec{\varepsilon} \in \mft(k)^{\times r}$ (where $\vec{\varepsilon} := \emptyset$ if $r =0$), we shall write
\begin{align}\label{eq35gg}
& \hspace{14mm}\mfM \overline{\mfO} \mfp^{^\mr{Zzz...}}_{\mfg, g,r, \vec{\varepsilon}} :=\mfM \overline{\mfO} \mfp^{^\mr{Zzz...}}_{\mfg, g,r}  \times_{\mfM \overline{\mfO} \mfp_{\mfg, g,r} } \mfM \overline{\mfO} \mfp_{\mfg, g,r, \vec{\varepsilon}} \\
&   \left(\text{resp.,} \ \mfM \overline{\mfO} \mfp^{^{p \text{-}\mr{nilp}}}_{\mfg, g,r, \vec{\varepsilon}} :=\mfM \overline{\mfO} \mfp^{^{p \text{-}\mr{nilp}}}_{\mfg, g,r} \times_{\mfM \overline{\mfO} \mfp_{\mfg, g,r} } \mfM \overline{\mfO} \mfp_{\mfg, g,r, \vec{\varepsilon}}\right). \notag 
\end{align}

\vspace{3mm}
\bde \label{Def46} \leavevmode\\
 \ \ \ 
Let $S$ be a $k$-scheme and $\mfX := (X/S, \{ \sigma_i \}_{i=1}^r)$ an $r$-pointed stable curve over $S$ of genus $g$.
We shall say that a generic  Miura $\mfg$-oper on $\mfX$ is {\bf dormant} (resp., {\bf $p$-nilpotent})
if it is classified by $\mfM \overline{\mfO} \mfp^{^\mr{Zzz...}}_{\mfg, g,r}$  (resp.,  $\mfM \overline{\mfO} \mfp^{^{p \text{-}\mr{nilp}}}_{\mfg,  g,r}$).
In other words, a dormant  (resp., $p$-nilpotent) generic Miura $\mfg$-oper is a  generic Miura $\mfg$-oper
whose  underlying $\mfg$-oper  is dormant (resp., $p$-nilpotent) (cf. ~\cite{Wak5}, Definition 3.6.1 and Definition 3.8.3).
\ede

\begin{rema} \label{fffg039}
\leavevmode\\
 \ \ \ 
 In ~\cite{Wak5} , \S\,3.12, 
the moduli stack
\begin{align} \label{eq35}
\overline{\mfO} \mfp_{\mfg, g,r,  \vec{\rho}}^{} \hspace{3mm}(\text{resp.,} \ \overline{\mfO} \mfp_{\mfg, g,r,  \vec{\rho}}^{^\mr{Zzz...}} ;\text{resp.,} \ \overline{\mfO} \mfp_{\mfg, g,r, \vec{\rho}}^{^{p \text{-} \mr{nilp}}} ),
\end{align}
 where $\vec{\rho} \in \mfc (k)^{\times r}$,
 classifying $r$-pointed genus $g$ stable curves over $k$ together with a  $\mfg$-oper  (resp.,  a dormant   $\mfg$-oper; resp., a $p$-nilpotent  $\mfg$-oper) of radii $\vec{\rho}$ was introduced.
(Notice that in {\it loc.\,cit.}, we used the notations $\mfO \mfp_{\mfg,  \vec{\rho},  g,r}$,  $\mfO \mfp_{\mfg, \vec{\rho}, g,r}^{^\mr{Zzz...}}$, and   $\mfO \mfp_{\mfg, \vec{\rho}, g,r}^{^{p \text{-} \mr{nilp}}}$ to denote these  moduli stacks respectively.
We refer to ~\cite{Wak5} for the study  of these moduli stacks.

Let $\vec{\varepsilon} :=(\varepsilon_i)_{i=1}^r \in \mft (k)^{\times r}$ (where $\vec{\varepsilon} := \emptyset$ if $r =0$), and 
suppose that each $\varepsilon_i$ ($i \in \{1, \cdots, r\}$) is regular.
According to the discussion in Remark \ref{f010}, the assignment from each generic Miura  $\mfg$-oper to its underlying $\mfg$-oper determines a morphism
\begin{align} 
\text{み}_{\mfg, g,r, \vec{\varepsilon}} : \left(\mfC \overline{\mfC} \mfo_{\mfg, g,r, \vec{\varepsilon}} \xrightarrow[\sim]{\Xi_{\mfg, g,r,  \vec{\varepsilon}, p_{-1}}} \right)  \ \mfM \overline{\mfO} \mfp_{\mfg, g,r, \vec{\varepsilon}} \migi   \overline{\mfO} \mfp_{\mfg, g,r, \chi (\vec{\varepsilon})} 
\end{align}
over $\overline{\mfM}_{g,r}$.
 This morphism restricts to morphisms
\begin{align}
\text{み}_{\mfg, g,r, \vec{\varepsilon}}^{^\mr{Zzz...}} & : \mfM \overline{\mfO} \mfp^{^\mr{Zzz...}}_{\mfg, g,r,  \vec{\varepsilon}} \migi \overline{\mfO} \mfp^{^\mr{Zzz...}}_{\mfg, g,r, \chi (\vec{\varepsilon})},  
\hspace{3mm}
\text{み}_{\mfg, g,r, \vec{\varepsilon}}^{^{p \text{-} \mr{nilp}}} & :  \mfM \overline{\mfO} \mfp^{^{p \text{-}\mr{nilp}}}_{\mfg, g,r,\vec{\varepsilon}} \migi \overline{\mfO} \mfp^{^{p \text{-}\mr{nilp}}}_{\mfg, g,r, \chi (\vec{\varepsilon})}.
\end{align}
We shall refer to $\text{み}_{\mfg, g,r, \vec{\varepsilon}}$ (resp.,  $\text{み}_{\mfg, g,r, \vec{\varepsilon}}^{^\mr{Zzz...}}$; resp.,  $\text{み}_{\mfg, g,r, \vec{\varepsilon}}^{^{p \text{-} \mr{nilp}}}$) as the {\bf universal Miura transformation} (resp., the {\bf universal dormant Miura transformation}; resp., the {\bf universal $p$-nilpotent Miura transformation}).
Moreover, these morphisms fit into the following commutative diagram:
 \begin{align}
\vcenter{\xymatrix{
\mfM \overline{\mfO} \mfp^{^\mr{Zzz...}}_{\mfg, g,r,  \vec{\varepsilon}} \ar[r] \ar[d]^{\text{み}_{\mfg, g,r, \vec{\varepsilon}}^{^\mr{Zzz...}}} & \mfM \overline{\mfO} \mfp^{^{p \text{-}\mr{nilp}}}_{\mfg, g,r,  \vec{\varepsilon}}  \ar[r] \ar[r] \ar[d]^{\text{み}_{\mfg, g,r, \vec{\varepsilon}}^{^{p \text{-} \mr{nilp}}}} & \mfM \overline{\mfO} \mfp_{\mfg, g,r, \vec{\varepsilon}}  \ar[d]^{\text{み}_{\mfg, g,r, \vec{\varepsilon}}} 
\\
\overline{\mfO} \mfp^{^\mr{Zzz...}}_{\mfg, g,r, \chi (\vec{\varepsilon})}  \ar[r] & \overline{\mfO} \mfp^{^{p \text{-}\mr{nilp}}}_{\mfg, g,r, \chi (\vec{\varepsilon})} \ar[r] &  \overline{\mfO} \mfp_{\mfg, g,r,  \chi (\vec{\varepsilon})}, 
}}
\end{align}
where all the  horizontal  arrows  are   closed immersions and   both sides of square diagrams are cartesian.
  \end{rema}

\vspace{3mm}
\bt  \label{th013} \leavevmode\\
\vspace{-5mm}
\begin{itemize}
\item[(i)]
   The isomorphism $\Xi_{\mfg, g,r, \vec{\varepsilon}, p_{-1}}$
   restricts to an isomorphism
 \begin{align} \label{ee280}
 \Xi_{\mfg, g,r, \vec{\varepsilon}, p_{-1}}^{^{p \text{-} \mr{nilp}}} : \mfC \overline{\mfC} \mfo^{^{p \text{-} \mr{nilp}}}_{\mfg, g,r, \vec{\varepsilon}} \isom \mfM \overline{\mfO} \mfp_{\mfg, g,r, \vec{\varepsilon}}^{^{p \text{-} \mr{nilp}}}.
 \end{align}
In particular, both $\mfM \overline{\mfO} \mfp_{\mfg, g,r, \vec{\varepsilon}}^{^{p \text{-} \mr{nilp}}}$ and $\mfM \overline{\mfO} \mfp_{\mfg, g,r, \vec{\varepsilon}}^{^\mr{Zzz...}}$ are finite over $\overline{\mfM}_{g,r}$.
\item[(ii)]
We have the decomposition
\begin{align}
\mfM \overline{\mfO} \mfp_{\mfg, g,r}^{^\mr{Zzz...}} = \coprod_{\vec{\varepsilon} \in \mft (\mbF_p)^{\times r}} \mfM \overline{\mfO} \mfp_{\mfg, g,r, \vec{\varepsilon}}^{^\mr{Zzz...}}.
\end{align}
In particular, $\mfM \overline{\mfO} \mfp_{\mfg, g,r, \vec{\varepsilon}}^{^\mr{Zzz...}}$ is empty unless $\vec{\varepsilon} \in \mft (\mbF_p)^{\times r}$ (or $\vec{\varepsilon} = \emptyset$).
\end{itemize}
 \et
\begin{proof}
Let us consider assertion (i).
By the various definitions involved, 
 the square diagram
\begin{align} \label{DD11}
\begin{CD}
\mfC \overline{\mfC} \mfo_{\mfg, g,r,  \vec{\varepsilon}} @> \Xi_{\mfg, g,r,  \vec{\varepsilon}, p_{-1}} > \sim > \mfM \overline{\mfO} \mfp_{\mfg, g, r,  \vec{\varepsilon}}
\\
@V \Psi_{\mfg, g,r,  \vec{\varepsilon}}
 VV @VV \rotatebox{50}{{\large $\kappa$}}_{\mfb, g,r, \vec{\varepsilon}} V
\\
\ooalign{$\bigoplus$ \cr $\bigotimes_{\mft, g,r}$} @<< \ooalign{$\oplus$ \cr $\otimes_{\mfb \migi \mft, g,r}$} < \ooalign{$\bigoplus$ \cr $\bigotimes_{\mfb, g,r}$}
\end{CD}
\end{align}
is commutative, 
where $\Psi_{\mfg, g,r,  \vec{\varepsilon}} := \Psi_{\mfg, g,r} |_{\mfC \overline{\mfC} \mfo_{\mfg, g,r,  \vec{\varepsilon}}}$ (cf. (\ref{eq010222})) and $\rotatebox{50}{{\large $\kappa$}}_{\mfb, g,r, \vec{\varepsilon}} := \rotatebox{50}{{\large $\kappa$}}_{\mfb, g,r} |_{\mfM \overline{\mfO} \mfp_{\mfg, g, r,  \vec{\varepsilon}}}$ (cf. (\ref{dormant2})).
On the other hand,  $ \mfC \overline{\mfC} \mfo^{^{p \text{-} \mr{nilp}}}_{\mfg, g,r, \vec{\varepsilon}}$ and  $\mfM \overline{\mfO} \mfp_{\mfg, g,r, \vec{\varepsilon}}^{^{p \text{-} \mr{nilp}}}$  are  the  inverse images of 
 $\ooalign{$\oplus$ \cr $\otimes^{-1}_{\mft \migi \mfc, g,r}$}([0]_{\mfc, g,r})$ (which is a closed substack of $\ooalign{$\bigoplus$ \cr $\bigotimes_{\mft, g,r}$}$) via $\Psi_{\mfg, g,r,  \vec{\varepsilon}}$ and 
 $\ooalign{$\oplus$ \cr $\otimes_{\mfb \migi \mft, g,r}$}  \circ \rotatebox{50}{{\large $\kappa$}}_{\mfb, g,r, \vec{\varepsilon}}$ respectively.
Hence,  by  the commutativity of (\ref{DD11}), the isomorphism $\Xi_{\mfg, g,r, \vec{\varepsilon}, p_{-1}}$ restricts to an  isomorphism $\Xi_{\mfg, g,r, \vec{\varepsilon}, p_{-1}}^{^{p \text{-} \mr{nilp}}}$ as of (\ref{ee280}).
Moreover, since $\mfM \overline{\mfO} \mfp_{\mfg, g,r, \vec{\varepsilon}}^{^\mr{Zzz...}}$ is a closed substack of $\mfM \overline{\mfO} \mfp_{\mfg, g,r, \vec{\varepsilon}}^{^{p \text{-} \mr{nilp}}}$, it follows from  the isomorphism $\Xi_{\mfg, \vec{a}, g,r, p_{-1}}^{^{p \text{-} \mr{nilp}}}$ and the result in  Corollary \ref{P230}  that both $\mfM \overline{\mfO} \mfp_{\mfg, g,r, \vec{\varepsilon}}^{^{p \text{-} \mr{nilp}}}$ and $\mfM \overline{\mfO} \mfp_{\mfg, g,r, \vec{\varepsilon}}^{^\mr{Zzz...}}$ are finite over $\overline{\mfM}_{g,r}$.
This completes the proof of assertion (i).

Assertion (ii) follows directly from assertion (i) and Corollary \ref{prff001}.
\end{proof}

\vspace{5mm}
\subsection{From Miura $\mfs \mfl_2$-opers to Miura $\mfg$-opers} \label{sc33gh5}
\leavevmode\\ \vspace{-4mm}

We shall construct a morphism from the moduli stack of generic Miura $\mfs \mfl_2$-opers to the moduli stack of  generic $\mfg$-opers.

Let us write 
\begin{align}
p_0^\circ := \begin{pmatrix} 1 &  0 \\ 0 & -1\end{pmatrix}, \hspace{5mm} p_1^\circ := \begin{pmatrix} 0 &  1 \\ 0 & 0\end{pmatrix}, \hspace{5mm}   p_{-1}^\circ := \begin{pmatrix} 0 &  0 \\ 1 & 0\end{pmatrix}.
\end{align}
The $\mfs \mfl_2$-triple $\{ p_{-1}, 2\check{\rho}, p_1 \}$ given in \S\,\ref{SeSe1}
 determines a Lie algebra homomorphism 
 \begin{align}
 \iota_\mfg : \mfs \mfl_2 \migi \mfg
 \end{align}
  in such a way that 
 $\iota_\mfg (p_0^\circ) = 2 \check{\rho}$, $\iota_\mfg (p_{1}^\circ) = p_1$, and $\iota_\mfg (p_{-1}^\circ) = p_{-1}$.
It follows from ~\cite{Wak5}, Proposition 2.5.1,  that  there exists 
 a unique homomorphism $\iota_\mbB : \mbB_{\mr{PGL}_2} \migi \mbB$ of algebraic groups 
 that is compatible, in the natural sense, with $\iota_\mfg$.
Then, $\iota_\mbB$ restricts to a homomorphism
 $\iota_{\mbT} : \mbT_{\mr{PGL}_2} \migi \mbT$.
 
 Let 
 $T^\mr{log}$ be  an fs log scheme  over $k$ and 
 $U^\mr{log}$ a log-curve over $T^\mr{log}$.
By the definition of $\mcE^\dagger_\mbT := \mcE^\dagger_{\mbT, U^\mr{log}/T^\mr{log}}$, 
 there exists  a canonical isomorphism $\mcE^\dagger_{\mbT_{\mr{PGL}_2}} \times^{\mbT_{\mr{PGL}_2}, \iota_\mbT} \mbT \isom \mcE^\dagger_{\mbT}$.
 Moreover,  the collection of data $(\mcE^\dagger_{\mbB_{\mr{PGL}_2}}, \mcE'^\dagger_{\mbB_{\mr{PGL}_2}}, \eta^\dagger_{\mcE_{\mbB_{\mr{PGL}_2}}})$
 becomes 
 $(\mcE^\dagger_{\mbB}, \mcE'^\dagger_{\mbB}, \eta^\dagger_{\mcE_\mbB})$
 after  change of structure group by $\iota_\mbT$.

Now, let us fix $\vec{\varepsilon} := (\varepsilon_{i})_{i=1}^r \in \mft_{\mr{PGL}_2} (k)^{\times r}$ if $r >0$  (resp.,  $\vec{\varepsilon} := \emptyset$ if $r =0$).
Write  $\iota_\mfg (\vec{\varepsilon}) := (\iota_\mfg (\varepsilon_i))_{i=1}^r \in \mft (k)^{\times r}$ (resp., $\iota_\mfg (\vec{\varepsilon}) := \emptyset$).
Suppose further that we are given a $p_{-1}^\circ$-special Miura $\mfg$-oper $\widehat{\mcE}^{\spadesuit \diamondsuit} := (\mcE_{\mbB_{\mr{PGL}_2}}^\dagger, \nabla_\mcE, {\mcE'}_{\mbB_{\mr{PGL}_2}}^\dagger, \eta^\dagger_{\mcE_{\mbB_{\mr{PGL}_2}}})$  on $U^\mr{log}/T^\mr{log}$ of exponents $\vec{\varepsilon}$.
According to ~\cite{Wak5}, the discussion following Definition 2.7.1,  
the underlying $\mfs \mfl_2$-oper $(\mcE^\dagger_{\mbB_{\mr{PGL}_2}}, \nabla_\mcE)$
 induces  a $\mfg$-oper $(\mcE^\dagger_\mbB, \iota_{\mfg *} (\nabla_\mcE))$ on $U^\mr{log}/T^\mr{log}$
 by  change of structure group.
Moreover, 
the collection of data
\begin{align}
\iota_{\mfg}^* (\widehat{\mcE}^{\spadesuit \diamondsuit}) := (\mcE^\dagger_\mbB, \iota_{\mfg *} (\nabla_\mcE), {\mcE'}^\dagger_\mbB, \eta^\dagger_{\mcE_\mbB})
\end{align}
 forms a $p_{-1}$-spacial Miura $\mfg$-oper on $\mfX$ of exponents $\iota_\mfg (\vec{\varepsilon})$.

\vspace{3mm}
\bpr \label{P01} \leavevmode\\
 \ \ \ 
Let $\vec{\varepsilon}$ be as above.
Then,
the assignment  $\widehat{\mcE}^{\spadesuit \diamondsuit} \mapsto \iota_\mfg^* (\widehat{\mcE}^{\spadesuit \diamondsuit})$ is compatible with base-change over $S$, and hence, determines a morphism
\begin{align} \label{E031}
\mfM \overline{\mfO} \mfp_{\mfs \mfl_2, g,r, \vec{\varepsilon}} \migi \mfM \overline{\mfO} \mfp_{\mfg, g,r,  \iota_\mfg (\vec{\varepsilon})}
\end{align}
over $\overline{\mfM}_{g,r}$, which  restricts to a morphism 
 \begin{align} \label{E031}
\mfM \overline{\mfO} \mfp_{\mfs \mfl_2, g,r, \vec{\varepsilon}}^{^\mr{Zzz...}} \migi \mfM \overline{\mfO} \mfp_{\mfg, g,r, \iota_\mfg (\vec{\varepsilon})}^{^\mr{Zzz...}}.
\end{align}
In particular, $\mfM \overline{\mfO} \mfp_{\mfg, g,r, \iota_\mfg (\vec{\varepsilon})}^{^\mr{Zzz...}}$ is nonempty if  $\mfM \overline{\mfO} \mfp_{\mfs \mfl_2, g,r, \vec{\varepsilon}}^{^\mr{Zzz...}}$ is nonempty (cf. Remark \ref{RRR040}, in which we will discuss the case where $\mfM \overline{\mfO} \mfp_{\mfs \mfl_2, g,r,  \vec{\varepsilon}}^{^\mr{Zzz...}}$ is nonempty).
\epr
\begin{proof}
The assertion follows from the above discussion and ~\cite{Wak5}, Proposition 3.2.3.
\end{proof}

\vspace{10mm}
\section{Miura $\mr{GL}_n$ opers on pointed stable curves} \vspace{3mm}

In this section,
we  describe Miura $\mfs \mfl_n$-opers in terms of vector bundles.
Let us fix  a positive integer $n$, and 
suppose that the characteristic $\mr{char} (k)$ of $k$ is either $0$ or a prime $p$ with $n <p$.
Given a vector bundle $\mcF$ on a scheme, an integer $l$ with $0 \leq l \leq n-1$,  and an $n$-step  decreasing filtration $\mpf := \{ \mcF^j \}_{j =0}^n$ 
 on $\mcF$, we shall write $\mr{gr}_\mpf^l := \mcF^l / \mcF^{l+1}$.
 Also, write $\mr{gr}_\mpf := \bigoplus_{l=0}^{n-1} \mr{gr}_\mpf^l$.

\vspace{5mm}
\subsection{$\mr{GL}_n$-opers on a log-curve} \label{se221}
\leavevmode\\ \vspace{-4mm}

We first recall the definition of a {\it $\mr{GL}_n$-oper} on a log-curve (cf. ~\cite{Wak5}, Definition 4.3.1).
Let  $T^\mr{log}$ be  an fs  log scheme over $k$, and $U^\mr{log}$ a log-curve over $T^\mr{log}$.
\vspace{3mm}
\bde \label{D07} \leavevmode\\
\vspace{-5mm}
\begin{itemize}
\item[(i)]
A {\bf $\mr{GL}_n$-oper} on $U^\mr{log}/T^\mr{log}$ is a collection of data
\begin{align}
\mcF^\heartsuit := (\mcF, \nabla_\mcF, \mathpzc{f}),
\end{align}
where
\begin{itemize}
\vspace{1mm}
\item[$\bullet$]
$\mcF$  denotes a vector bundle on $U$ of rank $n$;
\vspace{1mm}
\item[$\bullet$]
$\nabla_\mcF$ denotes a $T^\mr{log}$-connection $\mcF \migi \Omega_{U^\mr{log}/T^\mr{log}} \otimes \mcF$ on $\mcF$;
\vspace{1mm}
\item[$\bullet$]
$\mpf$ denotes an $n$-step  decreasing filtration $\{ \mcF^j \}_{j =0}^n$ on $\mcF$  consisting of  subbundles
\begin{align}
0 = \mcF^n \subseteq \mcF^{n-1} \subseteq \cdots \subseteq \mcF^0 = \mcF,
\end{align}
\end{itemize}
satisfying the following three conditions:
\begin{itemize}
\vspace{1mm}
\item[(1)]
The subquotients 
$\mr{gr}^j_\mpf \ (= \mcF^j /\mcF^{j+1})$
 ($0 \leq j \leq n-1$) are line bundles;\vspace{1mm}
\item[(2)]
$\nabla_\mcF (\mcF^j) \subseteq \Omega_{U^\mr{log}/T^\mr{log}} \otimes \mcF^{j-1}$ ($1 \leq j \leq n-1$);
\vspace{1mm}
\item[(3)]
The $\mcO_X$-linear morphism
\begin{align} \label{EqEq336}
\mcK \mcS_{\mcF^\heartsuit}^j : 
\mr{gr}^j_\mpf \migi \Omega_{U^\mr{log}/T^\mr{log}} \otimes \mr{gr}^{j-1}_\mpf
\end{align}
defined by assigning $\overline{a} \mapsto \overline{\nabla_\mcF (a)}$ for any local section $a \in \mcF^j$ (where $\overline{(-)}$'s denote the image in the respective quotients), which is well-deined by virtue of the second condition, is an isomorphism.
\end{itemize}
If $U^\mr{log}/T^\mr{log} = X^{\mfX \text{-}\mr{log}}/S^{\mfX \text{-}\mr{log}}$ for some pointed stable curve $\mfX := (X/S,$ $ \{ \sigma_i \}_{i=1}^r)$, then we shall refer to any $\mr{GL}_n$-oper on $X^{\mfX \text{-}\mr{log}}/S^{\mfX \text{-}\mr{log}}$ as a {\bf $\mr{GL}_n$-oper on $\mfX$}.
\item[(ii)]
Let $\mcF^\heartsuit := (\mcF, \nabla_\mcF, \mpf)$ and $\mcG^\heartsuit := (\mcG, \nabla_\mcG, \mpg)$ be $\mr{GL}_n$-opers on $U^\mr{log}/T^\mr{log}$.
An {\bf isomorphism of $\mr{GL}_n$-opers} from $\mcF^\heartsuit$ to $\mcG^\heartsuit$ is an isomorphism
 $(\mcF, \nabla_\mcF) \isom (\mcG, \nabla_\mcG)$ of flat bundles  compatible with the respective filtrations $\mpf$ and $\mpg$.
\end{itemize}
\ede

\vspace{5mm}
\subsection{Miura $\mr{GL}_n$-opers} \label{EqEq333}
\leavevmode\\ \vspace{-4mm}

Next, we define the notion of a {\it Miura $\mr{GL}_n$-oper}.
Let us consider the triple
\begin{align} \label{eq502}
(\mcF, \mpf, \mpf^{-})
\end{align}
consisting of a rank $n$ vector bundle $\mcF$ on $U$ and two $n$-step decreasing filtrations $\mpf := \{ \mcF^j \}_{j=0}^n$, $\mpf^{-} := \{ \mcF^{-, j} \}_{j=0}^n$ on $\mcF$ such that for each $j \in \{ 0, \cdots, n \}$, the composite 
$\alpha_j : \mcF^{-, j}  \migiincl  \mcF \migisurj \mcF /\mcF^{n-j}$ is an isomorphism, equivalently, the composite $\beta_j : \mcF^{j} \migiincl \mcF \migisurj \mcF / \mcF^{-, n-j}$ is an isomorphism.
 The composite 
 \begin{align} \label{EqEq111}
 \mcF^{-, j-1} \xrightarrow[\sim]{\alpha_{j-1}} \mcF/ \mcF^{n-j+1} \migisurj \mcF/\mcF^{n-j} \xrightarrow[\sim]{\alpha_{j}^{-1}} \mcF^{-, j}
 \end{align}
(where $j \in \{ 1, \cdots, n \}$)  determines a split surjection of the short exact sequence $0 \migi \mcF^{-, j} \migi \mcF^{-, j-1} \migi \mr{gr}_{\mpf^{-}}^{j-1} \migi 0$, which gives rise to a decomposition
  \begin{align} \label{EqEq112}
  \gamma_j : \mcF^{-, j-1} \isom \mcF^{-, j} \oplus \mr{gr}_{\mpf^{-}}^{j-1}.
  \end{align}
Also, this fact implies that the kernel of the surjection $\mcF / \mcF^{n-j+1} \migisurj \mcF / \mcF^{n-j}$ (appeared as  the second arrow in (\ref{EqEq111})) is isomorphic to $\mr{gr}_{\mpf^{-}}^{j-1}$.
That is to say, we have a canonical  isomorphism
\begin{align} \label{EqEq339}
\delta_{j, \mpf, \mpf^{-}} : \mr{gr}^{n-j}_\mpf \isom  \mr{gr}_{\mpf^{-}}^{j-1}.
\end{align}
Notice that  $\mcF$ decomposes 
 into the  direct sum of $n$ line bundles $\{ \mr{gr}_{\mpf^{-}}^l \}_{l=0}^{n-1}$ by means of the composite isomorphism
\begin{align}  \label{eq501}
\gamma_{(\mcF, \mpf, \mpf^{-})} :  \mcF \isom \mcF^{-, 1} \oplus \mr{gr}_{\mpf^{-}}^0 \isom  \mcF^{-, 2} \oplus \mr{gr}_{\mpf^{-}}^1 \oplus \mr{gr}_{\mpf^{-}}^0 \isom \cdots \isom \bigoplus_{l =0}^{n-1} \mr{gr}_{\mpf^{-}}^l \ \left(= \mr{gr}_{\mpf^{-}} \right),
\end{align}
where the  $j$-th  isomorphism (for each $j \in \{ 1, \cdots, n \}$) arises from $\gamma_j$.
In the following, we shall consider $\mcF$ as being equipped with a grading (indexed by $\{ 0, \cdots, n-1 \}$)  by means of $\gamma_{(\mcF, \mpf, \mpf^{-})}$,
   i.e., a grading whose $j$-th component 
   is $\mr{gr}_{\mpf^{-}}^{j}$.

Conversely, 
let us consider a rank $n$ vector bundle $\mcF$ on $U$ of the form $\mcF = \bigoplus_{l=0}^{n-1} \mcF_l$,  where each $\mcF_l$ ($l \in \{ 0, \cdots, n-1 \}$) is a line bundle.
For each subset $I$ of $\{ 0, \cdots, n-1 \}$, we shall consider  $\bigoplus_{l \in I} \mcF_l$ as an $\mcO_U$-submodule of $\mcF$  in  the evident manner.
For each $j \in \{ 0, \cdots, n-1 \}$, we shall write
\begin{align}
\mcF^j :=  \bigoplus_{l=0}^{n-1-j} \mcF_l, \hspace{10mm} {\mcF}^{-, j} := \bigoplus_{l= j}^{n-1} \mcF_l.
\end{align}
and write $\mcF^n := 0$, $\mcF^{-, n} :=0$.
The vector bundle $\mcF$ admits 
two decreasing filtrations 
$\mpf$, $\mpf^{-}$
 on $\mcF$ defined as follows:
\begin{align} \label{eq010}
\mpf := \{ \mcF^j\}_{j=0}^n, \hspace{5mm}  \mpf^{-} := \{ {\mcF}^{-, j}  \}_{j=0}^n.
\end{align}
 In particular,
the composite $\mcF^{-, j} \migiincl \mcF \migisurj \mcF/\mcF^{n-j}$ (for each $j \in \{ 0, \cdots, n \}$) is an isomorphism.
In this manner, \textit{ one may  identify each  triple $(\mcF, \mpf, \mpf^{-})$ as in (\ref{eq502}) with 
a rank  $n$ vector bundle decomposing into  a direct sum of $n$ line bundles}.

\vspace{3mm}
\bde \label{EqEq334}\leavevmode\\
 \ \ \ 
\vspace{-5mm}
\begin{itemize}
\item[(i)]
A {\bf (generic) Miura $\mr{GL}_n$-oper}  on $U^\mr{log}/T^\mr{log}$  is a quadruple
\begin{align}
\widehat{\mcF}^\heartsuit  := (\mcF, \nabla_\mcF, \mpf, \mpf^{-}),
\end{align}
where the triple $(\mcF, \mpf, \mpf^{-})$ is as in (\ref{eq502}) and $\nabla_\mcF$ denotes  a $T^\mr{log}$-connection on $\mcF$, such that
\begin{align} \label{eq017}
\nabla_{\mcF} ({\mcF}^{-, j}) \subseteq \Omega_{U^\mr{log}/T^\mr{log}} \otimes ({\mcF}^{-, j})
\end{align}
 for any $j \in \{ 0, \cdots, n \}$ and the triple
\begin{align} \label{f042}
\mcF^\heartsuit := (\mcF, \nabla_{\mcF}, \mpf)
\end{align}
  forms a $\mr{GL}_n$-oper on $U^\mr{log}/T^\mr{log}$.
We shall refer to $\mcF^\heartsuit$ as the {\bf underlying $\mr{GL}_n$-oper} of $\widehat{\mcF}^\heartsuit$.
If $U^\mr{log}/T^\mr{log} = X^{\mfX \text{-} \mr{log}}/ S^{\mfX \text{-} \mr{log}}$ for some pointed stable curve $\mfX := (X/S, \{ \sigma_i \}_{i=1}^r)$, then we shall refer to any (generic) Miura $\mr{GL}_n$-oper on $X^{\mfX \text{-} \mr{log}}/ S^{\mfX \text{-} \mr{log}}$ as a {\bf (generic) Miura $\mr{GL}_n$-oper on $\mfX$}.

 \item[(ii)]
 Let $\widehat{\mcF}^\heartsuit := (\mcF, \nabla_\mcF, \mpf, \mpf^{-})$ and $\widehat{\mcF}^{'\heartsuit} := (\mcF', \nabla_{\mcF'}, \mpf', \mpf'^{-})$
  be Miura $\mr{GL}_n$-opers on $U^\mr{log}/T^\mr{log}$.
 An {\bf isomorphism of (generic) Miura $\mr{GL}_n$-opers} from $\widehat{\mcF}^\heartsuit$ to $\widehat{\mcF}'^{\heartsuit}$ is an isomorphism 
 $\alpha : (\mcF, \nabla_\mcF) \isom (\mcF', \nabla_{\mcF'})$ of log flat bundles
 which restricts to   isomorphisms $\alpha  |_{\mcF^j} : \mcF^j \isom \mcF'^j$ and  $\alpha |_{\mcF^{-, j}} : \mcF^{-, j} \isom \mcF'^{-, j}$ for any $j \in \{ 0, \cdots, n\}$.

  \end{itemize}
 \ede

\bde \label{D7} \leavevmode\\
  \ \ \ 
 Assume  that $k$ has  characteristic $p$ with $n < p$.
We shall say that  a $\mr{GL}_n$-oper $\mcF^\heartsuit := (\mcF, \nabla_\mcF, \mpf)$  (resp., a Miura $\mr{GL}_n$-oper $\widehat{\mcF}^\heartsuit := (\mcF, \nabla_\mcF, \mpf, \mpf^{-})$) 
    on $\mfX$
     is {\bf dormant} 
   if   $\psi^{(\mcF, \nabla_\mcF)} = 0$.
\ede

\begin{rema} \label{TTgh}
\leavevmode\\
 \ \ \  
 If $\widehat{\mcF}^\heartsuit := (\mcF, \nabla_\mcF, \mpf, \mpf^{-})$ is a Miura $\mr{GL}_n$-oper on $U^\mr{log}/T^\mr{log}$.  
By the condition (\ref{eq017}), the $T^\mr{log}$-connection $\nabla_\mcF$ induces, for each $l \in \{ 0, \cdots, n-1 \}$, 
a $T^\mr{log}$-connection 
\begin{align}
\mr{gr}_{\mpf^{-}}^l (\nabla_{\mcF}) : \mr{gr}_{\mpf^{-}}^l \migi \Omega_{U^\mr{log}/T^\mr{log}} \otimes \mr{gr}_{\mpf^{-}}^l 
\end{align}
on the line bundle $\mr{gr}_{\mpf^{-}}^l$.
In particular, we  have a log flat bundle 
\begin{align} \label{eq0106}
\mr{gr}_{\mpf^{-}} (\widehat{\mcF}^\heartsuit) :=  (\mr{gr}_{\mpf^{-}},  \bigoplus_{l =0}^{n-1} \mr{gr}_{\mpf^{-}}^l (\nabla_{\mcF})) \ \left(= \bigoplus_{l=0}^{n-1} (\mr{gr}^l_{\mpf^{-}}, \mr{gr}_{\mpf^{-}}^l (\nabla_{\mcF})) \right)
\end{align}
on $U^\mr{log}/T^\mr{log}$.

Also, let  $\widehat{\mcF}'^{ \heartsuit} := (\mcF', \nabla_{\mcF'}, \mpf', \mpf'^{-})$ be another Miura $\mr{GL}_n$-oper
 on $U^\mr{log}/T^\mr{log}$ and $\alpha : {\widehat{\mcF}}^\heartsuit \isom \widehat{\mcF}'^\heartsuit$  an isomorphism of Miura $\mr{GL}_n$-opers.
 By taking the  gradings $\mr{gr}_{\mpf^{-}}$ and  $\mr{gr}_{\mpf'^{-}}$ of $\mpf^{-}$ and $\mpf'^{-}$ respectively, 
 we obtain an isomorphism
 \begin{align} \label{eq0105}
 \mr{gr}_{} (\alpha) : \mr{gr}_{\mpf^{-}} (\widehat{\mcF}^\heartsuit) \isom \mr{gr}_{\mpf'^{-}} (\widehat{\mcF}'^\heartsuit) 
 \end{align}
 of log flat bundles 
 which makes the  following  square diagram of vector bundles commute:
   \begin{align} \label{EqEq114}
 \begin{CD}
\mcF  @> \alpha > \sim >\mcF'
 \\
 @V \gamma_{(\mcF, \mpf, \mpf^{-})}V \wr V @V \wr V \gamma_{(\mcF', \mpf', \mpf'^{-})} V
 \\
 \mr{gr}_{\mpf^{-}} @> \sim > \mr{gr} (\alpha) > \mr{gr}_{\mpf'^{-}}.
 \end{CD}
 \end{align}
  \end{rema}

\vspace{5mm}
\subsection{Special Miura $\mr{GL}_n$-opers} \label{sc001}
\leavevmode\\ \vspace{-4mm}

Let $U^\mr{log}/T^\mr{log}$ be as before and $\mcN$  a line bundle on $U$.
Recall the rank $n$ vector bundle $\mcF^{[n]\dagger}_\mcN := \bigoplus_{l=0}^{n-1} \mcF^\dagger_{\mcN, l}$ (cf. (\ref{eq0F28})).
We shall write 
\begin{align}
(\mcF^{[n]\dagger}_\mcN, \mpf^\dagger_\mcN, \mpf^{-, \dagger}_\mcN)
\end{align}
for the triple corresponding (via the identification mentioned in the italicized comment preceding Definition \ref{EqEq334}) to the direct sum $\mcF^{[n]\dagger}_\mcN$ of $n$ line bundles (cf. (\ref{eq0F28})).
In particular, 
\begin{align}
\mr{gr}_{\mpf_\mcN^\dagger}^{j} \ (= \mcF^\dagger_{\mcN, n-j-1}) = \mcT^{\otimes (n-j-1)}_{U^\mr{log}/T^\mr{log}} \otimes \mcN,  \hspace{3mm}\mr{gr}_{\mpf_\mcN^{-, \dagger}}^{j} \ (= \mcF^\dagger_{\mcN, j}) = \mcT^{\otimes j}_{U^\mr{log}/T^\mr{log}} \otimes \mcN.
\end{align}

\vspace{3mm}
\bde \label{EqEq335} \leavevmode\\
 \ \ \ 
Let $\mcN^\flat := (\mcN, \nabla_\mcN)$ be a log flat line bundle on $U^\mr{log}/T^\mr{log}$.
We shall say that a Miura $\mr{GL}_n$-oper on $U^\mr{log}/T^\mr{log}$ is   {\bf $\mcN^\flat$-special}
if it is of the form
\begin{align}
\widehat{\mcF}^{\heartsuit \diamondsuit}_\mcN := (\mcF^{[n]\dagger}_{\mcN}, \nabla_\mcF, \mpf^\dagger_\mcN, \mpf^{-, \dagger}_\mcN)
\end{align}
and satisfies the following two conditions:
\begin{itemize}
\item[$\bullet$]
The equality $\mcN^\flat = (\mr{gr}^0_{\mpf^{-}}, \mr{gr}^0_{\mpf^{-}}(\nabla_\mcF))$ holds;
\item[$\bullet$]
For each $j \in \{0, \cdots, n-1\}$,
the automorphism
of $\mcT_{U^\mr{log}/T^\mr{log}}^{\otimes (n-j-1)}\otimes \mcN$ obtained from 
the isomorphism
\begin{align}
\mcK \mcS^j_{\mcF_\mcN^{\heartsuit}} : \mr{gr}_{\mpf_\mcN^\dagger}^j \ \left(= \mcT_{U^\mr{log}/T^\mr{log}}^{\otimes (n-j-1)}\otimes \mcN \right) \isom \Omega_{U^\mr{log}/T^\mr{log}} \otimes \mr{gr}_{\mpf_\mcN^\dagger}^{j-1} \ \left(= \mcT_{U^\mr{log}/T^\mr{log}}^{\otimes (n-j-1)}\otimes \mcN \right)
\end{align}
(cf. (\ref{EqEq336})) coincides with the identity morphism, where $\mcF_\mcN^{\heartsuit}$ denotes the underlying $\mr{GL}_n$-oper of $\widehat{\mcF}^{\heartsuit \diamondsuit}_\mcN$.
\end{itemize}
  \ede

\begin{rema} \label{KKKK}
\leavevmode\\
 \ \ \ 
 If $\widehat{\mcF}^{\heartsuit \diamondsuit} := (\mcF_{\mcN}^{[n] \dagger}, \nabla_\mcF, \mpf^\dagger_\mcN, \mpf^{-, \dagger}_{\mcN})$ is an $\mcN^\flat$-special Miura $\mr{GL}_n$-oper, then 
 the log flat bundle $\mr{gr}_{\mpf^{-, \dagger}_{\mcN}} (\widehat{\mcF}^{\heartsuit \diamondsuit})$  (cf. Remark \ref{TTgh}) forms a $(\mr{GL}_n, \mcN^\flat)$-Cartan connection.
 \end{rema}

\vspace{3mm}
\bpr \label{pr171}\leavevmode\\
 \ \ \ 
Let $\widehat{\mcF}^\heartsuit := (\mcF, \nabla_\mcF, \mpf, \mpf^{-})$ be 
a Miura $\mr{GL}_n$-oper  on $U^\mr{log}/T^\mr{log}$.
Write $\mcN :=  \mr{gr}_{\mpf^{-}}^{0}$
 and $\mcN^\flat := (\mcN, \mr{gr}^0_{\mpf^{-}} (\nabla_\mcF))$.
Then, there exists  a unique   pair $(\widehat{\mcF}^{\heartsuit \diamondsuit}, \iota_{\widehat{\mcF}^\heartsuit})$ consisting of 
\begin{itemize}
\item[$\bullet$]
 an $\mcN^\flat$-special Miura $\mr{GL}_n$-oper $\widehat{\mcF}^{\heartsuit \diamondsuit} := (\mcF^{[n]\dagger}_\mcN, \nabla_{\mcF^\dagger}, \mpf^\dagger_\mcN, \mpf^{-, \dagger}_\mcN)$ on $U^\mr{log}/T^\mr{log}$, and 
 \item[$\bullet$]
 an isomorphism  $\iota_{\widehat{\mcF}^\heartsuit} : \widehat{\mcF}^{\heartsuit} \isom \widehat{\mcF}^{\heartsuit \diamondsuit}$  of Miura $\mr{GL}_n$-opers
 such that the automorphism of $\mcN^\flat$ defined as the $0$-th component of $\mr{gr} (\iota_{\widehat{\mcF}^\heartsuit})$  coincides with the identity morphism.
\end{itemize}
 \epr
\begin{proof}
First, we shall consider the existence assertion.
Let us fix $j \in \{ 1, \cdots, n-1 \}$.
For each $l \in \{ 0, \cdots, j-1 \}$,  we shall write $\iota_{j}^l$ for the composite
\begin{align}
\iota_{j}^l : \mcT^{\otimes l}_{U^\mr{log}/T^\mr{log}} \otimes \mr{gr}_{\mpf^{-}}^{j-l} & \isom \mcT^{\otimes l}_{U^\mr{log}/T^\mr{log}} \otimes \mr{gr}_{\mpf}^{n-j+l -1} \\
& \isom \mcT_{U^\mr{log}/T^\mr{log}}^{\otimes (l+1)} \otimes \mr{gr}_{\mpf}^{n-j+l}\notag  \\
& \isom \mcT_{U^\mr{log}/T^\mr{log}}^{\otimes (l+1)} \otimes \mr{gr}_{\mpf^{-}}^{j-l-1}, \notag
\end{align}
where the first and third arrow denote $\mr{id}_{\mcT^{\otimes l}_{U^\mr{log}/T^\mr{log}}} \otimes \delta_{j -l +1, \mpf, \mpf^{-}}^{-1}$ and $ \mr{id}_{\mcT^{\otimes (l+1)}_{U^\mr{log}/T^\mr{log}}} \otimes \delta_{j-l, \mpf, \mpf^{-}}$ respectively and
 the second arrow denotes  $ \mr{id}_{\mcT^{\otimes (l+1)}_{U^\mr{log}/T^\mr{log}}} \otimes (\mcK \mcS_{\mcF^\heartsuit}^{n-j+l})^{-1}$.
Then, we obtain  the following  composite isomorphism
\begin{align}
\iota_j := \iota_{j}^{j-1} \circ  \iota_{j}^{j-2} \circ   \cdots   \circ \iota_{j}^0  : \mr{gr}_{\mpf^{-}}^j \isom \mcT_{U^\mr{log}/T^\mr{log}}^{\otimes j} \otimes \mr{gr}_{\mpf^{-}}^0 \  (= \mcF_{\mcN, j}^{\dagger}).
\end{align}
Write $\iota_0 := \mr{id}_{\mr{gr}_{\mpf^{-}}^0}$ and write
\begin{align} \label{EqEq360}
\iota_{\widehat{\mcF}^\heartsuit} := \left(\bigoplus_{j=0}^{n-1} \iota_j \right) \circ \gamma_{(\mcF, \mpf, \mpf^{-})} : \mcF \isom \mcF^{[n]\dagger}_\mcN.
\end{align}
Also,  write  $\nabla_{\mcF^\dagger}$ for  the $T^\mr{log}$-connection on $\mcF^{[n]\dagger}_\mcN$ corresponding, via $\iota_{\widehat{\mcF}^\heartsuit}$, to $\nabla_\mcF$.
It follows from  the definitions of $\iota_{\widehat{\mcF}^\heartsuit}$ and  $\nabla_{\mcF^\dagger}$ that
the quadruple 
\begin{align}
\widehat{\mcF}^{\heartsuit \diamondsuit} := (\mcF^{[n]\dagger}_{\mcN}, \nabla_{\mcF^\dagger}, \mpf^{\dagger}_\mcN, \mpf^{-, \dagger}_\mcN)
\end{align} 
specifies a Miura $\mr{GL}_n$-oper  on $U^\mr{log}/T^\mr{log}$ and 
$\iota_{\widehat{\mcF}^\heartsuit}$ specifies  an isomorphism $\widehat{\mcF}^{\heartsuit} \isom \widehat{\mcF}^{\heartsuit \diamondsuit}$ of Miura $\mr{GL}_n$-opers.  
Moreover, if $\mcF^{\heartsuit \diamondsuit}$ denotes the underlying $\mr{GL}_n$-oper of $\widehat{\mcF}^{\heartsuit \diamondsuit}$, then
  for each $j \in \{1, \cdots, n-1 \}$ the following sequence of equalities hold:
\begin{align}
 & \ \ \  \ \, \mcK \mcS^j_{\mcF^{\heartsuit \diamondsuit}} \\
& =  (\mr{id}_{\Omega_{U^\mr{log}/T^\mr{log}}} \otimes\iota_{n-j})\circ (\mr{id}_{\Omega_{U^\mr{log}/T^\mr{log}}}\otimes \delta_{n-j+1, \mpf, \mpf^{-}})\circ \mcK \mcS^j_{\mcF^\heartsuit} \circ \delta_{n-j, \mpf, \mpf^{-}}^{-1}\circ \iota^{-1}_{n-j-1} \notag \\
& = (\mr{id}_{\Omega_{U^\mr{log}/T^\mr{log}}} \otimes\iota^{n-j-1}_{n-j}) \circ \cdots \circ
(\mr{id}_{\Omega_{U^\mr{log}/T^\mr{log}}} \otimes\iota^0_{n-j}) \circ 
(\mr{id}_{\Omega_{U^\mr{log}/T^\mr{log}}} \otimes (\iota^0_{n-j})^{-1})\notag  \\
& \ \ \ \  \circ (\iota_{n-j-1}^0)^{-1}  \circ \cdots \circ  (\iota_{n-j-1}^{n-j-2})^{-1} \notag \\
& = (\mr{id}_{\Omega_{U^\mr{log}/T^\mr{log}}} \otimes\iota^{n-j-1}_{n-j}) \circ \cdots \circ
(\mr{id}_{\Omega_{U^\mr{log}/T^\mr{log}}} \otimes\iota^1_{n-j}) \circ (\iota_{n-j-1}^0)^{-1}  \circ \cdots \circ  (\iota_{n-j-1}^{n-j-2})^{-1} \notag \\
& = (\mr{id}_{\Omega_{U^\mr{log}/T^\mr{log}}} \otimes\iota^{n-j-1}_{n-j}) \circ \cdots \circ
(\mr{id}_{\Omega_{U^\mr{log}/T^\mr{log}}} \otimes\iota^2_{n-j}) \circ (\iota_{n-j-1}^1)^{-1}  \circ \cdots \circ  (\iota_{n-j-1}^{n-j-2})^{-1} \notag \\
& = \cdots \notag \\
& = \mr{id}_{\mr{gr}^{n-j-1}_{f^{-, \dagger}}}, \notag
\end{align}
where the second equality follows from the equality 
\begin{align}
\mr{id}_{\Omega_{U^\mr{log}/T^\mr{log}}} \otimes (\iota^0_{n-j})^{-1} = (\mr{id}_{\Omega_{U^\mr{log}/T^\mr{log}}}\otimes \delta_{n-j+1, \mpf, \mpf^{-}})\circ \mcK \mcS^j_{\mcF^\heartsuit} \circ \delta_{n-j, \mpf, \mpf^{-}}^{-1}
\end{align}
 and  all the equalities after the third  equality  follow from the equalities 
$(\mr{id}_{\Omega_{U^\mr{log}/T^\mr{log}}} \otimes\iota^{l+1}_{n-j}) = \iota_{n-j-1}^l$ ($l = 0, \cdots, n-j-2$).
This implies that the Miura $\mr{GL}_n$-oper $\widehat{\mcF}^{\heartsuit \diamondsuit}$ 
is  $\mcN^\flat$-special, and hence, completes the existence assertion.

The uniqueness assertion  
may be immediately verified.
Indeed, 
let  $(\widehat{\mcF}^{\heartsuit \diamondsuit}, \iota_{\widehat{\mcF}^\heartsuit})$ (where $\widehat{\mcF}^{\heartsuit \diamondsuit} := (\mcF^{[n]\dagger}_\mcN, \nabla_{\mcF^\dagger}, \mpf^\dagger, \mpf^{-, \dagger})$) be a desired  pair.
Since $\mcK \mcS_{\mcF^\heartsuit}^j$ and $\mcK \mcS_{\mcF^{\heartsuit \diamondsuit}}^j$ (for any $j$) are compatible, $\iota_{\widehat{\mcF}^\heartsuit}$ turns out to coincide with the isomorphism ``$\iota_{\widehat{\mcF}^\heartsuit}$" as defined  by (\ref{EqEq360}).
Also, $\nabla_{\mcF^\dagger}$ must be the $T^\mr{log}$-connection corresponding to $\nabla_\mcF$ via $\iota_{\widehat{\mcF}^\heartsuit}$. 
Thus, $(\widehat{\mcF}_\mcN^{\heartsuit \diamondsuit}, \iota_{\widehat{\mcF}^\heartsuit})$ is uniquely determined.
This completes the proof of Proposition \ref{pr171}.
\end{proof}

\vspace{3mm}
\bde \label{D70ff7} \leavevmode\\
  \ \ \ 
  For each Miura $\mr{GL}_n$-oper $\widehat{\mcF}^\heartsuit$ on $U^\mr{log}/T^\mr{log}$, we shall refer to the pair $(\widehat{\mcF}^{\heartsuit \diamondsuit}, \iota_{\widehat{\mcF}^\heartsuit})$
  obtained by applying Proposition \ref{pr171} to $\widehat{\mcF}^\heartsuit$ as the {\bf specialization} of $\widehat{\mcF}^\heartsuit$.
  \ede
\vspace{3mm}

Next,  let us assume  that $n>1$ and  $r >0$.
Let  $S$ be a $k$-scheme, $\mfX := (X/S, \{ \sigma_i \}_{i=1}^r)$  an $r$-pointed stable curve of genus $g$ over $S$,
 and $\widehat{\mcF}^\heartsuit := (\mcF, \nabla_\mcF, \mpf, \mpf^{-})$  a Miura $\mr{GL}_n$-oper on $\mfX$.
 Denote by $(\widehat{\mcF}^{\heartsuit \diamondsuit}, \iota_{\widehat{\mcF}^\heartsuit})$ the specialization of $\widehat{\mcF}^\heartsuit$, where $\widehat{\mcF}^{\heartsuit \diamondsuit} := (\mcF^{[n]\dagger}_\mcN, \nabla_{\mcF^\dagger}, \mpf^{\dagger}_\mcN, \mpf^{-, \dagger}_\mcN)$.
Then, for each $i \in \{1, \cdots, r\}$, we obtain an element
\begin{align} \label{EqEq361}
\epsilon_i^{\widehat{\mcF}^{\heartsuit}} : = \mu_i^{(\mr{gr}_{\mpf^{-, \dagger}_\mcN}, \mr{gr}_{\mpf_\mcN^{-, \dagger}} (\nabla_{\mcF^\dagger}))} \in \Gamma (S, \mcO_S)^{\oplus n} \ (= \mft_{\mr{GL}_n}(S)).
\end{align}

\vspace{3mm}
\bde \label{D70f7} \leavevmode\\
  \ \ \ 
Let $\vec{\epsilon} := (\epsilon_i)_{i=1}^r\in  \mft_{\mr{GL}_n}(S)^{\times r}$.
We shall say that a Miura $\mr{GL}_n$-oper $\widehat{\mcF}^\heartsuit$ is {\bf of exponents $\vec{\epsilon}$} if 
the equality 
$\epsilon_i^{\widehat{\mcF}^{\heartsuit}} = \epsilon_i$ holds  for any $i \in \{ 1, \cdots, r \}$.
If $r =0$, then we shall refer to any Miura $\mr{GL}_n$-oper as being {\bf of exponent $\emptyset$}.
  \ede

\vspace{5mm}
\subsection{From  Miura $\mr{GL}_n$-opers to Miura $\mfs \mfl_n$-opers}
\leavevmode\\ \vspace{-4mm}

In what follows (cf. Theorem \ref{P004}),
we 
construct bijections  between the set of special Miura $\mr{GL}_n$-opers and other objects defined so far.
Let $S$ be a $k$-scheme and $\mfX := (X/S, \{ \sigma_i \}_{i=1}^r)$ an $r$-pointed stable curve over $S$ of genus $g$.
In this subsection, let $\mbG := \mr{PGL}_n$ and denote by $\mbB$ the Borel subgroup of $\mbG$ defined to be the image of the upper triangular matrices via the quotient $\mr{GL}_n \migisurj \mr{PGL}_n$.

First, let $\widehat{\mcF}^\heartsuit := (\mcF, \nabla_\mcF, \mpf, \mpf^{-})$ be an $\mcO_X^\flat$-special  Miura $\mr{GL}_n$-oper  on  $\mfX$.
The underlying $\mr{GL}_n$-oper   $(\mcF, \nabla_\mcF, \mpf)$ induces naturally   an $\mfs \mfl_n (= \mfp \mfg \mfl_n)$-oper $(\mcE_\mbB, \nabla_\mcE)$ on $\mfX$ (cf. ~\cite{Wak5}, the discussion following Lemma 4.4.1).
Also,  the filtration $\mpf^{-}$  determines  a $\mbB$-reduction of $\mcE_\mbG := \mcE_\mbB \times^\mbB \mbG$, that is to say, a   $\mbB$-torsor $\mcE'_\mbB$  together with an isomorphism $\eta_\mcE : \mcE'_\mbB \times^\mbB \mbG  \isom \mcE_\mbG$.
It follows from the definition of a Miura $\mr{GL}_n$-oper that
the quadruple 
\begin{align}
\widehat{\mcF}^{\heartsuit \spadesuit} := (\mcE_\mbB, \nabla_\mcE, \mcE'_\mbB, \eta_\mcE)
\end{align}
forms a generic  Miura $\mfs \mfl_n$-oper on $\mfX$.
Hence, we obtain an assignment $\widehat{\mcF}^\heartsuit \mapsto \widehat{\mcF}^{\heartsuit \spadesuit}$ from 
each
 $\mcO_X^\flat$-special Miura $\mr{GL}_n$-opers on $\mfX$ to 
 a generic Miura $\mfs \mfl_n$-opers on $\mfX$.

Next, let us consider the natural decomposition
\begin{align} \label{eq700}
\mr{End}_{\mcO_X} (\mcF^{[n]\dagger}_{\mcO_X}, \Omega_{X^\mr{log}/S^\mr{log}}\otimes \mcF^{[n]\dagger}_{\mcO_X}) = \bigoplus_{l, l' =0}^{n-1} \mr{End}_{\mcO_X} (\mcF^\dagger_{\mcO_X, l}, \Omega_{X^\mr{log}/S^\mr{log}}\otimes \mcF^\dagger_{\mcO_X, l'}).
\end{align}
The $(l'', l''+1)$-th   component (for each $l'' \in \{ 0, \cdots, n-2\}$) in the right-hand side of (\ref{eq700}) admits 
a natural identification
\begin{align} \label{eq701}
&  \ \ \ \  \mr{End}_{\mcO_X} (\mcF^\dagger_{\mcO_X, l''}, \Omega_{X^\mr{log}/S^\mr{log}}\otimes \mcF^\dagger_{\mcO_X, l''+1}) \\
&= \mr{End}_{\mcO_X} (\mcT_{X^\mr{log}/S^\mr{log}}^{\otimes l''}, \Omega_{X^\mr{log}/S^\mr{log}} \otimes \mcT_{X^\mr{log}/S^\mr{log}}^{\otimes (l''+1)}) \notag \\
&  \cong \Gamma (X, \mcO_X) \notag \\
&  = k. \notag
\end{align}
Let  
\begin{align}
(e_{l, l'})_{l, l'}  \in \bigoplus_{l, l' =0}^{n-1} \mr{End}_{\mcO_X} (\mcF^\dagger_{\mcO_X, l}, \Omega_{X^\mr{log}/S^\mr{log}}\otimes \mcF^\dagger_{\mcO_X, l'})
\end{align}
be the element defined  in such a way that 
 $e_{l, l'} = 1$ (via the identification (\ref{eq701})) if 
 $l + 1 = l'$
  and $e_{l, l'} =0$ if otherwise.
The element $(e_{l, l'})_{l, l'}$ corresponds, via  (\ref{eq700}),
to an $\mcO_X$-linear morphism
\begin{align}
\nabla_0 : \mcF^{[n]\dagger}_{\mcO_X} \migi \Omega_{X^\mr{log}/S^\mr{log}}\otimes \mcF^{[n]\dagger}_{\mcO_X}.
\end{align}

Now, let $\widehat{\mcF}^\clubsuit := (\mcF_{\mcO_X}^{[n]\dagger}, \bigoplus_{l =0}^{n-1} \nabla_l)$
   be
a $(\mr{GL}_n, \mcO_X^\flat)$-Cartan connection on $\mfX$.
The sum $\nabla_0 + \bigoplus_{l =0}^{n-1}\nabla_l$ specifies   an $S^\mr{log}$-connection on $\mcF^{[n]\dagger}_{\mcO_X}$, and moreover, 
the quadruple 
\begin{align} \label{EEE123}
\widehat{\mcF}^{\clubsuit \heartsuit} := (\mcF^{[n]\dagger}_{\mcO_X}, \nabla_0 + \bigoplus_{l=0}^{n-1} \nabla_l, \mpf_{\mcO_X}^\dagger, \mpf_{\mcO_X}^{-, \dagger})
\end{align}
forms 
an $\mcO_X^\flat$-special Miura $\mr{GL}_n$-oper on $\mfX$.
One verifies immediately that 
the  assignment  $\widehat{\mcF}^\clubsuit \mapsto \widehat{\mcF}^{\clubsuit \heartsuit}$ 
determines a functorial (with respect to $S$) bijection from the set of 
$(\mr{GL}_n, \mcO_X^\flat)$-Cartan connections on $\mfX$ to the set of isomorphism classes of   $\mcO_X^\flat$-special Miura $\mr{GL}_n$-opers on $\mfX$.
Indeed, the inverse map is given by assigning, to each $\mcO_X^\flat$-special Miura $\mr{GL}_n$-oper 
$\widehat{\mcF}^\heartsuit$, the $(\mr{GL}_n, \mcO_X^\flat)$-Cartan connection $\mr{gr}_{\mpf_{\mcO_X}^{-, \dagger}} (\widehat{\mcF}^{\heartsuit \diamondsuit})$ (cf. Remark \ref{KKKK}), where let  $(\widehat{\mcF}^{\heartsuit \diamondsuit}, \iota_{\widehat{\mcF}^\heartsuit})$ be  the specialization of $\widehat{\mcF}^\heartsuit$.

Here, for each $\vec{\varepsilon} := ((\varepsilon_{i1}, \varepsilon_{i2}, \cdots, \varepsilon_{i (n-1)}))_{i=1}^r  \in (\Gamma (S, \mcO_S)^{\oplus (n-1)})^{\times r}$,
we write
\begin{align} \label{ee201}
\vec{\varepsilon}^{\, +0} := ((0, \varepsilon_{i1}, \varepsilon_{i2}, \cdots, \varepsilon_{i (n-1)}))_{i=1}^r \in (\Gamma (S, \mcO_S)^{\oplus n})^{\times r} = \mft_{\mr{GL}_n} (S)^{\times r}
\end{align}
and 
\begin{align} \label{ee480}
[\vec{\varepsilon}\,] := \pi (\vec{\varepsilon}^{\, +0})
\end{align}
(cf. (\ref{ee300})  for the definition of $\pi$).
 If $r =0$ and $\vec{\varepsilon} := \emptyset$, then we write $\vec{\varepsilon}^{\,+0} := \emptyset$ and $[\vec{\varepsilon}\,] := \emptyset$.

\vspace{3mm}
\bt \label{P004} \leavevmode\\
 \ \ \ 
Let  $\vec{\varepsilon}  \in (\Gamma (S, \mcO_S)^{\oplus (n-1)})^{\times r}$ (where $\vec{\varepsilon} := \emptyset$ if $r=0$).
Then, the following square diagram is commutative (cf. Remark \ref{RRRH} for the construction of the lower horizontal morphism): 
   \begin{align} \label{efe0991}
\begin{CD}
\begin{pmatrix}
\text{the set of $(\mr{GL}_n, \mcO_X^\flat)$-} \\
\text{Cartan connections} \\
\text{on $\mfX$ of monodromies  $\vec{\varepsilon}^{\,+0}$} 
\end{pmatrix}
@>  \widehat{\mcF}^\clubsuit \mapsto \widehat{\mcF}^{\clubsuit \heartsuit}> \sim >
\begin{pmatrix}
\text{the set  of isomorphism classes} \\
\text{of $\mcO_X^\flat$-special Miura $\mr{GL}_n$-opers} \\
\text{on $\mfX$ of exponents  $\vec{\varepsilon}^{\,+0}$}
\end{pmatrix}
 \\
@V  V \wr V @V  V \widehat{\mcF}^\heartsuit \mapsto  \widehat{\mcF}^{\heartsuit \spadesuit}V
\\
\begin{pmatrix}
\text{the set  of $\mfs \mfl_n$-Cartan} \\
\text{connections on $\mfX$ of} \\
\text{monodromies $[\vec{\varepsilon}\,]$} \\
\end{pmatrix}
 @> \sim >  (\ref{equicat2}) > \begin{pmatrix}
\text{the set of isomorphism } \\
\text{classes of Miura $\mfs \mfl_n$-opers}  \\
\text{on $\mfX$ of exponents $[\vec{\varepsilon}\,]$}
\end{pmatrix},
\end{CD}
\end{align}
  where 
the left-hand vertical arrow denotes the second arrow in (\ref{ee290}).
In particular, the right-hand vertical arrow  
$\widehat{\mcF}^\heartsuit \mapsto  \widehat{\mcF}^{\heartsuit \spadesuit}$ turns out to be bijective.
\et
\begin{proof}
The assertion follows from the various definitions of maps involved.
\end{proof}
\vspace{3mm}

\vspace{10mm}
\section{Pre-Tango structures on a log-curve} \label{fdaa} \vspace{3mm}

In this section, we recall the notion of a {\it Tango structures}  on a smooth curve (cf., e.g., ~\cite{Tan}, ~\cite{TakYok} for various discussions concerning
Tango structures) and  study  its generalization, i.e., a {\it pre-Tango structure on a pointed stable curve}.
Our purpose is  to understand the relationship between pre-Tango structures and generic Miura opers.
The main result of this section (cf. Theorem  \ref{Td01}) asserts that 
there exists a canonical bijective correspondence between the set of pre-Tango structures (of  monodromies $\vec{\varepsilon}$)  and the set of isomorphism classes of dormant generic  Miura $\mfs \mfl_2$-opers (of exponents $\vec{\varepsilon}$).

\vspace{5mm}
\subsection{Tango structures on a smooth curve} \label{SeSe444}
\leavevmode\\ \vspace{-4mm}

In this section, we assume that $k$ has characteristic $p > 2$.
We first recall  the definition of a Tango structure on a smooth curve.
Let $T$ be a $k$-scheme and $U$ a smooth curve over $T$.
Denote by $\mcB_{U/T}$ ($\subseteq \Omega_{U/T}$) the sheaf of locally exact $1$-forms on $U$ relative to $T$.
The direct image $F_{U/T*} (\mcB_{U/T})$ (cf. \S\,\ref{sc091} for the definition of $F_{U/T}$) forms a vector bundle on $U^{(1)}_T$ of rank $p-1$.

Now, let $\mcL$ be a line subbundle  of  $F_{U/T*}(\mcB_{U/T})$.
Consider the $\mcO_{U^{(1)}_T}$-linear  composite 
\begin{align}
\mcL \migiincl F_{U/S*} (\mcB_{U/S}) \migiincl  F_{U/S*} (\Omega_{U/S}),
\end{align}
where the first arrow denotes the natural inclusion 
 and the  second arrow denotes  the morphism obtained by applying the functor $ F_{U/S*} (-)$ to the natural inclusion $\mcB_{U/S} \migiincl \Omega_{U/S}$.
This composite corresponds to a morphism
\begin{align} \label{E101}
\xi_\mcL : F^*_{U/T} (\mcL) \migi \Omega_{U/T}
\end{align}
 via   the adjunction relation ``$F_{U/T}^*(-) \dashv F_{U/T*}(-)$".

\vspace{3mm}
\bde \label{D13} \leavevmode\\
 \ \ \ 
We shall say that 
$\mcL$ is 
a {\bf Tango structure} on $U/T$ if the morphism 
$\xi_\mcL$ is an isomorphism.
 \ede

\begin{rema} \label{R001}
\leavevmode\\
 \ \ \ 
Suppose that $X$ is a  geometrically connected, proper,  and smooth curve over $k$ of genus $g$ ($> 1$).
Also, suppose that $X/k$ admits a Tango structure $\mcL$ ($\subseteq F_{X/k*} (\mcB_{X/k})$).
Then, since $F^*_{X/k} (\mcL) \cong \Omega_{X/k}$, 
the following equalities hold:
\begin{align}
\mr{deg} (\mcL) = \frac{1}{p} \cdot \mr{deg} (F_{X/S}^* (\mcL)) = \frac{1}{p} \cdot \mr{deg} (\Omega_{X/S}) = \frac{2 (g-1)}{p}.
\end{align} 
In particular, {\it  if $p \nmid g-1$, then there is no Tango structure on $X/k$}.
\end{rema}
\vspace{3mm}

\begin{rema} \label{RR44}
\leavevmode\\
 \ \ \ Let $X$ be as in Remark \ref{R001} above.
Then, {\it a Tango structure $\mcL \subseteq F_{X/k *} (\mcB_{X/k})$ on $X/k$ is completely determined  the isomorphism class  $[\mcL]$ of the underlying  line bundle $\mcL$}.
Indeed, if $\mcL \subseteq F_{X/k *} (\mcB_{X/k})$ and   $\mcL' \subseteq F_{X/k *} (\mcB_{X/k})$ are  Tango structures on $X/k$ admitting an isomorphism $\mcL \isom \mcL'$.
The map of sets $\mr{Isom}_{\mcO_{X^{(1)}_k}} (\mcL, \mcL') \migi \mr{Isom}_{\mcO_{X}} (F_{X/k}^*(\mcL), F_{X/k}^*(\mcL'))$ induced by pull-back via $F_{X/k}$ is bijective because it may be identified with the $p$-power map on $k^\times$ (under an identification  $\mcL \cong \mcL'$).
Hence, the isomorphism $\xi^{-1}_{\mcL'} \circ \xi_\mcL : F_{X/k}^*(\mcL) \isom  F_{X/k}^*(\mcL')$ comes from a unique isomorphism $\xi : \mcL \isom \mcL'$,  which is verified (by  the adjunction relation ``$F_{X/k}^*(-) \dashv F_{X/k*}(-)$") to be compatible with the respective inclusions $\mcL \subseteq F_{X/k *} (\mcB_{X/k})$ and  $\mcL' \subseteq F_{X/k *} (\mcB_{X/k})$.
This implies that $\mcL$ and $\mcL'$ specify the same Tango structure. 
\end{rema}

\vspace{5mm}
\subsection{From Tango structures to dormant  Miura $\mr{GL}_2$-opers} \label{SeSe445}
\leavevmode\\ \vspace{-4mm}

In what follows,
 we shall discuss a construction of dormant Miura $\mr{GL}_2$-oper on $U/T$ by means of a Tango structure.
Let $\mcL$  ($\subseteq F_{U/T*}(\mcB_{U/T})$) be a Tango structure on $U/T$.
Denote by $\mcG_\mcL$ the inverse image of $\mcL$ via the surjection $F_{U/T*} (\mcO_U) \migi F_{U/T*}(\mcB_{U/T})$ induced from the universal derivation $d : \mcO_U \migi \Omega_{U/T}$.
We have an inclusion of short exact sequences:
\begin{align} \label{E2}
\vcenter{\xymatrix{
0 \ar[r] & \mcO_{U^{(1)}_T} 
 \ar@{}[d]|{| |} 
 \ar[r]  & \mcG_\mcL \ar@{}[d]|{\bigcap} \ar[r]  & \mcL \ar@{}[d]|{\bigcap}  \ar[r]  & 0
\\
0 \ar[r] & \mcO_{U^{(1)}_T} \ar[r] & F_{U/T*}(\mcO_U) \ar[r] & F_{U/T*}(\mcB_{U/T}) \ar[r] & 0.
}}
\end{align}
Consider  the following inclusion of short exact sequences defined to be the pull-back of
 (\ref{E2}) via 
  $F_{U/T}$:
\begin{align} \label{EE2}
\vcenter{\xymatrix{
0 \ar[r] & \mcO_U 
 \ar@{}[d]|{| |} 
 \ar[r]  & F^*_{U/T}(\mcG_\mcL) \ar@{}[d]|{\bigcap} \ar[r]  & F^*_{U/T}(\mcL) \ar@{}[d]|{\bigcap}  \ar[r]  & 0
\\
0 \ar[r] & \mcO_U \ar[r] & F^*_{U/T}(F_{U/T*}(\mcO_U)) \ar[r] & F^*_{U/T} (F_{U/T*}(\mcB_{U/T})) \ar[r] & 0.
}}
\end{align}
Then, the morphism 
\begin{align} \label{EE10}
F^*_{U/T}(F_{U/T*} (\mcO_U)) \migi \mcO_U
\end{align}
 corresponding, via the adjunction relation ``$F_{U/T}^*(-) \dashv F_{U/T*}(-)$", to the identity morphism $F_{U/T*} (\mcO_U) \migi F_{U/T*} (\mcO_U)$ is verified to be
a split surjection of the lower sequence in (\ref{EE2}).
This split surjection determines 
 a decomposition
\begin{align} \label{EE04}
F^*_{U/T}(F_{U/T*}(\mcO_U))  \isom F^*_{U/T} (F_{U/T*}(\mcB_{U/T})) \oplus  \mcO_U,
\end{align}
which restricts to an isomorphism
\begin{align} \label{EqEq1111}
F^*_{U/T} (\mcG_\mcL)  \isom  F^*_{U/T} (\mcL) \oplus \mcO_U.
\end{align}
Hence, we obtain a composite isomorphism
\begin{align} \label{EE05}
F^*_{U/T} (\mcG_\mcL \otimes \mcL^\vee)  & \isom F^*_{U/T} (\mcG_\mcL)  \otimes F^*_{U/T} (\mcL^\vee)  \\
& \isom  (F^*_{U/T} (\mcL) \oplus \mcO_U) \otimes  F^*_{U/T} (\mcL^\vee) \notag  \\
 &  \isom \mcO_U \oplus F^*_{U/T} (\mcL^\vee)  \notag \\
  & \isom \mcO_U \oplus  \mcT_{U/T}  \notag \\
 & \stackrel{\varsigma}{\isom} \mcO_U \oplus  \mcT_{U/T}  \left(= \mcF^{[2]\dagger}_{\mcO_U}\right), \notag
\end{align}
where the second and fourth  arrows arise from (\ref{EqEq1111}) and $\xi_\mcL$ respectively and the fifth 
arrow  $\varsigma$
arises from  both $\mr{id}_{\mcO_U}$ and the automophism of  $\mcT_{U/T}$ determined by multiplication by $(-1)$.
Denote by $\nabla_\mcL$ the $T$-connection on $\mcF^{[2]\dagger}_{\mcO_U}$ corresponding to the $T$-connection $\nabla^\mr{can}_{\mcG_\mcL \otimes \mcL^\vee}$ (resp., (\ref{eeqq1})) via the composite isomorphism  (\ref{EE05}).
Thus, we obtain a quadruple
\begin{align}
\widehat{\mcT}\hspace{-1mm}{\it an}^{\heartsuit \diamondsuit}_{\mcL} := (\mcF_{\mcO_U}^{[2]\dagger}, \nabla_\mcL, \mpf^\dagger_{\mcO_U}, \mpf^{-, \dagger}_{\mcO_U}).
\end{align}

\vspace{3mm}
\bpr \label{PP020} \leavevmode\\
 \ \ \ 
The quadruple $\widehat{\mcT}\hspace{-1mm}{\it an}^{\heartsuit \diamondsuit}_{\mcL}$ 
forms a dormant $\mcO_X^\flat$-special  Miura $\mr{GL}_2$-oper on $U/T$. 
Moreover, let us  write 
\begin{align}
\widehat{\mcF}^\clubsuit_\mcL := (\mcF_{\mcO_U}^{[2]\dagger}, d \oplus \xi^{\vee *}_\mcL (\nabla^\mr{can}_{\mcL^\vee})),
\end{align}
which specifies  
a $(\mr{GL}_2, \mcO_U^\flat)$-Cartan connection on $U/T$.
Then,   
there exists a natural isomorphism
\begin{align}
\mr{gr}_{\mpf^{-, \dagger}_{\mcO_U}} (\widehat{\mcT}\hspace{-1mm}{\it an}^{\heartsuit \diamondsuit}_{\mcL}) \isom \widehat{\mcF}^{\clubsuit}_\mcL.
\end{align}
In particular, $\widehat{\mcT}\hspace{-1mm}{\it an}^{\heartsuit \diamondsuit}_{\mcL}$ is isomorphic to $\widehat{\mcF}^{\clubsuit \heartsuit}_\mcL$ (cf.  (\ref{EEE123})).
\epr
\begin{proof}
Let us prove the first assertion.
Denote by $\widehat{\mcT}\hspace{-1mm}{{\it an}_\mcL^{\heartsuit}}'$ the $\mr{GL}_n$-oper defined as
$\widehat{\mcT}\hspace{-1mm}{\it an}^{\heartsuit}_{\mcL}$ (i.e., the underlying $\mr{GL}_n$-oper  of $\widehat{\mcT}\hspace{-1mm}{\it an}^{\heartsuit \diamondsuit}_{\mcL}$)  tensored with
the log flat line bundle  $(\Omega_{U/T}, \xi_{\mcL *} (\nabla^\mr{can}_\mcL))$ (where $\xi_{\mcL *} (\nabla^\mr{can}_\mcL)$ denotes the $T$-connection on $\Omega_{U/T}$ corresponding to $\nabla^\mr{can}_\mcL$ via $\xi_\mcL$).
 In the following, let us use the various notations (e.g., $\alpha$, $\beta$, $\gamma, \cdots$)  defined in Lemma \ref{Le020} below.
 By the definition of $\widehat{\mcT}\hspace{-1mm}{{\it an}_\mcL^{\heartsuit \diamondsuit}}$,
the following sequence of equalities holds:
 \begin{align} \label{EqEq2345}
 \mcK \mcS^1_{\widehat{\mcT}\hspace{-1mm}{{\it an}_\mcL^{\heartsuit}}'} & =  \epsilon^{-1} \circ \beta |_{\Omega_{U/T}\otimes F^*_{U/T} (\mcG_\mcL)}\circ \nabla^\mr{can}_{\mcG_\mcL}\circ  \alpha |_{F_{U/T}^* (\mcL)} \circ \xi_\mcL^{-1} \\
 & = \epsilon^{-1} \circ \beta \circ \nabla^\mr{can}_{F_{U/T*}(\mcO_U)} \circ \alpha \circ \mr{incl} \circ \xi_\mcL^{-1} \notag  \\
 & = \delta \circ \gamma \circ \mr{incl} \circ \xi_\mcL^{-1} \notag \\
 & =  \mr{id}_{\Omega_{U/T}}, \notag 
 \end{align}
where $\mr{incl}$ denotes the inclusion $F_{U/T}^* (\mcL) \migiincl F_{U/T}^* (F_{U/T*}(\mcB_{U/T}))$,
the third equality follows from Lemma \ref{Le020} below,  and the last equality follows from the definition of $\xi_\mcL$.
 This implies the equality  $\mcK \mcS_{\widehat{\mcT}\hspace{-1mm}{{\it an}_\mcL^{\heartsuit}}}^1 = \mr{id}_{\mcO_U}$, and hence, completes the first  assertion.
 
The second  assertion follows  directly from  the definition of  $\widehat{\mcT}\hspace{-1mm}{\it an}^{\heartsuit \diamondsuit}_{\mcL}$.
Also, the third assertion follows from the bijectivity of the upper horizontal arrow in (\ref{efe0991}).
\end{proof}
\vspace{3mm}

The following lemma was used to prove the above proposition.

\vspace{3mm}
\ble \label{Le020} \leavevmode\\
 \ \ \ 
The following diagram is commutative:
\begin{align} \label{Edd2}
\vcenter{\xymatrix@C=38pt{
F^*_{U/T} (F_{U/T*} (\mcB_{U/T}))  \ar[r]^{\alpha}  \ar[d]_{\gamma} & F^*_{U/T} (F_{U/T*} (\mcO_U))  \ar[r]^-{\nabla^\mr{can}_{F_{U/T*}(\mcO_U)}} &  \Omega_{U/T} \otimes F^*_{U/T} (F_{U/T*} (\mcO_U))  \ar[d]^{\beta} \\
 F^*_{U/T} (F_{U/T*} (\Omega_{U/T}))  \ar[r]_{\hspace{10mm}\delta} & \Omega_{U/T}  \ar[r]_{\hspace{-15mm}\epsilon}  & \Omega_{U/T} \otimes \mcO_U \ \left(= \Omega_{U/T} \right),
}}
\end{align}
where
\begin{itemize}
\item[$\bullet$]
$\alpha$ denotes the split injection of the lower sequence in (\ref{EE2}) corresponding to
the split surjection (\ref{EE10});
\item[$\bullet$]
$\beta$ denotes the tensor product of (\ref{EE10}) and the identity morphism of $\Omega_{U/T}$;
\item[$\bullet$]
$\gamma$ denotes the injection arising from the natural inclusion $\mcB_{U/T} \migiincl \Omega_{U/T}$;
\item[$\bullet$]
$\delta$ denotes  the morphism corresponding to the identity morphism $F_{U/S*} (\Omega_{U/T})$ $\migi$ $F_{U/S*} (\Omega_{U/T})$
via the adjunction relation ``$F_{U/T}^*(-) \dashv F_{U/T*}(-)$";
\item[$\bullet$]
$\epsilon$ denotes the  automorphism given by multiplication by $(-1)$.
\end{itemize}

\ele
\begin{proof}
Since $U$ is flat over $T$, it suffices to prove Lemma \ref{Le020} after restricting $U$ to its fibers over various geometric  points of $T$. 
Hence, one may assume that $k$ is algebraically closed and $T = \mr{Spec} (k)$.
Moreover, as the statement is of local nature, 
 one may replace $U$ by $\mr{Spec} (k[[x]])$ and $\Omega_{U/T}$ by  the $\mcO_U$-module given by  the $k[[x]]$-module  $k[[x]] dx$.
Observe that for each  $\mcO_U$-module $\mcM$  obtained from  some  $k[[x]]$-module $M$,  we have 
\begin{align}
\Gamma (U, F^*_{U/T}(F_{U/T*} (\mcM))) = k [[x]] \otimes_{k[[x^p]]} M.
\end{align}
Then, 
 $F^*_{U/T}(F_{U/T*} (\mcB_{U/T}))$ ($\subseteq F^*_{U/T}(F_{U/T*} (\Omega_{U/T}))$) 
 corresponds to the $k [[x]]$-modules
$\bigoplus_{l=0}^{p-2} k [[x]] \cdot  1 \otimes x^l dx$,  and 
  the kernel of (\ref{EE10}) corresponds to
 $\bigoplus_{l=1}^{p-1} $ $k [[x]](1 \otimes x^l - x^l \otimes 1)$.
The morphisms $\alpha$ and $\nabla^\mr{can}_{F_{U/T*}(\mcO_U)}$  are  given by assigning  $1 \otimes x^l dx \mapsto \frac{1}{l+1} \cdot (1 \otimes x^{l+1} - x^{l+1} \otimes 1)$ and $a \otimes b \mapsto da \otimes b$ (for any $a$, $b \in k[[x]]$) respectively.
Hence,  the following equalities hold:
\begin{align}
& \ \ \ \   \beta (\nabla^\mr{can}_{F_{U/T*}(\mcO_U)}(\alpha (1 \otimes x^l dx)))  \\
& = \beta \left(\nabla^\mr{can}_{F_{U/T*}(\mcO_U)} \left(\frac{1}{l+1}
 \cdot \left(1 \otimes x^{l+1} - x^{l+1} \otimes 1\right) \right) \right)  \notag  \\
 & = \beta (-x^l dx \otimes 1) \notag  \\
 & = - x^l dx.\notag
\end{align}
On the other hand, we have
\begin{align}
\epsilon (\delta (\gamma (1 \otimes x^l))) & = \epsilon (\delta (1 \otimes x^l dx))  = \epsilon (x^l dx) =  -x^l dx. 
\end{align}
Thus, the equality $\beta \circ \nabla^\mr{can}_{F_{U/T*}(\mcO_U)} \circ \alpha = \epsilon \circ \delta \circ \gamma$ holds.
 This completes the proof of Lemma \ref{Le020}.
\end{proof}
\vspace{3mm}

\vspace{5mm}
\subsection{Pre-Tango structures on a log-curve}
\leavevmode\\ \vspace{-4mm}

In this section, we shall discuss the definition of a pre-Tango structure on a pointed stable curve.
Let $T^\mr{log}$ be an fs log scheme over $k$ and $U^\mr{log}$ a log-curve over $T^\mr{log}$.

Here, recall the {\it Cartier operator} 
\begin{align} \label{GGh}
C_{U^\mr{log}/T^\mr{log}} : F_{U/T*} (\Omega_{U^\mr{log}/T^\mr{log}}) \migi \Omega_{U^{(1) \mr{log}}_T/T^\mr{log}}
\end{align}
  of $U^\mr{log}/T^\mr{log}$.
  That is to say, $C_{U^\mr{log}/T^\mr{log}}$ is  a unique  $\mcO_{U^{(1)}_T}$-linear morphism   whose composite with the inclusion 
  $\Omega_{U^{(1) \mr{log}}_T/T^\mr{log}} \migi \Omega_{U^{(1) \mr{log}}_T/T^\mr{log}} \otimes F_{U/T*} (\mcO_U)$ induced by the natural injection $\mcO_{U^{(1)}_T} \migi F_{U/T*} (\mcO_U)$ coincides with the Cartier operator associated with $(\mcO_U, d)$ in the sense of ~\cite{Og}, Proposition 1.2.4.

\vspace{3mm}
\bde \label{D113} \leavevmode\\
  \ \ \  A {\bf pre-Tango structure} on $U^\mr{log}/T^\mr{log}$ is a $T^\mr{log}$-connection $\nabla_\Omega$
 on $\Omega_{U^\mr{log}/T^\mr{log}}$ with vanishing $p$-curvature
satisfying that $F_{U/T*}(\mr{Ker} (\nabla_\Omega)) \subseteq \mr{Ker} (C_{U^\mr{log}/T^\mr{log}})$.
 If $U^\mr{log}/T^\mr{log} = X^{\mfX \text{-} \mr{log}}/S^{\mfX \text{-} \mr{log}}$ for a pointed stable curve $\mfX := (X/S, \{ \sigma_i \}_{i=1}^r)$, then we shall refer to a pre-Tango structure on $X^{\mfX \text{-} \mr{log}}/S^{\mfX \text{-} \mr{log}}$ as a {\bf pre-Tango structure on $\mfX$}.
\ede
\vspace{3mm}

If the curve $U^\mr{log}/T^\mr{log}$ under consideration is a non-logarithmic smooth curve, then the notion of a pre-Tango structure is equivalent to the notion of a Tango structure defined in  Definition \ref{D13}, as verified  in the following proposition.

\vspace{3mm}
\bpr \label{Pff110} \leavevmode\\
 \ \ \ 
Let $U/T$ be as in \S\S\,\ref{SeSe444}-\ref{SeSe445}.
Then, the assignment $\nabla_\Omega \mapsto F_{U/T*} (\mr{Ker} (\nabla_\Omega))$ determines a bijection of sets
  \begin{equation} \label{EE09ff8}
\begin{pmatrix}
\text{the set of pre-Tango} \\
\text{structures on $U/T$}
\end{pmatrix}
\isom
\begin{pmatrix}
\text{the set  of Tango} \\
\text{structures on $U/T$}\\
\end{pmatrix}.
\end{equation}
\epr
\begin{proof}
First, we shall prove the claim that the assignment $\nabla_\Omega \mapsto F_{U/T*} (\mr{Ker} (\nabla_\Omega))$ defines a map from the set of pre-Tango structures on $U/T$ to the set of Tango structures on $U/T$.
Recall from  ~\cite{K}, \S\,5, p.\,190, Theorem 5.1, that the assignments $\mcV \mapsto (F^*_{U/T}(\mcV), \nabla^\mr{can}_\mcV)$ and $(\mcF, \nabla_\mcF) \mapsto F_{U/T*} (\mr{Ker} (\nabla_\mcF))$ determines an equivalence of categories 
  \begin{equation} \label{ED09ff8}
\begin{pmatrix}
\text{the category of} \\
\text{vector bundles on $U^{(1)}_T$}
\end{pmatrix}
\isom
\begin{pmatrix}
\text{the category of} \\
\text{flat bundles on $U/T$}\\
\text{with vanishing $p$-curvature}
\end{pmatrix}.
\end{equation}
Now, let
 $\nabla_\Omega$
be a pre-Tango structure on $U/T$.
The equivalence of categories recalled above
 implies that $F_{U/T*} (\mr{Ker}(\nabla_\Omega))$ is a line bundle on $U^{(1)}_T$ and the morphism 
\begin{align} \label{EqEqEq1}
F^*_{U/T} (F_{U/T*} (\mr{Ker} (\nabla_\Omega))) \migi \Omega_{U/T}
\end{align}
 corresponding to the inclusion $F_{U/T*} (\mr{Ker}(\nabla_\Omega)) \migiincl F_{U/T*} (\Omega_{U/T})$ (via the adjunction relation ``$F_{U/T}^*(-) \dashv F_{U/T*}(-)$") is an isomorphism.
Also, since the sequence 
\begin{align}\label{E0001}
0 \migi F_{U/T*}(\mcB_{U/T}) \migi F_{U/T*}(\Omega_{U/T}) \xrightarrow{C_{U/T}} \Omega_{U_T^{(1)}/T} \migi 0
\end{align}
is exact, the inclusion $F_{U/T*}(\mr{Ker} (\nabla_\Omega)) \migiincl F_{U/T*}(\Omega_{U/T})$
factors through  the inclusion $F_{U/T*}(\mcB_{U/T}) \migi F_{U/T*}(\Omega_{U/T})$.
By taking account of  the resulting injection $F_{U/T*}(\mr{Ker} (\nabla_\Omega)) \migiincl F_{U/T*} (\mcB_{U/T})$, one may   regard $F_{U/T*}(\mr{Ker} (\nabla_\Omega))$ as an $\mcO_{U_T^{(1)}}$-submodule  of 
$F_{U/T*} (\mcB_{U/T})$.
Now, consider the natural  exact sequence
 \begin{align} \label{Www6}
 0 \migi F_{U/T*}(\Omega_{U/T}) / F_{U/T*}(\mr{Ker} (\nabla_\Omega)) \migi
F_{U/T*}(\Omega_{U/T}^{\otimes 2}) \migi F_{U/T*}(\mr{Coker} (\nabla_\Omega)) \migi 0,
 \end{align}
 where
 the second arrow arises from $\nabla_\Omega : \Omega_{U/T} \migi \Omega_{U/T}^{\otimes 2}$.
  By  ~\cite{Wak5}, Proposition 6.8.3, $F_{U/T*}(\mr{Coker} (\nabla_\Omega)) $ turns out to be  a vector bundle on $U_T^{(1)}$.
 Since $F_{U/T*}(\Omega_{U/T}^{\otimes 2})$ is  a vector bundle,
the exactness of  (\ref{Www6}) implies that  the quotient $F_{U/T*}(\Omega_{U/T}) / F_{U/T*}(\mr{Ker} (\nabla_\Omega)) $  is
 a vector bundle.
 Moreover, let us consider the short exact sequence
 \begin{align}
 0 \migi F_{U/T*}(\mcB_{U/T})/F_{U/T*}(\mr{Ker} (\nabla_\Omega))   \migi F_{U/T*}(\Omega_{U/T}) / F_{U/T*}(\mr{Ker} (\nabla_\Omega))
 \\ \migi
  \Omega_{U_T^{(1)}/T} \migi 0\notag
 \end{align}
induced from (\ref{E0001}) via taking quotients by $F_{U/T*}(\mr{Ker} (\nabla_\Omega)) $.
 Since both $\Omega_{U_T^{(1)}/T}$ and $F_{U/T*}(\Omega_{U/T}) / F_{U/T*}(\mr{Ker} (\nabla_\Omega))$  are vector bundles, 
 $F_{U/T*}(\mcB_{U/T})/F_{U/T*}(\mr{Ker} (\nabla_\Omega))$ is verified to be   a vector bundle, namely,
 $F_{U/T*}(\mr{Ker} (\nabla_\Omega))$ is a line subbundle of $F_{U/T*}(\mcB_{U/T})$.
  Thus, $F_{U/T*}(\mr{Ker} (\nabla_\Omega))$ specifies  a Tango structure on $U/T$.
  This completes the proof of the claim.

 Next, let $\mcL$ ($\subseteq F_{X/S*} (\mcB_{T/T})$) be a Tango structure on $U/T$.
 Denote by $\xi_{\mcL*} (\nabla^\mr{can}_\mcL)$ the $T$-connection on $\Omega_{U/T}$ corresponding to $\nabla^\mr{can}_\mcL$ via $\xi_\mcL$, which has vanishing $p$-curvature.
By the equivalence of categories displayed in (\ref{ED09ff8}), $F_{U/T*} (\mr{Ker} (\xi_{\mcL*} (\nabla^\mr{can}_\mcL)))$ ($\subseteq F_{U/T*} (\Omega_{U/T})$) turns out  to coincide with the image of 
the composite $\mcL \migiincl F_{U/T*} (\mcB_{U/T}) \migiincl F_{U/T*} (\Omega_{U/T})$.
 But, since (\ref{E0001}) is exact,   the composite  
 \begin{align}
 F_{U/T*} (\mr{Ker} (\xi_{\mcL*} (\nabla^\mr{can}_\mcL))) \migiincl  F_{U/T*} (\Omega_{U/T}) \xrightarrow{C_{U/T}} \Omega_{U^{(1)}_T/T}
 \end{align}
  is the zero map.
 Hence, $\xi_{\mcL*} (\nabla^\mr{can}_\mcL)$ specifies a pre-Tango structure on $U/T$.
 One verifies that  the assignment $\mcL \mapsto \xi_{\mcL*} (\nabla^\mr{can}_\mcL)$ determines the  inverse to
  the map $\nabla_\Omega \mapsto F_{U/T*}(\mr{Ker} (\nabla_\Omega))$ discussed above.
Consequently, we obtain  the desired bijection.
\end{proof}
\vspace{3mm}

\begin{exa} \label{E0140}
\leavevmode\\
 \ \ \ We shall  consider pre-Tango structures on an (ordinary) elliptic curve.
Let $X$ be a geometrically connected, proper, and  smooth curve over $k$ of genus $1$. 
Suppose that $X$ is ordinary, i.e., the $p$-rank of its Jacobian is maximal. 
One may find 
 an invariant differential $\delta \in \Gamma (X, \Omega_{X/k})$
 with $C_{X/k} (\delta) = \delta$, which induces  an identification $\mcO_X \isom \Omega_{X/k}$ (given by assigning $s \mapsto s \cdot \delta$ for any local section $s \in \mcO_X$).
 Now, let $\nabla$ be a $k$-connection on $\Omega_{X/k}$ with vanishing $p$-curvature.
 By means of the above identification,
 $\nabla$ may be considered as a $k$-connection on $\mcO_X$, and hence, expressed as
 $\nabla = d + u \cdot \delta$ (for some $u \in k = \Gamma (X, \mcO_X)$).
 It follows from ~\cite{K2}, Proposition 7.1.2, that $\psi^{(\mcO_X, \nabla)} =0$ implies the equality
 $(\mr{id}_X \times F_{\mr{Spec}(k)})^*(u \cdot \delta) = C_{X/k} (u \cdot  \delta)$ ($= u^{\frac{1}{p}} \cdot  C_{X/k} (\delta)$).
 Hence, we have $u = u^{\frac{1}{p}}$, or equivalently, $u \in \mbF_p$.

Here, we shall suppose that $u \neq 0$.
 If $h$ is a local section of $\mr{Ker} (\nabla)$ ($\subseteq \mcO_X$), i.e., a local function satisfying  the equality $dh = -u h \delta$, 
 then $C_{X/k} (h \cdot \delta) = - u^{-\frac{1}{p}} \cdot C_{X/k} (dh) = - u^{-\frac{1}{p}} \cdot 0 = 0$.
 This implies that $F_{X/k*}(\mr{Ker}(\nabla)) \subseteq \mr{Ker} (C_{X/k})$, that is to say, $\nabla$ forms
a pre-Tango structure on $X/k$.
Next, suppose that $u =0$, i.e., $\nabla = d$.
Then, $F_{X/k*}(\mr{Ker}(\nabla))$ coincides with $\mcO_{X^{(1)}_k}\cdot \delta$.
For any nonzero  local section $v$ of $\mcO_X$, $C_{X/k} (v^p \cdot  \delta) = v \cdot C_{X/k} (\delta) = v\delta \neq 0$.
Hence,  $F_{X/k*}(\mr{Ker}(\nabla)) \nsubseteq \mr{Ker} (C_{X/k})$, and $\nabla$ is not a pre-Tango structure.

Consequently,   we have obtained the fact that  {\it the set of pre-Tango structures on an ordinary elliptic curve $X/k$ is in bijection with
the set of $k$-connections $\nabla$ on $\mcO_X$ with vanishing $p$-curvature which is not equal to the universal derivation $d$.}
In particular, the number of pre-Tango structures on $X/k$ is exactly $p-1$ ($= \sharp (\mbF_p \setminus \{ 0 \})$).
\end{exa}

\vspace{5mm}
\subsection{Pre-Tango structures vs. dormant Miura $\mfs \mfl_2$-opers} \label{ScSc3}
\leavevmode\\ 
\vspace{-4mm}

Denote by
\begin{align}
\overline{\mfT} \mfa \mfn_{g,r}
\end{align}
the set-valued contravariant  functor on $\mfS \mfc \mfh_{/\overline{\mfM}_{g,r}}$ which, 
to any object $S \migi \overline{\mfM}_{g,r}$ classifying a pointed stable curve 
  $\mfX$,
assigns the set of pre-Tango structures on  $\mfX$.
Also, for each $\vec{\varepsilon} \in k^{\times r}$ (where we take $\vec{\varepsilon} := \emptyset$ if $r =0$), 
we shall write
\begin{align} \label{GGk}
\overline{\mfT} \mfa \mfn_{g,r, \vec{\varepsilon}}
\end{align}
for the subfunctor of $\overline{\mfT} \mfa \mfn_{g,r} $ classifying pre-Tango structures of monodromies  $\vec{\varepsilon}$.
$\overline{\mfT} \mfa \mfn_{ g,r}$ and $\overline{\mfT} \mfa \mfn_{g,r, \vec{\varepsilon}}$
may be represented by closed substacks of $\overline{\mfC} \mfo_{g,r}$ and $\overline{\mfC} \mfo_{g,r, \vec{\varepsilon}}$ respectively.
One verifies (from an argument similar to the argument in the proof of the  non-resp'd  assertion in Proposition \ref{pr0011}) that  $\overline{\mfT} \mfa \mfn_{g,r, \vec{\varepsilon}}$ is empty unless $\vec{\varepsilon}$ lies $\mbF_p^{\times r}$ (or $\vec{\varepsilon} = \emptyset$),  and  $\overline{\mfT} \mfa \mfn_{g,r}$ decomposes into the disjoint union
\begin{align} \label{GG2}
\overline{\mfT} \mfa \mfn_{g,r} = \coprod_{\vec{\varepsilon} \in \mbF_p^{\times r}} \overline{\mfT} \mfa \mfn_{g,r, \vec{\varepsilon}}.
\end{align}
Moreover,  the following Theorem  \ref{Td01} holds; by this proposition, 
{\it the notion of a dormant generic Miura $\mfg$-oper may be thought of as a generalization of the notion of a Tango structure}.

\vspace{3mm}
\bt[Theorem \ref{ThD}, (i)] \label{Td01}
 \leavevmode\\
 \ \ \ 
 Let $\vec{\varepsilon} := (\varepsilon_i)_{i=1}^r \in k^{\times r}$,
  and 
 write $- \vec{\varepsilon} := (-\varepsilon_i)_{i=1}^r$ (where $\vec{\varepsilon}  := \emptyset$ and $- \vec{\varepsilon} := \emptyset$ if $r =0$).
 Then, the composite isomorphism 
 \begin{align}  \label{EEE12}
 \overline{\mfC} \mfo_{g,r, - \vec{\varepsilon}} \xrightarrow[\sim]{(\ref{eq071})} \mfC \overline{\mfC} \mfo_{\mfs \mfl_2, g,r,  [\vec{\varepsilon}\,]} \xrightarrow[\sim]{\Xi_{\mfs \mfl_2, g,r, [\vec{\varepsilon}], p^\circ_{-1}}} \mfM \overline{\mfO} \mfp_{\mfs \mfl_2, g,r, [\vec{\varepsilon}\,]}
 \end{align}
  restricts  to an isomorphism
\begin{align} \label{E0010}
\overline{\mfT} \mfa \mfn_{g,r, -\vec{\varepsilon}} \isom \mfM \overline{\mfO} \mfp_{\mfs \mfl_2, g,r, [\vec{\varepsilon}\,]}^{^{\mr{Zzz...}}}
\end{align}
 over $\overline{\mfM}_{g,r}$.
\et
\begin{proof}
Let  $S$ be a $k$-scheme and  $\mfX := (X/S, \{ \sigma_i \}_{i=1}^r)$ an $r$-pointed stable curve over $S$ of genus $g$.
Also, let $\nabla$ be an $S^\mr{log}$-connection on $\Omega_{X^\mr{log}/S^\mr{log}}$ classified by $\overline{\mfC} \mfo_{g,r, - \vec{\varepsilon}}$.
 Denote by $\widehat{\mcF}^{\heartsuit \diamondsuit} := (\mcF_{\mcO_X}^{[2]\dagger}, \nabla_\mcF, \mpf^{\dagger}_{\mcO_X}, \mpf^{-, \dagger}_{\mcO_X})$
 the $\mcO_X^\flat$-special Miura $\mr{GL}_2$-oper on $\mfX$ determined by $\nabla$ via the composite of (\ref{EEE12}) and the inverse of  the right-hand vertical arrow in  (\ref{efe0991}).
In order to complete the proof,
it suffices to prove the claim that $\nabla$ specifies a pre-Tango structure on $\mfX$ if and only if 
$\widehat{\mcF}^{\heartsuit \diamondsuit}$ is dormant.

To begin with,
denote by $U$ the smooth locus of  $X \setminus \mr{Supp} (D_\mfX)$ relative to  $S$, where $D_\mfX$ denotes the 
effective relative divisor  on $X$ defined to be the union of the image of the marked points $\sigma_i$ ($i =1, \cdots, r$).
We equip $U$ with the  log  structure  pulled-back from $X^\mr{log}$ via the open immersion $U\migiincl X$.
Denote by $U^\mr{log}$ the resulting log scheme.
Since the natural projection $U^\mr{log} \migi S^\mr{log}$ is strict (cf. ~\cite{ILL2}, \S\,1.2), 
 each  pre-Tango structure  (resp.,  Miura $\mr{GL}_2$-oper) on $U^\mr{log}/S^\mr{log}$ may be identified with a pre-Tango structure (resp., Miura $\mr{GL}_2$-oper) on $U/S$.

Now, suppose that $\nabla$ specifies a pre-Tango structure.
By Proposition \ref{Pff110}, the restriction $\nabla |_U$ corresponds (via  (\ref{EE09ff8})) to a
Tango structure $\mcL$ ($\subseteq F_{U/S*}(\mcB_{U/S})$) on $U/S$. 
Moreover, by  the last assertion of Proposition \ref{PP020},  the restriction  $\widehat{\mcF}^{\heartsuit \diamondsuit} |_U$ of $\widehat{\mcF}^{\heartsuit \diamondsuit}$ to $U$  is isomorphic to the dormant  Miura $\mr{GL}_2$-oper  $\widehat{\mcT} \hspace{-1mm} an_\mcL^{\heartsuit \diamondsuit} := (\mcF^{[2]\dagger}_{\mcO_U}, \nabla_\mcL, \mpf^\dagger_{\mcO_U}, \mpf^{-, \dagger}_{\mcO_U})$. 
Hence, since 
 $U$ is scheme-theoretically dense in $X$,  $\widehat{\mcF}^{\heartsuit \diamondsuit}$ itself  turns out to be dormant, as desired.

Conversely, suppose  that $\widehat{\mcF}^{\heartsuit \diamondsuit}$ is dormant.
The dormant Miura $\mr{GL}_2$-oper on $U/S$ defined as the restriction $\widehat{\mcF}^{\heartsuit \diamondsuit} |_U$ of $\widehat{\mcF}^{\heartsuit \diamondsuit}$ to $U$
comes from a Tango-structure $\mcL$ (cf. Proposition \ref{PP020}), i.e., isomorphic to $\widehat{\mcT} \hspace{-1mm} an_\mcL^{\heartsuit \diamondsuit}$.
By the various definitions involved, the pre-Tango structure on $U/S$ corresponding to $\mcL$ via  (\ref{EE09ff8}) 
 coincides with $\nabla |_U$ (cf. Proposition \ref{Pff110}).
In particular, both the $p$-curvature of $\nabla$ and the composite
\begin{align} \label{EEE2234}
F_{X/S*} (\mr{Ker} (\nabla)) \migi F_{X/S*} (\Omega_{X^\mr{log}/S^\mr{log}}) \xrightarrow{C_{X^\mr{log}/S^\mr{log}}} \Omega_{X^\mr{log}/S^\mr{log}}
\end{align}
vanishes  on $U$.
But, since $U$ is scheme-theoretically dense in $X$, both  the $p$-curvature of $\nabla$ and the composite
(\ref{EEE2234})   vanishe identically on $X$.
That is to say, $\nabla$ specifies a pre-Tango structure on $\mfX$.
This completes the proof of Theorem \ref{Td01}.
\end{proof}

\vspace{10mm}
\section{Deformations of (dormant) Miura opers} \vspace{3mm}

In this section, we  describe the deformation space of a given (dormant) generic Miura $\mfg$-oper in terms of de Rham cohomology of complexes (cf. Propositions \ref{P0f7} and \ref{P0337}). 
By applying these descriptions, one may prove (cf. Theorem  \ref{T01}) that the moduli stack of dormant generic Miura $\mfs \mfl_2$-opers  is  smooth.

\vspace{5mm}
\subsection{The tangent bundle of $\mfM \overline{\mfO} \mfp_{\mfg, g,r}$.}
\leavevmode\\ \vspace{-4mm}

First, assume that either one of the  three  conditions (Char)$_{0}$,  $(\mr{Char})_{p}$,  and $(\text{Char})_{p}^{\mfs \mfl}$ (cf. \S\,\ref{SSQ1}) is  satisfied.
The  moduli stack $\mfM \overline{\mfO} \mfp_{\mfg,   g,r}$
  admits a  log structure
 pulled-back from the log structure of $\overline{\mfM}_{g,r}^\mr{log}$;
we denote the resulting log stack by
\begin{align}\label{univ1p}
\mfM \overline{\mfO} \mfp_{\mfg, g,r}^\mr{log}.
\end{align}
The forgetting morphism  $\mfM \overline{\mfO} \mfp_{\mfg, g,r} \migi \overline{\mfM}_{g,r}$ extends to a morphism
$\mfM \overline{\mfO} \mfp_{\mfg, g,r}^\mr{log} \migi \overline{\mfM}_{g,r}^\mr{log}$.

Let $S$ be a scheme over $k$ and 
$\mfX := (f : X \migi S, \{ \sigma_i \}_{i=1}^r)$  an $r$-pointed stable curve of genus $g$ over $S$.
If  $\nabla : \mcK^0 \migi \mcK^1$ is a morphism of sheaves of abelian groups on $X$, then 
it may be thought of as a complex concentrated in  degrees $0$ and $1$; we denote this complex by 
\begin{equation}\mcK^\bullet[\nabla]
\end{equation}
(where $\mcK^i [\nabla] := \mcK^i$ for $i =0$, $1$).
Also, for $i = 0, 1, \cdots$, we obtain the sheaf
\begin{equation}  \mbR^i f_*(\mcK^\bullet [\nabla])  \end{equation}
on $S$, where $ \mbR^i f_*(-)$ is the $i$-th hyper-derived functor of $\mbR^0f_*(-)$ (cf. ~\cite{K}, (2.0)).
In particular, $\mbR^0f_*(\mcK^\bullet [\nabla]) = f_*(\mr{Ker}(\nabla))$.

Now, let
$\widehat{\mcE}^{\spadesuit \diamondsuit} := (\mcE_\mbB^\dagger, \nabla_\mcE, \mcE'^\dagger_\mbB, \eta_\mcE^\dagger)$ be a $p_{-1}$-special Miura $\mfg$-oper on $\mfX$.
We shall write $(\mcE'_\mbB, \nabla_{\mcE'_\mbB})$ for the log flat $\mbB$-torsor associated with $\widehat{\mcE}^{\spadesuit \diamondsuit}$ and 
\begin{align}
\nabla_{\mcE'_\mbB}^\mr{ad} : \mfb_{\mcE'^\dagger_\mbB} \migi \Omega_{X^\mr{log}/S^\mr{log}}\otimes \mfb_{\mcE'^\dagger_\mbB}
\end{align}
for the $S^\mr{log}$-connection on  the vector bundle $\mfb_{\mcE'^\dagger_\mbB}$ induced from $\nabla_{\mcE'_\mbB}$ via the adjoint representation $\mbB \migi \mr{GL} (\mfb)$.
If we write $\mfg_j := \mfg^j / \mfg^{j+1}$ (for each $j \in \mbZ$),
then it has a $\mbT$-action induced from  the $\mbT$-action on $\mfg^j$, and hence, 
we have an $\mcO_X$-module $(\mfg_j)_{\mcE^\dagger_\mbT}$.
The decomposition displayed in  (\ref{22224})
 gives rise to  a  decomposition 
\begin{align}
\mfg_{\mcE^\dagger_{\mbT}} \isom \bigoplus_{j \in \mbZ} (\mfg_j)_{\mcE^\dagger_\mbT}
\end{align}
(which restricts to a decomposition $\mfb_{\mcE'^\dagger_\mbT} \isom \bigoplus_{j \leq 0} (\mfg_j)_{\mcE^\dagger_\mbT}$)
on $\mfg_{\mcE^\dagger_\mbT}$.
We shall write
\begin{align}
\mfg^{-1/1}_{\mcE^\dagger_\mbT} := (\mfg_{-1})_{\mcE^\dagger_\mbT} \oplus (\mfg_{0})_{\mcE^\dagger_\mbT} \ \left(= \mfg_{\mcE_\mbT^\dagger}^{-1} / \mfg^1_{\mcE_\mbT^\dagger} \right)
\end{align}
and regard it
 as an  $\mcO_X$-submodule of $\mfb_{\mcE'^\dagger_\mbT}$.
Consider the $f^{-1} (\mcO_S)$-linear morphism
\begin{align}
\widetilde{\nabla}_{\mcE'_\mbB}^{\mr{ad}} : \widetilde{\mcT}_{\mcE'^{\dagger\mr{log}}_\mbB/S^\mr{log}} \migi \Omega_{X^\mr{log}/S^\mr{log}} \otimes \widetilde{\mcT}_{\mcE'^{\dagger \mr{log}}_\mbB/S^\mr{log}} 
\end{align}
determined uniquely by the condition that
\begin{align}
\langle \partial, \widetilde{\nabla}^{\mr{ad}}_{\mcE'_\mbB} (s) \rangle = [\nabla_\mcE (\partial), s] - \nabla_\mcE ([\partial, \mfa^\mr{log}_{\mcE'^\dagger_\mbB} (s)]),
\end{align}
where
\begin{itemize}
\item[$\bullet$]
$s$ and $\partial$ denote arbitrary  local sections of $\widetilde{\mcT}_{\mcE'^{\dagger \mr{log}}_{\mbB}/S^\mr{log}}$ and $\mcT_{X^\mr{log}/S^\mr{log}}$ respectively;
\item[$\bullet$]
$\langle -, - \rangle$ denotes the $\mcO_X$-bilinear pairing
\begin{align}
\mcT_{X^\mr{log}/S^\mr{log}} \times (\Omega_{X^\mr{log}/S^\mr{log}} \otimes \widetilde{\mcT}_{\mcE'^{\dagger \mr{log}}_{\mbB}}) \migi \widetilde{\mcT}_{\mcE'^{\dagger \mr{log}}_{\mbB}/S^\mr{log}}
\end{align}
induced by the natural pairing $\mcT_{X^\mr{log}/S^\mr{log}}  \times \Omega_{X^\mr{log}/S^\mr{log}}\migi \mcO_X$;
\item[$\bullet$]
$[-, -]$'s denote the Lie bracket operators in the respective tangent bundles.
\end{itemize}
One verifies that the restriction of $\widetilde{\nabla}^{\mr{ad}}_{\mcE'_\mbB}$ to $\mfb_{\mcE'^{\dagger}_\mbB}$ ($\subseteq \widetilde{\mcT}_{\mcE'^{\dagger \mr{log}}_{\mbB}/S^\mr{log}}$) coincides with $\nabla_{\mcE'_\mbB}^\mr{ad}$.
Moreover,  the following assertion holds.

\vspace{3mm}
\ble \label{LEE0} \leavevmode\\
 \ \ \ 
The image of $\widetilde{\nabla}^{\mr{ad}}_{\mcE'_\mbB}$ is contained in $\Omega_{X^\mr{log}/S^\mr{log}} \otimes \mfg_{\mcE^\dagger_{\mbT,}}^{-1/1}$.
\ele
\begin{proof}
For any local sections $s$ and  $\partial$  of $\widetilde{\mcT}_{\mcE'^{\dagger \mr{log}}_\mbB/S^\mr{log}}$ and  $\mcT_{X^\mr{log}/S^\mr{log}}$ respectively,    the following sequence of equalities holds:
\begin{align}
& \  \langle \partial, (\mr{id}_{\Omega_{X^\mr{log}/S^\mr{log}}} \otimes \mfa^\mr{log}_{\mcE'^\dagger_\mbB})( \widetilde{\nabla}^\mr{ad}_{\mcE'_\mbB} (s)) \rangle \\
= & \ \mfa^\mr{log}_{\mcE'^\dagger_\mbB} (\langle \partial, \widetilde{\nabla}^\mr{ad}_{\mcE'_\mbB} (s)\rangle) \notag  \\
=  & \ \mfa^\mr{log}_{\mcE'^\dagger_\mbB} ([\nabla_\mcE (\partial), s]) - (\mfa^\mr{log}_{\mcE'^\dagger_\mbB} \circ \nabla_\mcE) ([\partial, \mfa^\mr{log}_{\mcE'^\dagger_\mbB} (s)]) \notag \\
 = & \  [(\mfa^\mr{log}_{\mcE'^\dagger_\mbB} \circ \nabla_\mcE)(\partial), \mfa^\mr{log}_{\mcE'^\dagger_\mbB} (s)] - \mr{id}_{\mcT_{X^\mr{log}/S^\mr{log}}} ([\partial, \mfa^\mr{log}_{\mcE'^\dagger_\mbB} (s)]) \notag \\
 = & \ [\partial, \mfa^\mr{log}_{\mcE'^\dagger_\mbB} (s)] - [\partial, \mfa^\mr{log}_{\mcE'^\dagger_\mbB} (s)] \notag \\
 = & \ 0. \notag
\end{align}
This implies that the image  of $\widetilde{\nabla}^{\mr{ad}}_{\mcE'_\mbB}$ is contained in $\Omega_{X^\mr{log}/S^\mr{log}} \otimes \mfg_{\mcE^\dagger_\mbT}$ ($= \mr{Ker} (\mr{id}_{\Omega_{X^\mr{log}/S^\mr{log}}}$ $\otimes \mfa^\mr{log}_{\mcE'^\dagger_\mbB})$).
Moreover, by the definition of a $\mfg$-oper,
the image of $\widetilde{\nabla}^{\mr{ad}}_{\mcE'_\mbB}$ is contained in 
$\Omega_{X^\mr{log}/S^\mr{log}} \otimes \widetilde{\mcT}^{-1}_{\mcE_\mbG^{\dagger\mr{log}}/S^\mr{log}}$, and hence,  in $\Omega_{X^\mr{log}/S^\mr{log}} \otimes \mfg_{\mcE^\dagger_{\mbT}}^{-1/1}$ ($= \Omega_{X^\mr{log}/S^\mr{log}} \otimes (\widetilde{\mcT}_{\mcE'^{\dagger \mr{log}}_{\mbB}/S^\mr{log}} \cap \mfg_{\mcE^\dagger_\mbT} \cap \widetilde{\mcT}^{-1}_{\mcE_\mbG^{\dagger\mr{log}}/S^\mr{log}})$).
This completes the proof of Lemma \ref{LEE0}.
\end{proof}
\vspace{3mm}

Because of  the above lemma, $\widetilde{\nabla}^{\mr{ad}}_{\mcE'_\mbB}$ restricts to an $f^{-1}(\mcO_S)$-linear morphism
\begin{align}
\widetilde{\nabla}^{\mr{ad}, -1/1}_{\mcE_\mbT} : \widetilde{\mcT}_{\mcE^{\dagger \mr{log}}_{\mbT}/S^\mr{log}} \migi  \Omega_{X^\mr{log}/S^\mr{log}} \otimes \mfg_{\mcE^\dagger_{\mbT}}^{-1/1}.
\end{align}
The pair of the  morphism  $\mfa^\mr{log}_{\mcE^\dagger_\mbT} : \widetilde{\mcT}_{\mcE^{\dagger \mr{log}}_{\mbT}/S^\mr{log}} \migi \mcT_{X^\mr{log}/S^\mr{log}}$ and the zero map $\Omega_{X^\mr{log}/S^\mr{log}} \otimes \mfg_{\mcE^\dagger_{\mbT}}^{-1/1} \migi 0$  specifies  a morphism
\begin{align} \label{FGH}
\mcK^\bullet [\widetilde{\nabla}^{\mr{ad}, -1/1}_{\mcE_\mbT}] \migi \mcT_{X^\mr{log}/S^\mr{log}} [0]
\end{align}
of complexes, where for any abelian sheaf $\mcF$ we shall write $\mcF [0]$ for $\mcF$ considered as 
a complex concentrated at degree $0$.

\vspace{3mm}
\bpr \label{P0f7} \leavevmode\\
 \ \ \ 
Let $c_{\mfX, \widehat{\mcE}^{\spadesuit\diamondsuit}} : S \migi \mfM \overline{\mfO} \mfp_{\mfg,  g, r}$  (resp., $c_{\mfX} : S \migi \overline{\mfM}_{g,r}$) be the  $S$-rational point of $\mfM \overline{\mfO} \mfp_{\mfg, g, r}$ (resp., $\overline{\mfM}_{g,r}$) classifying the pair $(\mfX, \widehat{\mcE}^{\spadesuit\diamondsuit})$ (resp., $\mfX$).
 Then, there exists a natural  isomorphism
 \begin{align} \label{FGH2}
c_{\mfX, \widehat{\mcE}^{\spadesuit\diamondsuit}}^*(\mcT_{\mfM \overline{\mfO} \mfp^\mr{log}_{\mfg, g, r}/k}) \isom  \mbR^1 f_* (\mcK^\bullet [\widetilde{\nabla}^{\mr{ad},  -1/1}_{\mcE_\mbT}])
\end{align}
  of $\mcO_S$-modules, which makes the square diagram
  \begin{align}
  \vcenter{\xymatrix{
c_{\mfX, \widehat{\mcE}^{\spadesuit\diamondsuit}}^*(\mcT_{\mfM \overline{\mfO} \mfp^\mr{log}_{\mfg, g, r}/k})  \ar[r]^-{(\ref{FGH2})}_-{\sim} \ar[d] & \mbR^1 f_* (\mcK^\bullet [\widetilde{\nabla}^{\mr{ad},  -1/1}_{\mcE_\mbT}])  \ar[d] \\
c_{\mfX}^*(\mcT_{\overline{\mfM}^\mr{log}_{g, r}/k})   \ar[r]^{\sim} & \mbR^1f_* (\mcT_{X^\mr{log}/S^\mr{log}})
  }}
  \end{align}
  commute, where the right-hand vertical arrow denotes the morphism obtained by applying the functor $\mbR^1 f_* (-)$ to  (\ref{FGH}), the left-hand vertical arrow denotes the morphism arising from the forgetting morphism  $\mfM \overline{\mfO}\mfp_{\mfg, g,r}^{\mr{log}} \migi \overline{\mfM}_{g,r}^{\mr{log}}$, and the lower horizontal arrow denotes the isomorphism defined as  the Kodaira-Spencer map of $\mfX$.
  \epr
\begin{proof}
The assertion follows from 
an argument (in the case where $\mfX$ is a pointed stable curve over an  arbitrary scheme $S$) similar to the argument (in the case where $\mfX$ is an unpointed smooth curve over $\mbC$) given in ~\cite{Chen}.
Indeed, by the explicit description of the hyperchomology sheaf $ \mbR^1f_*(\mcK^\bullet [\widetilde{\nabla}_{\mcE_\mbT}^{\mr{ad}, -1/1}])$ in terms of the \v{C}ech double complex associated with  $\mcK^\bullet [\widetilde{\nabla}_{\mcE_\mbT}^{\mr{ad}, -1/1}]$, one verifies that $ \mbR^1f_*(\mcK^\bullet [\widetilde{\nabla}_{\mcE_\mbT}^{\mr{ad}, -1/1}])$ may be naturally identified with the deformation space (relative to $S$) of $(\mfX, \widehat{\mcE}^{\spadesuit\diamondsuit})$, i.e., 
the pair consisting of the pointed stable curve $\mfX$ and  (the isomorphism class of)  the 
generic Miura $\mfg$-oper $\widehat{\mcE}^{\spadesuit\diamondsuit}$ on it.
That is to say, $ \mbR^1f_*(\mcK^\bullet [\widetilde{\nabla}_{\mcE_\mbT}^{\mr{ad}, -1/1}])$ is canonically isomorphic to the $\mcO_S$-module 
$c^{*}_{\mfX, \widehat{\mcE}^{\spadesuit \diamondsuit}}(\mcT_{\mfM \overline{\mfO} \mfp^\mr{log}_{\mfg, \hslash, g,r}/k})$.
In particular, for completing the proof of Proposition \ref{P0f7}, we refer to ~\cite{Chen}, Propositions 4.1.3 and  4.3.1.
\end{proof}

\vspace{5mm}
\subsection{The tangent bundle of $\mfM \overline{\mfO} \mfp^{^\mr{Zzz...}}_{\mfg, g,r}$.}
\leavevmode\\ \vspace{-4mm}

Next, suppose further that $\mr{char} (k) = p >0$, i.e.,   either one of the two  conditions $(\mr{Char})_{p}$,  and $(\text{Char})_{p}^{\mfs \mfl}$  is  satisfied.
Let us consider the natural composite
\begin{align}
\omega : \Omega_{X^\mr{log}/S^\mr{log}} \otimes \mfg_{\mcE^\dagger_{\mbT}}^{-1/1} \migiincl  \Omega_{X^\mr{log}/S^\mr{log}} \otimes \mfb_{\mcE'^\dagger_{\mbT}} \migisurj
\mr{Coker} (\nabla_{\mcE'_\mbB}^\mr{ad}).
\end{align}
By Lemma  \ref{L07gg}  below, 
$\widetilde{\nabla}^{\mr{ad}, -1/1}_{\mcE^{\dagger}_\mbT}$ induces (by restricting its codomain)
an $f^{-1} (\mcO_S)$-linear morphism 
\begin{align}
\widetilde{\nabla}^{\mr{ad},   \omega}_{\mcE_\mbT} : \widetilde{\mcT}_{\mcE^{\dagger \mr{log}}_{\mbT}/S^\mr{log}} \migi  \mr{Ker} (\omega).
\end{align}

\vspace{3mm}
\ble \label{L07gg} \leavevmode\\
 \ \ \ 
The image of $\widetilde{\nabla}^{\mr{ad}, -1/1}_{\mcE_\mbT} $ is contained in $\mr{Ker} (\omega)$.
 \ele
\begin{proof}
Let us consider the $\mcO_X$-linear endomorphism $\zeta$ of
$\widetilde{\mcT}_{\mcE'^{\dagger \mr{log}}_\mbB/S^\mr{log}}$ determined by assigning
$s \mapsto s - \nabla_\mcE \circ \mfa^\mr{log}_{\mcE'^\dagger_\mbB} (s)$ for any local section $s \in \widetilde{\mcT}_{\mcE'^{\dagger \mr{log}}_\mbB/S^\mr{log}}$.
Then,
  for any local section $s$ of $\widetilde{\mcT}_{\mcE'^{\dagger \mr{log}}_\mbB/S^\mr{log}}$, the following equalities hold:
\begin{align}
 \mfa^\mr{log}_{\mcE'^\dagger_\mbB}(\zeta (s))  
& = \mfa^\mr{log}_{\mcE'^\dagger_\mbB} (s - \nabla_\mcE \circ \mfa^\mr{log}_{\mcE'^\dagger_\mbB} (s))  \\
& = \mfa^\mr{log}_{\mcE'^\dagger_\mbB} (s) - (\mfa^\mr{log}_{\mcE'^\dagger_\mbB} \circ \nabla_\mcE) (\mfa^\mr{log}_{\mcE'^\dagger_\mbB} (s)) \notag  \\
& = \mfa^\mr{log}_{\mcE'^\dagger_\mbB} (s) -\mfa^\mr{log}_{\mcE'^\dagger_\mbB} (s) \notag \\
& = 0. \notag 
\end{align}
Hence, we have  $\mr{Im}(\zeta) \subseteq \mr{Ker} (\mfa^\mr{log}_{\mcE'^\dagger_\mbB})$ ($= \mfb_{\mcE'^\dagger_\mbB}$).
Moreover, let us observe the following equalities:
\begin{align} \label{F3456}
& \ \langle \partial, \nabla_{\mcE'_\mbB}^\mr{ad} (\zeta (s))\rangle \\
= & \ 
\langle \partial, \widetilde{\nabla}_{\mcE'_\mbB}^{\mr{ad}, -1/1} (\zeta (s))\rangle \notag  \\
= &  \ [\nabla_\mcE (\partial), \zeta (s)] - \nabla_\mcE ([\partial, \mfa^\mr{log}_{\mcE'^\dagger_\mbB} (\zeta (s))]) \notag \\
= & \ [\nabla_\mcE (\partial), s] - [\nabla_\mcE (\partial), \nabla_\mcE \circ \mfa^\mr{log}_{\mcE'^\dagger_\mbB} (s)] - \nabla_\mcE ([\partial, \mfa^\mr{log}_{\mcE'^\dagger_\mbB} (s)]) \notag \\ 
&   +  \nabla_\mcE ([\partial, \mfa^\mr{log}_{\mcE'^\dagger_\mbB} (\nabla_\mcE \circ \mfa^\mr{log}_{\mcE'^\dagger_\mbB} (s))]) \notag  \\
= & \ [\nabla_\mcE (\partial), s] - \nabla_\mcE ([\partial,  \mfa^\mr{log}_{\mcE'^\dagger_\mbB} (s)]) - \nabla_\mcE ([\partial, \mfa^\mr{log}_{\mcE'^\dagger_\mbB} (s)]) +  \nabla_\mcE ([\partial,  \mfa^\mr{log}_{\mcE'^\dagger_\mbB} (s)]) \notag \\
= & \ \langle \partial, \widetilde{\nabla}_{\mcE'_\mbB}^\mr{ad} (s)\rangle. \notag 
\end{align}
This implies that the image of $\widetilde{\nabla}_{\mcE'_\mbB}^\mr{ad}$ is contained in the image of $\nabla_{\mcE'_\mbB}^\mr{ad}$, which coincides tautologically  with the kernel of the quotient $\Omega_{X^\mr{log}/S^\mr{log}} \otimes \mfb_{\mcE'^\dagger_{\mbT}} \migisurj
\mr{Coker} (\nabla_{\mcE'_\mbB}^\mr{ad})$.
Thus, this  observation and the definition of 
$\widetilde{\nabla}_{\mcE_\mbT}^{\mr{ad}, -1/1}$ implies  the validity  of Lemma \ref{L07gg}.
\end{proof}
\vspace{3mm}

Here, we shall denote by 
\begin{align} \label{FGH7}
\mfM \overline{\mfO} \mfp_{\mfg, g,r}^{^\mr{Zzz...} \mr{log}}
\end{align}
 the log stack defined to be 
the stack $\mfM \overline{\mfO} \mfp_{\mfg, g,r}^{^\mr{Zzz...}}$ equipped with the log structure pulled-back from the log structure of $\overline{\mfM}_{g,r}^\mr{log}$ via the forgetting morphism  $\mfM \overline{\mfO} \mfp_{\mfg, g,r}^{^\mr{Zzz...}} \migi \overline{\mfM}_{g,r}$.

\vspace{3mm}
\bpr \label{P0337} \leavevmode\\
 \ \ \ 
Let $c_{\mfX, \widehat{\mcE}^{\spadesuit \diamondsuit}}$ and $c_\mfX$ be as in Proposition \ref{P0f7}.
Suppose further  that
$\widehat{\mcE}^{\spadesuit \diamondsuit}$ is dormant, i.e.,  $c_{\mfX, \widehat{\mcE}^{\spadesuit \diamondsuit}}$ factors through the closed immersion $\mfM \overline{\mfO} \mfp_{\mfg, g,r}^{^\mr{Zzz...}} \migi \mfM \overline{\mfO} \mfp_{\mfg, g,r}$.
Denote by $\overline{c}_{\mfX, \widehat{\mcE}^{\spadesuit \diamondsuit}} : S \migi \mfM \overline{\mfO} \mfp_{\mfg, g,r}^{^\mr{Zzz...}}$ the resulting $S$-rational point of $\mfM \overline{\mfO} \mfp_{\mfg, g,r}^{^\mr{Zzz...}}$.
 Then, there exists a natural  isomorphism
 \begin{align} \label{FGH4}
 \overline{c}_{\mfX, \widehat{\mcE}^{\spadesuit \diamondsuit}}^* (\mcT_{\mfM \overline{\mfO} \mfp^{^{\mr{Zzz...}}\mr{log}}_{\mfg, g, r}/k}) \isom \mbR^1 f_* (\mcK^\bullet [\widetilde{\nabla}^{\mr{ad},   \omega}_{\mcE_\mbT}])
 \end{align}
which makes the square diagram
\begin{align}
\vcenter{\xymatrix{
\overline{c}_{\mfX, \widehat{\mcE}^{\spadesuit \diamondsuit}}^* (\mcT_{\mfM \overline{\mfO} \mfp^{^{\mr{Zzz...}}\mr{log}}_{\mfg, g, r}/k})
\ar[r]_{\sim}^{(\ref{FGH4})} \ar[d] &  \mbR^1 f_* (\mcK^\bullet [\widetilde{\nabla}^{\mr{ad},   \omega}_{\mcE_\mbT}])
 \ar[d]
\\
c_{\mfX, \widehat{\mcE}^{\spadesuit \diamondsuit}}^* (\mcT_{\mfM \overline{\mfO} \mfp^\mr{log}_{\mfg, g, r}/k})
\ar[r]^{(\ref{FGH2})}_{\sim} &  
\mbR^1 f_* (\mcK^\bullet [\nabla^{\mr{ad}, -1/1}_{\mcE_\mbT}]) 
}}
\end{align}
commute, where  the right-hand vertical arrow denotes the morphism obtained by applying the functor $\mbR^1 f_* (-)$ to   the natural inclusion
$\mcK^\bullet [\widetilde{\nabla}^{\mr{ad},  \omega}_{\mcE_\mbT}] \migi\mcK^\bullet [\widetilde{\nabla}^{\mr{ad}, -1/1}_{\mcE_\mbT}]$
and the left-hand vertical arrow arises from 
the closed immersion 
$\mfM \overline{\mfO} \mfp^{^\text{Zzz...}\mr{log}}_{\mfg, g,r} \migi \mfM \overline{\mfO} \mfp^\mr{log}_{\mfg, g,r}$.

  \epr
\begin{proof}
The assertion follows from ~\cite{Wak5},  Proposition 6.8.1, and the various definitions involved.
\end{proof}

\vspace{5mm}
\subsection{The smoothness of  $\mfM \overline{\mfO} \mfp^{^\mr{Zzz...}}_{\mfs \mfl_2, g,r,   \vec{\varepsilon}}$.}
\leavevmode\\ \vspace{-4mm}

In what follows, let us consider the case where $\mfg = \mfs \mfl_2$ and  prove the smoothness of the moduli stack $\mfM \overline{\mfO} \mfp^{^\mr{Zzz...}}_{\mfs \mfl_2, g,r,  \vec{\varepsilon}}$ of dormant generic Miura $\mfs \mfl_2$-opers.

\vspace{3mm}
\ble \label{L07} \leavevmode\\
 \ \ \ 
Let us keep the above notation, and suppose further that $\mfg = \mfs \mfl_2$ and   $S = \mr{Spec} (k)$.
Denote by $\mbH^j (X, \mcK^\bullet [\widetilde{\nabla}^{\mr{ad}, \omega}_{\mcE_\mbT}])$
the $j$-th  hypercohomology group of the complex  $\mcK^\bullet [\widetilde{\nabla}^{\mr{ad},  \omega}_{\mcE_\mbT}]$.
\begin{itemize}
\item[(i)]
The following  equalities   hold:
\begin{align} \label{ee320}
\mbH^0 (X, \mcK^\bullet [\widetilde{\nabla}^{\mr{ad}, \omega}_{\mcE_\mbT}]) = \mbH^2 (X, \mcK^\bullet [\widetilde{\nabla}^{\mr{ad}, \omega}_{\mcE_\mbT}]) =0.
 \end{align} 
In particular, any $1$-st order deformation of $(\mfX, \widehat{\mcE}^{\spadesuit \diamondsuit})$ (i.e.,  the pair consisting of the pointed stable curve $\mfX$ and the dormant generic Miura $\mfs \mfl_2$-oper $\widehat{\mcE}^{\spadesuit \diamondsuit}$ on it) is unobstructed.
That is to say,  $\mfM \overline{\mfO} \mfp_{\mfs \mfl_2, g,r}^{^\mr{Zzz...}}$ is smooth at the point $\overline{c}_{\mfX, \widehat{\mcE}^{\spadesuit \diamondsuit}}$.
\item[(ii)]
Let  $\vec{\varepsilon} := (\varepsilon_i)_{i=1}^r \in \mbF_p^{\times r}$ (where $\vec{\varepsilon} := \emptyset$ if $r =0$).
 Assume  further that  $X$ is smooth over $k$ and   $\widehat{\mcE}^{\spadesuit \diamondsuit}$ is of exponents $[\vec{\varepsilon}\,]$. 
Then, 
  $\mbH^1 (X, \mcK^\bullet [\widetilde{\nabla}^{\mr{ad},  \omega}_{\mcE_\mbT}])$ is a $k$-vector space of dimension 
$2g-2 + \frac{2 g-2 + r + \sum_{i=1}^r \tau^{-1} (-\varepsilon_i)}{p}$ (cf. (\ref{eqeq45}) for the definition of $\tau$).

\end{itemize}
 \ele
\begin{proof}
First, we shall prove assertion (i). 
Since the latter assertion follows  from the former assertion (and Proposition \ref{P0337}), it suffices to prove   the former assertion (i.e., the equalities in (\ref{ee320})).

Let  $\overline{\zeta} :  \widetilde{\mcT}_{\mcE^{\dagger \mr{log}}_\mbT/S^\mr{log}} \migi \mfb_{\mcE'^\dagger_\mbB}$ be  the morphism
obtained from $\zeta$ (defined in the proof of Lemma \ref{L07gg}) by restricting its domain and codomain.
One verifies from the definition of a Miura $\mfg$-oper (in particular, of $\nabla_\mcE$) that 
$\overline{\zeta}$ is surjective.
Both $\widetilde{\mcT}_{\mcE^{\dagger \mr{log}}_\mbT/S^\mr{log}}$ and $\mfb_{\mcE'^\dagger_\mbB}$ are rank $2$ vector bundles, and hence, $\overline{\zeta}$ turns out to be  an isomorphism.
Since $\mr{Im} (\zeta) \subseteq \mfb_{\mcE'^\dagger_\mbB} = \mfg^{-1/1}_{\mcE^\dagger_\mbT}$ and the composite equality (\ref{F3456}) holds,
the pair of $\overline{\zeta}$ and 
the identity morphism of $\Omega_{X^\mr{log}/S^\mr{log}} \otimes \mfb_{\mcE'^\dagger_\mbB}$
specifies  an isomorphism
$\mcK^\bullet [\widetilde{\nabla}^{\mr{ad}, -1/1}_{\mcE_\mbT}] \isom \mcK^\bullet [\nabla^{\mr{ad}}_{\mcE'_\mbB}]$.
This isomorphism restricts to an isomorphism  
\begin{align}
\mcK^\bullet [\widetilde{\nabla}^{\mr{ad}, \omega}_{\mcE_\mbT}] \isom \mcK^\bullet [\nabla^{\mr{ad}, \omega}_{\mcE'_\mbB}],
\end{align}
where $\nabla^{\mr{ad}, \omega}_{\mcE'_\mbB}$ denotes the morphism
$\mfb_{\mcE'_\mbB} \migi \mr{Im} (\nabla^{\mr{ad}}_{\mcE'_\mbB})$ obtained  by restricting the codomain of  $\mcK^\bullet [\nabla^{\mr{ad}}_{\mcE'_\mbB}]$.
 $\mcK^\bullet [\nabla^{\mr{ad}, \omega}_{\mcE'_\mbB}]$
is quasi-isomorphic to  $\mr{Ker} (\nabla^\mr{ad}_{\mcE'_\mbB}) [0]$ via the natural inclusion
$\mr{Ker} (\nabla^\mr{ad}_{\mcE'_\mbB}) [0] \migi \mcK^\bullet [\nabla^{\mr{ad}, \omega}_{\mcE'_\mbB}]$.
Hence, we have  the composite isomorphism
\begin{align} \label{ko090}
\mbH^j (X, \mcK^\bullet [\widetilde{\nabla}^{\mr{ad}, \omega}_{\mcE_\mbT}] )
\isom
\mbH^j (X, \mcK^\bullet [\nabla^{\mr{ad}, \omega}_{\mcE'_\mbB}])
 \isom H^j (X, \mr{Ker} (\nabla^\mr{ad}_{\mcE'_\mbB})).
\end{align}
In particular, (since $\mr{dim} (X) =1$)  the equality $\mbH^2 (X, \mcK^\bullet [\widetilde{\nabla}^{\mr{ad}, \omega}_{\mcE_\mbT}] )= 0$ holds.

Next, let us prove that ($\mbH^0 (X, \mcK^\bullet [\widetilde{\nabla}^{\mr{ad}, \omega}_{\mcE_\mbT}] )\cong$) $H^0 (X, \mr{Ker} (\nabla^\mr{ad}_{\mcE'_\mbB})) = 0$.
Let 
\begin{align}
\kappa : F^*_{X/k} (F_{X/k*}(\mr{Ker} (\nabla^\mr{ad}_{\mcE'_\mbB}))) \migi \mfb_{\mcE'^\dagger_\mbB}
\end{align}
be the $\mcO_X$-linear morphism corresponding, via the adjunction relation ``$F_{X/k}^*(-) \dashv F_{X/k*}(-)$",  to the natural inclusion $F_{X/k*}(\mr{Ker} (\nabla^\mr{ad}_{\mcE'_\mbB})) \migiincl F_{X/k*} (\mfb_{\mcE'_\mbB})$.
Since $\nabla_{\mcE'_\mbB}^\mr{ad}$ has vanishing $p$-curvature, $\kappa$ becomes  an isomorphism when restricted to  $X \setminus \bigcup_{i=1}^r \{ \mr{Im} (\sigma_i) \}$
(cf. ~\cite{K}, \S\,5, p.\,190, Theorem 5.1).
Also,   $\kappa$ is compatible with the respective $k$-connections $\nabla^\mr{can}_{F_{X/k*}(\mr{Ker} (\nabla^\mr{ad}_{\mcE'_\mbB}))}$ and $\nabla^\mr{ad}_{\mcE'_\mbB}$.
Now, suppose that $H^0 (X, \mr{Ker} (\nabla^\mr{ad}_{\mcE'_\mbB})) \  (= H^0 (X^{(1)}_k, F_{X/k*}(\mr{Ker} (\nabla^\mr{ad}_{\mcE'_\mbB}))))\neq  0$, i.e., 
that there exists a nonzero element of $H^0 (X^{(1)}_k, F_{X/k*}(\mr{Ker} (\nabla^\mr{ad}_{\mcE'_\mbB})))$.
This element determines an injection $\mcO_{X^{(1)}_k} \migiincl F_{X/k*}(\mr{Ker} (\nabla^\mr{ad}_{\mcE'_\mbB}))$.
By pulling-back via $F_{X/k}$, we obtain an injection $\mcO_X \migiincl F^*_{X/k} (F_{X/k*}(\mr{Ker} (\nabla^\mr{ad}_{\mcE'_\mbB})))$ compatible with  the respective $k$-connections $d$ and $\nabla^\mr{can}_{F_{X/k*}(\mr{Ker} (\nabla^\mr{ad}_{\mcE'_\mbB}))}$.
Denote by $\kappa' : \mcO_X \migiincl \mfb_{\mcE'_\mbB}$ 
the composite of that injection and $\kappa$.
Here, recall that $\mfb_{\mcE'^\dagger_\mbB} = \mfg^{0}_{\mcE^\dagger_\mbT} \oplus \mfg^{-1}_{\mcE^\dagger_\mbT} \isom \mcO_X \oplus \mcT_{X^\mr{log}/k}$ (cf. (\ref{isomoper})).
Since $\mr{deg} (\mcT_{X^\mr{log}/k}) < 0$, 
the composite of $\kappa'$ and the second projection $\mfb_{\mcE'^\dagger_\mbB}  \ (= \mfg^{0}_{\mcE^\dagger_\mbT} \oplus \mfg^{-1}_{\mcE^\dagger_\mbT})\migisurj  \mfg^{-1}_{\mcE^\dagger_\mbT}$ is zero.
Hence, the image of $\kappa'$ lies in the first component $\mfg^{0}_{\mcE^\dagger_\mbT}  \subseteq \mfb_{\mcE'^\dagger_\mbB}$.
But,   $\mfg^{0}_{\mcE^\dagger_\mbT}$ is (by the definition of $\nabla_\mcE$) not closed under $\nabla^\mr{ad}_{\mcE'_\mbB}$.
This is a contradiction,  and consequently, the equality  $H^0 (X^{(1)}_k, F_{X/k*}(\mr{Ker} (\nabla^\mr{ad}_{\mcE'_\mbB}))) = 0$ holds.
This completes the proof of assertion (i).

Next, let us consider assertion (ii).
Let $X$ and $\widehat{\mcE}^{\spadesuit \diamondsuit}$ be as assumed in the statement of (ii).
By the definition of a Miura $\mfg$-oper,
$\nabla^\mr{ad}_{\mcE'_\mbB}$ may restricts to a $k$-connection $\nabla_{-1}$ on $\mfg^{-1}_{\mcE_\mbT^\dagger}$.
Also, we obtain a $k$-connection  $\nabla_0$ on $\mfg^{0}_{\mcE_\mbT^\dagger}$ induced from $\nabla^\mr{ad}_{\mcE'_\mbB}$ via the projection $\mfb_{\mcE_\mbT^\dagger} \migisurj (\mfb_{\mcE_\mbT^\dagger} / \mfg^{-1}_{\mcE_\mbT^\dagger}=) \ \mfg^{0}_{\mcE_\mbT^\dagger}$.
For each $i \in \{ 1, \cdots ,r \}$, the monodromy  of $(\mcE'^\dagger_{\mbB},  \nabla_{\mcE'_\mbB})$ at $\sigma_i$
coincides with $\begin{pmatrix} -\frac{\varepsilon_i}{2} & 0 \\ 1 & \frac{\varepsilon_i}{2}\end{pmatrix} \in \Gamma (S, \sigma^*_i (\mfb_{\mcE'^\dagger_\mbB})) = \mfb^{-} (k)$ ($\subseteq  \mfs \mfl_2 (k)$).
It follows that the monodromies at $\sigma_i$ of $(\mfg^{-1}_{\mcE_\mbT^\dagger}, \nabla_{-1})$ and $(\mfg^{0}_{\mcE_\mbT^\dagger}, \nabla_0)$  are $\varepsilon_i$ and $0$ respectively.
(Indeed, the monodromy of $\nabla_{\mcE'_\mbB}^\mr{ad}$ at $\sigma_i$ may be expressed, by means of the basis  $\Big\langle \begin{pmatrix} 1 & 0 \\ 0 & -1\end{pmatrix},  \begin{pmatrix} 0 & 0 \\ 1 & 0\end{pmatrix} \Big\rangle$ of $\mfb^{-} (k)$,  as the matrix $\begin{pmatrix} 0 & 0 \\ 2 & \varepsilon_i\end{pmatrix}$.)
Hence, 
if $\kappa_{-1} :  F^*_{X/k} (F_{X/k*} (\mr{Ker} (\nabla_{-1}))) \migiincl \mfg^{-1}_{\mcE_\mbT^\dagger}$ and $\kappa_{0} : F^*_{X/k} (F_{X/k*} (\mr{Ker} (\nabla_{0}))) \migiincl \mfg^{0}_{\mcE_\mbT^\dagger}$ are the injections obtained in the same manner as $\kappa$,
then $\kappa_0$ is an isomorphism and 
$\mr{Coker} (\kappa_{-1}) \cong \bigoplus_{i=1}^r \Lambda_i$, where each $\Lambda_i$ ($i =1, \cdots, r$) is an $\mcO_X$-module supported on $\mr{Im} (\sigma_i)$ of length $\tau^{-1}(-\varepsilon_i) \in \widetilde{\mbF}_p$.
On the other hand, by  ~\cite{Og3}, Corollary 3.2.2 (and the fact that $F_{X/k}$ is flat), 
the following natural sequence  turns out to be  exact:
\begin{align}\label{ee340}
0 \migi F^*_{X/k} (F_{X/k*} (\mr{Ker} (\nabla_{-1}))) & \migi  F^*_{X/k} (F_{X/k*} (\mr{Ker} (\nabla_{\mcE'_\mbB}^\mr{ad}))) \\
& \migi  F^*_{X/k} (F_{X/k*} (\mr{Ker} (\nabla_{0}))) \migi 0 \notag
\end{align}
Thus,
  the following sequence of equalities holds:
\begin{align} \label{gqf}
& \ \  \ \ \mr{deg} (F_{X/k*}(\mr{Ker} (\nabla^\mr{ad}_{\mcE'_\mbB}))) \\
& = \frac{1}{p} \cdot \mr{deg} (F_{X/k}^* (F_{X/k*}(\mr{Ker} (\nabla^\mr{ad}_{\mcE'_\mbB})))) \notag  \\
& = \frac{1}{p} \cdot \left(\mr{deg}(F^*_{X/k} (F_{X/k*} (\mr{Ker} (\nabla_{0})))) + \mr{deg}(F^*_{X/k} (F_{X/k*} (\mr{Ker} (\nabla_{-1}))))\right) \notag \\
& = \frac{1}{p} \cdot \left(\mr{deg} (\mfg^0_{\mcE^\dagger_\mbT}) + \left(\mr{deg} (\mfg^{-1}_{\mcE^\dagger_\mbT}) - \sum_{i=1}^r \mr{length} (\Lambda_i)\right)\right) \notag \\
& = \frac{1}{p} \cdot \left(\mr{deg} (\mcO_X) + \mr{deg} (\mcT_{X^\mr{log}/k}) - \sum_{i=1}^r \tau^{-1}(-\varepsilon_i) \right) \notag \\
& =  - \frac{2g-2  +  r +  \sum_{i=1}^r \tau^{-1}(-\varepsilon_i)}{p}.  \notag
\end{align}
Since $X$ was assumed to be  smooth (hence $F_{X/k*}(\mr{Ker} (\nabla^\mr{ad}_{\mcE'_\mbB}))$ is a rank $2$ vector bundle), the following sequence of equalities holds:
\begin{align}
&  \ \ \ \ \mr{dim} (H^1 (X, \mr{Ker} (\nabla^\mr{ad}_{\mcE'_\mbB})))  \\
& =  \mr{dim} (H^1 (X, F_{X/k*}(\mr{Ker} (\nabla^\mr{ad}_{\mcE'_\mbB})))) \notag  \\
& =  \mr{dim} (H^0 (X, F_{X/k*}(\mr{Ker} (\nabla^\mr{ad}_{\mcE'_\mbB}))))-\mr{rk} (F_{X/k*}(\mr{Ker} (\nabla^\mr{ad}_{\mcE'_\mbB}))) \cdot (1-g) \notag \\
& \ \ \ - \mr{deg} (F_{X/k*}(\mr{Ker} (\nabla^\mr{ad}_{\mcE'_\mbB})))  \notag \\
& = 0 - 2 \cdot (1-g) +  \frac{2g-2  + r +  \sum_{i=1}^r \tau^{-1}(-\varepsilon_i)}{p}  \notag  \\
& = 2g-2 + \frac{2g-2  + r +  \sum_{i=1}^r \tau^{-1}(-\varepsilon_i)}{p}, \notag
\end{align} 
where the second equality follows from  the Riemann-Roch theorem and  the third equality follows from assertion (i) and (\ref{gqf}).
This completes the proof of assertion (ii).
\end{proof}
\vspace{3mm}

We shall write
\begin{align} \label{ee330}
\mfM \mfO \mfp^{^\mr{Zzz...}}_{\mfs \mfl_2, g,r,  \vec{\varepsilon}} := \mfM \mfO \mfp^{^\mr{Zzz...}}_{\mfs \mfl_2, g,r, \vec{\varepsilon}} \times_{\overline{\mfM}_{g,r}} \mfM_{g,r}
\end{align}
(where $\vec{\varepsilon} \in \mft_{\mr{PGL}_2} (\mbF_p)^{\times r}$ or $\vec{\varepsilon} = \emptyset$).
By Lemma \ref{L07} (and Theorem \ref{th013} (ii)), the following assertion holds.

\vspace{3mm}
\bt  \label{T01} \leavevmode\\
 \ \ \ 
Let $\vec{\varepsilon} \in \mbF_p^{\times r}$ (where  $\vec{\varepsilon} := \emptyset$ if $r =0$).
Then,  the stack $\mfM \overline{\mfO} \mfp^{^\mr{Zzz...}}_{\mfs \mfl_2, g,r, [\vec{\varepsilon}\,]}$ ($\stackrel{(\ref{E0010})}{\cong} \overline{\mfT} \mfa \mfn_{g,r, -\vec{\varepsilon}}$)
is
  a (possibly empty) smooth proper   Deligne-Mumford stack over $k$.
 If $2g-2 + \frac{2g-2 + r +  \sum_{i=1}^r \tau^{-1}(-\varepsilon_i)}{p} < 0$,
then $\mfM \mfO \mfp^{^\mr{Zzz...}}_{\mfs \mfl_2, g,r, [\vec{\varepsilon}\,]}$ is empty.
Moreover, if  $2g-2 + \frac{2g-2 + r +  \sum_{i=1}^r \tau^{-1}(-\varepsilon_i)}{p} \in \mbZ_{\geq  0}$, then any irreducible component $\mcN$ of $\mfM \overline{\mfO} \mfp_{\mfs \mfl_2, g,r, [\vec{\varepsilon}\,]}^{^\mr{Zzz...}}$
   with $\mcN \times_{\overline{\mfM}_{g,r}} \mfM_{g,r} \neq \emptyset$ is  equidimensional of dimension $2g-2 + \frac{2g-2 +r + \sum_{i=1}^r \tau^{-1} (-\varepsilon_i)}{p}$. 
     \et
\vspace{3mm}

In particular, the above corollary (of the case where $r =0$) and Remark \ref{RRR040} (i) below  imply Theorem B.

\vspace{3mm}
\begin{rema} \label{RRR040}
\leavevmode\\
 \ \ \ In this remark, let us mention two  facts concerning the non-emptiness  of $\mfM \mfO \mfp^{^\mr{Zzz...}}_{\mfs \mfl_2, g,r}$ deduced from the previous works.
\begin{itemize}
\item[(i)]
Let $l$ be an integer with $lp \geq 4$.
According to ~\cite{Ray}, Example, it is verified that
\begin{align} \label{uuu4}
\mfM \mfO \mfp^{^\mr{Zzz...}}_{\mfs \mfl_2, \frac{(lp-1)(lp-2)}{2},0} \neq \emptyset.
\end{align}
In order to prove this, (under the assumption that $k$ is an algebraically closed field) let us consider the smooth projective  curve $X$ in $\mbP^2$ ($= \mr{Proj} (k[x, y, z])$)  defined by the equation $x^{lp} - xy^{lp-1} - yz^{lp-1}=0$.
The genus of $X$ is given by $\frac{(lp-1)(lp-2)}{2}$ ($\geq 2$).
One verifies  that $\Omega_{X/k} = \mcO_X ( lp (lp -3) \cdot P_\infty)$, where $P_\infty =[0 : 0: 1]$,
and 
the line bundle $(\mr{id}_X \times F_k)^*(\mcO_{X} (l (lp-3) \cdot P_\infty))$ on $X^{(1)}_k$ 
  specifies a Tango structure.
The existence of this Tango structure  implies (\ref{uuu4}).
\item[(ii)]
Suppose that $g \geq p$.
Then,  there exists an element $\vec{\varepsilon} \in \mbF_p^{\times r}$ (or $\vec{\varepsilon} = \emptyset$)
such  that
$\mfM \mfO \mfp_{\mfs \mfl_2, g,r, [\vec{\varepsilon}\,]}^{^\mr{Zzz...}} \neq  \emptyset$ 
and the composite projection 
\begin{align}
\mfM \mfO \mfp_{\mfs \mfl_2, g,r, [\vec{\varepsilon}\,]}^{^\mr{Zzz...}}
 \migi \mfM_{g,r} \migi \mfM_{g, r-1} \migi \cdots \migi \mfM_{g, 0}
\end{align}
 has dense image.
Indeed, let us take a geometric generic point   $c : \mr{Spec} (K) \migi \mfM_{g, 0}$ of $\mfM_{g,0}$, where $K$ denotes some  algebraically closed field over $k$.
Denote by $\mfX  := (X/K, \{ \sigma_i \}_{i=1}^r)$ the (unpointed) proper smooth  curve over $K$ classified by $c$.
It follows from  ~\cite{TakYok},  Corollary 1.5,  that
there exist a (possibly empty) collection of $K$-rational points $\sigma_{1}, \cdots, \sigma_r$ of $X$, 
  an element  $\vec{\varepsilon}$ in $\mbF_p^{\times r}$, and a pre-Tango structure $\nabla$ on $\mfX$ of monodromies $-\vec{\varepsilon}$.
 $\nabla$ specifies a $K$-rational point of  $\mfM \mfO \mfp_{\mfs \mfl_2, g,r,  [\vec{\varepsilon}\,]}^{^\mr{Zzz...}}$ over  $c$.
 This implies  that the image of the composite projection   $\mfM \mfO \mfp_{\mfs \mfl_2, g,r,  [\vec{\varepsilon}\,]}^{^\mr{Zzz...}} \migi \mfM_{g, 0}$ is dense.
\end{itemize}
\end{rema}

\vspace{5mm}
\subsection{The case of $g=1$}
\leavevmode\\ \vspace{-4mm}

In this last subsection, we study the moduli stack $\mfM \mfO \mfp_{\mfs \mfl_2,  g,r, \vec{\varepsilon}}^{^\mr{Zzz...}}$ of the case where $g =1$.
To begin with, we shall observe the following assertion.

\vspace{3mm}
\bpr  \label{dT01} \leavevmode\\
 \ \ \ 
 Let $\vec{\varepsilon} := (\varepsilon_i)_{i=1}^r \in \mbF_p^{\times r}$ (where  $\vec{\varepsilon} := \emptyset$ if $r =0$).
For any positive integer $s$,  we shall write
$\vec{\varepsilon}_{+ 1 \times s} := (\varepsilon_1,  \cdots, \varepsilon_r, 1,1, \cdots, 1) \in \mbF_p^{\times (r +s)}$, where the last $s$ factors are all $1$.
Then, there exists a canonical isomorphism
\begin{align} \label{wert}
 \mfM \mfO \mfp_{\mfs \mfl_2, g,r, [\vec{\varepsilon}\,]}^{^\mr{Zzz...}} \times_{\mfM_{g,r}} \mfM_{g, r +s}  \isom  \mfM \mfO \mfp_{\mfs \mfl_2, g, r+s, [\vec{\varepsilon}_{+1 \times s}\,]}^{^\mr{Zzz...}}
\end{align}
over $\mfM_{g, r +s}$.
In particular, 
\begin{align}
``\, \mfM \mfO \mfp_{\mfs \mfl_2, g,r, [\vec{\varepsilon}\,]}^{^\mr{Zzz...}} \neq \emptyset\,"  \Longleftrightarrow ``\, \mfM \mfO \mfp_{\mfs \mfl_2, g, r+s, [\vec{\varepsilon}_{+1 \times s}\,]}^{^\mr{Zzz...}} \neq \emptyset\,".
\end{align}

  \epr
\begin{proof}
Let $S$ be a $k$-scheme and $\mfX := (X/S, \{ \sigma_i \}_{i=1}^{r+s})$ an $(r+s)$-pointed smooth  curve over $S$ of genus $g$.
Denote by $\overline{\mfX} := (X/S, \{ \sigma_i \}_{i=1}^r)$ the pointed smooth  curve obtained from $\mfX$ by forgetting the last $s$ marked points $\{ \sigma_{i} \}_{i=r+1}^s$ (hence, $S^{\mfX \text{-} \mr{log}} = S^{\overline{\mfX} \text{-} \mr{log}} = S$).
Suppose that we are given a pre-Tango structure  $\nabla$ on $\overline{\mfX}$ of monodromies $-\vec{\varepsilon}$.
The pair $(\nabla, \mfX)$ specifies an $S$-rational point of $\overline{\mfT} \mfa \mfn_{g,r, - \vec{\varepsilon}} \times_{\overline{\mfM}_{g,r}} \mfM_{g, r+s}$.
One may find   a unique 
Tango structure $\nabla_+$ on $\mfX$
whose restriction to $\Omega_{X^{\overline{\mfX} \text{-}\mr{log}}/S}$  ($= \Omega_{X^{\mfX \text{-}\mr{log}}/S} (-\sum_{i=r+1}^r \sigma_i)$) coincides with $\nabla$.
Moreover, the monodromy of $\nabla_+$ at $\sigma_i$ (for each $i = r +1, \cdots, r +s$) coincides with $-1$.
That is to say,
$(\nabla_+, \mfX)$ specifies an $S$-rational point of $\overline{\mfT} \mfa \mfn_{g, r+s, -(\vec{\varepsilon}_{+ 1 \times s})} \times_{\overline{\mfM}_{g,r+s}} \mfM_{g, r+s}$.
One verifies immediately  that the resulting  morphism
\begin{align} \label{ffde}
\overline{\mfT} \mfa \mfn_{g,r, - \vec{\varepsilon}} \times_{\overline{\mfM}_{g,r}} \mfM_{g, r+s} \migi \overline{\mfT} \mfa \mfn_{g, r+s, -(\vec{\varepsilon}_{+ 1 \times s})} \times_{\overline{\mfM}_{g,r+s}} \mfM_{g, r+s}
\end{align}
(i.e., the morphism given by $(\nabla, \mfX) \mapsto (\nabla_+, \mfX)$) determines an isomorphism  over $\mfM_{g, r+s}$.
Hence,
 the desired isomorphism may be obtained from this isomorphism and (\ref{E0010}).
\end{proof}

\vspace{3mm}
\bco  \label{dT01c} \leavevmode\\
 \ \ \ 
Suppose that $r >0$, and write $\vec{1}_r := (1, 1, \cdots, 1) \in \mbF_p^{\times r}$. 
Then, $\mfM \mfO \mfp_{\mfs \mfl_2, 1, r, \vec{1}_r}^{^\mr{Zzz...}}$ is a  nonempty, geometrically connected, and  smooth  Deligne-Mumford stack over $k$ of dimension $r$.
Moreover, the forgetting morphism   $\mfM \mfO \mfp_{\mfs \mfl_2, 1, r, \vec{1}_r}^{^\mr{Zzz...}} \migi \mfM_{1,r}$
is finite, surjective, and generically \'{e}tale of degree $p-1$.
  \eco
\begin{proof}
According to Example \ref{E0140},
there exists a pre-Tango structure $\nabla$  on any ordinary elliptic curve $X/k$ (without marked points).
It follows from Proposition \ref{dT01} (more precisely, from  the discussion in the proof of that proposition) that
$\nabla$ gives rise to  a pre-Tango structure  of monodromies $-\vec{1}_r$ on $X/k$ with $r$ marked points.
This implies that 
$\mfM \mfO \mfp_{\mfs \mfl_2, 1, r, \vec{1}_r}^{^\mr{Zzz...}}$ ($\stackrel{}{\cong} \overline{\mfT} \mfa \mfn_{g,r, -\vec{1}_r} \times_{\overline{\mfM}_{1,r}} \mfM_{1, r}$ by (\ref{E0010}))
 is nonempty.
By Theorem  \ref{T01},  we can verify the smoothness 
 and the calculation of the dimension, and moreover, the fact that 
 the forgetting morphism  $\mfM \mfO \mfp_{\mfs \mfl_2, 1, r, \vec{1}_r}^{^\mr{Zzz...}} \migi \mfM_{1,r}$
is finite and  surjective.
Also, according to the discussion in Example \ref{E0140},
$\mfM \mfO \mfp_{\mfs \mfl_2, 1, r, \vec{1}_r}^{^\mr{Zzz...}}$ is isomorphic to 
$\mfC \mfo^{\psi =0}_{\mcO_{\mfC_{1, r}}, 1, r, (0, 0, \cdots, 0)} \setminus \mcI$, where $\mcI$ denotes 
the component
 classifying
  the trivial connections (i.e., the universal derivations) on curves.
 In particular, the forgetting morphism  $\mfM \mfO \mfp_{\mfs \mfl_2, 1, r, \vec{1}_r}^{^\mr{Zzz...}}   \migi \mfM_{1, r}$ is generically \'{e}tele and of degree $p-1$.
(More precisely, the forgetting morphism  is \'{e}tale at the points classifying pointed smooth  curves whose underlying curves are  ordinary.)
Its fiber 
 over the point classifying any supersingular  curve
  consists precisely of one point.
Hence, 
$\mfM \mfO \mfp_{\mfs \mfl_2, 1, r, \vec{1}_r}^{^\mr{Zzz...}}$ ($\cong \mfC \mfo^{\psi =0}_{\mcO_{\mfC_{1, r}}, 1, r, (0, 0, \cdots, 0)} \setminus \mcI$)
is geometrically connected.
This completes the proof of Corollary \ref{dT01c}.
\end{proof}
\vspace{3mm}

\vspace{10mm}
\section{Pathology in positive characterisitc} \label{GHTY}\vspace{3mm}

In this last section, we study the pathology of algebraic geometry in positive  characteristic, which is of certain interest, since pathology reveals some completely different geometric phenomena from those in complex geometry.

\vspace{5mm}
\subsection{Generalized Tango curves.}
\leavevmode\\ \vspace{-4mm}

In what follows, suppose that $k$ is an algebraically closed field of characteristic $p >2$.

\vspace{0mm}
\bde \label{Yff19}[cf. ~\cite{Tak1}, \S\,3]\leavevmode\\
 \ \ \  
Let $S$ be a $k$-scheme and $l$ a positive integer.
\begin{itemize}
\item[(i)]
A  {\bf generalized Tango curve of index $l$} over $S$ is a quadruple
\begin{align}
\widehat{\mfX} := (X,  \mcL, \mcN, \nu),
\end{align}
where $X$ denotes   a geometrically connected, proper, and  smooth curve $X$  over $S$, $\mcL$ denotes a Tango structure on $X/S$, $\mcN$ a line bundle  on $X^{(1)}_S$, and $\nu$ denotes an isomorphism
$\nu : \mcN^{\otimes (lp -1)} \isom \mcL$. 
(Notice that this  definition is equivalent to the definition of a generalized Tango curve of index $lp-1$ in the sense of ~\cite{Tak1}, \S\,3.)
\item[(ii)]
Let $\widehat{\mfX} := (X,  \mcL, \mcN, \nu)$ and $\widehat{\mfX}' := (X',  \mcL', \mcN', \nu')$
be generalized Tango curves of index $l$ over $S$.
An {\bf isomorphism of generalized Tango curves} from $\widehat{\mfX}$ to $\widehat{\mfX}'$
is a pair $(h_X, h_\mcN)$ consisting of  an isomorphism $h_X : X \isom X'$ over $S$ (where we shall write $h_{X^{(1)}_S}$ for the isomorphism $X^{(1)}_S \isom X'^{(1)}_S$ obtained  from $h_X$ via base-change by $F_S$)  and an isomorphism
$h_\mcN : h_{X^{(1)}_S}^* (\mcN') \isom \mcN$
satisfying the following conditions:
\begin{itemize}
\item[$\bullet$]
The isomorphism 
$h^*_{X^{(1)}_S} (F_{X'/S*}(\Omega_{X'/S})) \ (\cong F_{X/S*}(h^*_X (\Omega_{X'/S}))) \isom F_{X/S*}(\Omega_{X/S})$ induced by $h_X$ restricts to an isomorphism
\begin{align} \label{ee400}
h_{X^{(1)}_S}^*(\mcL') \isom \mcL.
\end{align}
\item[$\bullet$]
The  square diagram
\begin{align}
\begin{CD}
h_{X^{(1)}_S}^*(\mcN')^{\otimes (lp-1)}@> h_\mcN^{\otimes (lp-1)}> \sim > \mcN^{\otimes (lp-1)}
\\
@V h_{X^{(1)}_S}^*(\nu') V \wr V @V \wr V \nu V
\\
h_{X^{(1)}_S}^*(\mcL') @>\sim >(\ref{ee400})> \mcL,
\end{CD}
\end{align}
 is commutative.
\end{itemize}
\end{itemize}
   \ede
\vspace{3mm}

Let $g$ be an integer with $g > 1$.
Since the pull-back of a generalized Tango curve by a morphism $T \migi S$ of $k$-schemes can be defined in a natural manner, we obtain a stack in groupoids
\begin{align}
\mfG \mfT \mfa \mfn_{g}^l
\end{align}
 over $\mfS \mfc \mfh_{/\mr{Spec} (k)}$ whose category of sections
over an object $S$ of $\mfS \mfc \mfh_{/\mr{Spec} (k)}$ is the groupoid of generalized Tango curves of index $l$ over $S$.
According to ~\cite{Lang} and ~\cite{Muk}, 
  one may construct a family 
   of algebraic varieties parametrized by $\mfG \mfT \mfa \mfn_{g}^l$ each of whose fiber forms   a counterexample to the Kodaira  vanishing theorem (cf. \S\,\ref{wwww} below).

\vspace{3mm}
\begin{rema} \label{Rrg040}
\leavevmode\\
 \ \ \ {\it $\mfG \mfT \mfa \mfn_{g}^l$ may be represented by either the empty stack or  an equidimensional  smooth  Deligne-Mumford over $k$  of dimension $2g-2 + \frac{2g-2}{p}$,  and the forgetting morphism $\mfG \mfT \mfa \mfn_{g}^l \migi \mfT \mfa \mfn_{g}$ (i.e., the morphism given by  $(\mfX, \mcL, \mcN, \nu) \mapsto (X, \mcL)$) is finite and \'{e}tale}.
Indeed, for each integer $d$, let us denote by $\mfP \mfi \mfc^d_g$ the Picard stack for the universal family of curves $f_{\mft \mfa \mfu} : \mfC_{g} \migi \mfM_g$ classifying line bundles  of relative degree $d$.
For each pair of integers  $(d, e)$ with $e \geq 1$, denote by $\mu_{d, e}$ the morphism $\mfP \mfi \mfc^d_g \migi \mfP \mfi \mfc_g^{de}$ over $\mfM_g$ given by  assigning $\mcL \mapsto \mcL^{\otimes e}$ (for any line bundle $\mcL$ of relative degree $d$).
Then, we obtain a canonical isomorphism 
\begin{align}  \label{gggh1}
\mfG \mfT \mfa \mfn_{g}^l \hspace{3mm}& \isom \mfT \mfa \mfn_g \times_{\mfP \mfi \mfc_g^{N(lp-1)}, \mu_{N, lp-1}} \mfP \mfi \mfc_g^{N}. \notag \\
\vin \hspace{9mm}& \hspace{30mm}\vin  \\
(X, \mcL, \mcN, \nu) &  \mapsto \hspace{15mm}  ((X, \mcL), \mcN) \notag 
\end{align}
over $\mfM_g$, where $N :=\frac{2g-2}{p (lp-1)}$.
Since $p \nmid lp-1$, the morphism  $\mu_{N, lp-1}  : \mfP \mfi \mfc_g^{N} \migi \mfP \mfi \mfc_g^{N(lp-1)}$ is finite, surjective,  and \'{e}tale.
Hence, by  Theorem  \ref{T01},  the fiber product $\mfT \mfa \mfn_g \times_{\mfP \mfi \mfc_g^{\frac{2g-2}{p}}} \mfP \mfi \mfc_g^{N}$ ($\stackrel{(\ref{gggh1})}{\cong} \mfG \mfT \mfa \mfn_{g, l}$) is an equidimensional  smooth  Deligne-Mumford over $k$  of dimension $2g-2 + \frac{2g-2}{p}$, as desired.
\end{rema}

\vspace{5mm}
\subsection{Generalized Raynaud surfaces.} \label{wwww}
\leavevmode\\ \vspace{-4mm}

In what follows, we construct, by means of a certain closed substack of  $\mfG \mfT \mfa \mfn_{g}^l$, a family of algebraic surfaces (of general type) which is parametrized by a high dimensional variety and  each of whose fiber has the  automorphism  group scheme which is not reduced.

First, recall the {\it generalized Raynaud surface} associated with a Tango structure.
Let $S$  be a $k$-scheme, $X$ a proper smooth curve over $S$ of genus $g$, 
and $\mcL$ ($\subseteq F_{X/S*}(\mcB_{X/S})$) a Tango structure on $X/S$.
Suppose that $l p (p-1) = 2g-2$ for some positive integer $l$, and that there exists a line bundle $\mcN$ on $X^{(1)}_S$ of relative degree $l$ (relative to $S$) admitting an isomorphism  $\nu : \mcN^{\otimes (p-1)} \cong \mcL$.
In particular, the quadruple $\widehat{\mfX} := (X, \mcL, \mcN, \nu)$ specifies a generalized Tango curve of index $1$  over $S$.

Let us take an affine open covering $\{ U_{\alpha} \}_{\alpha \in I}$ of $X$ 
such that there exists an isomorphism $\xi_{\alpha} :\mcO_{(U_\alpha)_S^{(1)}} \isom  \mcN |_{(U_\alpha)_S^{(1)}} $.
Write $I_2 := \{ (\alpha, \beta) \in I \times I \ | \ U_{\alpha \beta}  \neq \emptyset \}$ (where $U_{\alpha \beta} := U_\alpha \cap U_\beta$).
For each pair $(\alpha, \beta) \in I_2$,
the automorphism $\xi_\beta |_{(U_{\alpha \beta})^{(1)}_S} \circ \xi_{\alpha}^{-1} |_{(U_{\alpha \beta})^{(1)}_S}$ of $ \mcN |_{(U_{\alpha \beta})_S^{(1)}} $ is given by 
multiplication by some element  $u_{\alpha \beta} \in \Gamma ((U_{\alpha \beta})_S^{(1)}, \mcO^\times_{(U_{\alpha \beta})^{(1)}_S})$.
For each $\alpha \in I$, the composite of $F^*_{X/S}(\xi_\alpha^{\otimes (p-1)}) : \mcO_{U_\alpha} \isom F^*_{X/S} (\mcN |_{(U_\alpha)_S^{(1)}} ^{\otimes (p-1)})$ and the restriction of the
composite isomorphism $F^*_{X/S} (\mcN^{\otimes (p-1)})$ $\xrightarrow{F^*_{X/S} (\nu)} F^*_{X/S} (\mcL) \xrightarrow{\xi_\mcL} \Omega_{X/S}$ to $U_\alpha$ 
determines an element of $\Gamma (U_\alpha, \Omega_{X/S})$;
after possibly replacing $\{ U_\alpha \}_\alpha$ with its refinement,
this  element may be expressed as $d t_\alpha$ for some $t_\alpha \in \Gamma (U_\alpha, \mcO_{U_\alpha})$
(hence, we have $d t_\alpha = F^*_{U_{\alpha \beta}/S}(u_{\alpha \beta}^{p-1})dt_\beta$ for each $(\alpha, \beta) \in I_2$).
Moreover, one may find,  after possibly replacing $\{ U_\alpha \}_\alpha$ with its refinement,   a collection 
$\{ r_{\alpha \beta} \}_{(\alpha,  \beta) \in I_2}$ (where each $r_{\alpha \beta}$ is an element of $\Gamma ((U_{\alpha \beta})^{(1)}_S, \mcO_{(U_{\alpha \beta})^{(1)}_S})$) satisfying that
$t_\alpha = F^*_{U_{\alpha \beta}/S}(u_{\alpha \beta}^{p-1}) t_\beta - F^*_{U_{\alpha \beta}/S} (r_{\alpha \beta})$.
Write  $P_\alpha$ ($\alpha \in I$)
for the closed subscheme of $\mbP^2 \times (U_\alpha)_S^{(1)}$ (where $\mbP^2 := \mr{Proj} (k [x_\alpha, y_\alpha, z_\alpha])$, i.e.,   the projective plane over $k$) defined by the equation 
\begin{align}
y_\alpha^{p-1}z_\alpha = x_\alpha^p + (\mr{id}_{U_\alpha} \times  F_S)^*(t_\alpha) z_\alpha^{p}.
\end{align}
One may  glue together $\{ P_\alpha \}_{\alpha \in I}$  by means of the relations
\begin{align} 
x_\alpha & = u_{\alpha \beta}^{p-1} x_\beta +  r_{\alpha \beta}z_\beta \\
y_\alpha & = u_{\alpha \beta}^p y_\beta \notag  \\
z_\alpha & = z_\beta. \notag 
\end{align}
Denote by $P$ the resulting algebraic surface and refer to  it  as the  {\bf  generalized Raynaud surface} associated with the generalized Tango curve $\widehat{\mfX}$.
Denote by $\Psi : P \migi X^{(1)}_S$ the natural projection.
For each point $s$ of $S$, we shall denote by $P_s$ and  $X^{(1)}_{s}$ the  fibers  over $s$  of the composite projection $P \stackrel{\Psi}{\migi} X^{(1)}_S \migi S$ and  the projection $X^{(1)}_S \migi S$ respectively.
According to ~\cite{Ray} (or, ~\cite{Tak3}, Theorem 3.1),  the fiber $P_s$ (for any $k$-rational point $s$ of $S$) is a proper smooth algebraic surface over $k$,   and it is 
 of general type  (resp., a quasi-elliptic surface) if $p >3$ (resp., $p =3$).
 Moreover, by ~\cite{Tak2}, Theorem 2.1, there exists   an isomorphism
\begin{align} \label{ee401}
\Gamma (P_s, \mcT_{P_s/k}) \cong \Gamma (X^{(1)}_s, \mcN |_{X^{(1)}_s}).
\end{align}

Here, we shall assume  that $p >3$  and $\Gamma (X^{(1)}_s, \mcN |_{X^{(1)}_s}) \neq 0$.
Since $X_s^{(1)}$ is a surface of general type, the automorphism group scheme $\mr{Aut}_k (X_s^{(1)})$   of $X_s^{(1)}$ is finite.
Hence,  because of (\ref{ee401}), $\mr{Aut}_k (X_s)$ is not reduced; this fact   may be thought of as a  pathological phenomenon  (relative to zero characteristic) of algebraic geometry in positive characteristic.

Denote by
\begin{align}
\mfG \mfT \mfa \mfn^1_{g,  \Gamma \neq 0}
\end{align}
the closed substack of $\mfG \mfT \mfa \mfn_{g}^1$  ($\stackrel{(\ref{gggh1})}{\cong} \mfT \mfa \mfn_g \times_{\mfP \mfi \mfc_g^{\frac{2g-2}{p}}} \mfP \mfi \mfc_g^{\frac{2g-2}{p (p-1)}}$)
classifying generalized Tango curves   $\widehat{\mfX} := (X/k, \mcL, \mcN, \nu)$
with $\Gamma (X^{(1)}_k, \mcN) \neq 0$.
The tautological family of generalized Tango curves over $\mfG \mfT \mfa \mfn_{g, \Gamma \neq 0}^1$
induces, by means of the above discussion, a family
\begin{align} \label{sssqw}
\mfP \migi \mfG \mfT \mfa \mfn_{g, \Gamma \neq 0}^1
\end{align}
 of proper smooth surfaces of general type  parametrized by $\mfG \mfT \mfa \mfn_{g, \Gamma \neq 0}^1$  all of whose fibers have the automorphism group schemes being  non-reduced. 

\vspace{3mm}
\bt \label{P0fgh7} \leavevmode\\
 \ \ \ 
Suppose that $p >3$, $p (p-1) | 2g-2$, and  $4 | p-3$.
Then, 
$\mfG \mfT \mfa \mfn_{g, \Gamma \neq 0}^1$ is a nonempty  closed substack of $\mfG \mfT \mfa \mfn_{g}^1$ of dimension $\geq  g-2 + \frac{2g-2}{p-1}$.
  \et
\begin{proof}
It follows from ~\cite{Tak3}, Theorem 4.1,  that
 $\mfG \mfT \mfa \mfn_{g, \Gamma \neq 0}^1$ is nonempty.
In the following, we shall 
 calculate the dimension of $\mfG \mfT \mfa \mfn_{g,  \Gamma \neq 0}^1$.
Let us take  a $k$-scheme $T$ together with an \'{e}tale surjective morphism $T \migi \mfG \mfT \mfa \mfn_{g}^1$;
this morphism classifies a  generalized Tango curve  $(Y, \mcL_Y, \mcN_Y, \nu_Y)$
 of index $1$ over $T$.
Since $\mfG \mfT \mfa \mfn_{g}^1$, as well as $T$, has dimension  $2g-2 + \frac{2g-2}{p}$ (cf. Remark \ref{Rrg040}),  the relative Picard scheme $\mr{Pic}_{Y/T}^{N}$   of $Y/T$ (cf. the proof of Proposition \ref{prfg001tt}), where $N:= \frac{2g-2}{p(p-1)}$,  has dimension 
 ($2g-2 + \frac{2g-2}{p} + g=$) $3g-2 + \frac{2g-2}{p}$.
Denote by 
\begin{align}
\mr{Pic}_{Y/T, \Gamma \neq 0}^{N} \ \  \ \left(\text{resp.,}  \ \mfP \mfi \mfc_{Y/T, \Gamma \neq 0}^{N}\right)
\end{align}
 the closed subscheme of $\mr{Pic}_{Y/T}^{N}$ (resp., the closed substack of $\mfP \mfi \mfc_{Y/T}^{N} := \mfP \mfi \mfc_{g}^{N} \times_{\mfM_g} T$) classifying line bundles admitting a nontrivial global section. 
The isomorphism displayed in  (\ref{gggh1}) (of the case where $l =1$) restricts to an isomorphism
\begin{align}
\mfG \mfT \mfa \mfn_{g, \Gamma \neq 0}^1 \isom \mfT \mfa \mfn_g \times_{\mfP \mfi \mfc_g^{N(p-1)}, \mu_{N, p-1}} \mfP \mfi \mfc_{Y/T, \Gamma \neq 0}^{N}.
\end{align}
Thus, it suffices to calculate the dimension of $\mfT \mfa \mfn_g \times_{\mfP \mfi \mfc_g^{N(p-1)}} \mfP \mfi \mfc_{Y/T, \Gamma \neq 0}^{N}$.
By a well-known fact of the Brill-Noether theory (cf. ~\cite{G}, ~\cite{L}), 
the closed subscheme $\mr{Pic}_{Y/T, \Gamma \neq 0}^{N}$ of $\mr{Pic}_{Y/T}^{N}$ is of codimension $\leq g - N$.
On the other hand, the closed subscheme $T \times_{[\mcL_Y], \mr{Pic}_{Y/T}^{N(p-1)}, \overline{\mu}_{N, p-1}} \mr{Pic}_{Y/T}^{N}$ of $\mr{Pic}_{Y/T}^{N}$  is of codimension $g$, where $[\mcL_Y] : T \migi \mr{Pic}_{Y/T}^{N(p-1)}$ denotes the classifying morphism of $\mcL_Y$ and $\overline{\mu}_{N, p-1} : \mr{Pic}_{Y/T}^{N} \migi \mr{Pic}_{Y/T}^{N(p-1)}$ denotes the finite  \'{e}tale morphism given by assigning $\mcL' \mapsto \mcL'^{\otimes (p-1)}$ (for any line bundle $\mcL'$).
Hence, 
their intersection
\begin{align}
& \ \ \ \  \ \ \ T \times_{[\mcL_Y], \mr{Pic}_{Y/T}^{N(p-1)}
 } \mr{Pic}_{Y/T, \Gamma \neq 0}^{N} \\
&  \left(\cong \Big(T \times_{[\mcL_Y], \mr{Pic}_{Y/T}^{N(p-1)}, \overline{\mu}_{N, p-1}} \mr{Pic}_{Y/T}^{N} \Big) \times_{\mr{Pic}_{Y/T}^{N}} \mr{Pic}_{Y/T, \Gamma \neq 0}^{N} \right) \notag
\end{align}
 is of codimension $\leq g + (g-N) = 2g - N$.
That is to say, 
it has  dimension $\geq 3g-2 + \frac{2g-2}{p} - (2g-N) = g-2 + \frac{2g-2}{p-1}$.
Since the natural projections in the diagram
\begin{align}
\begin{CD}
T \times_{\mcL_Y, \mfP \mfi \mfc_{Y/T}^{N(p-1)} } \mfP \mfi \mfc_{Y/T, \Gamma \neq 0}^{N} @>>> T \times_{[\mcL_Y], \mr{Pic}_{Y/T}^{N(p-1)} } \mr{Pic}_{Y/T, \Gamma \neq 0}^{N} \\
@VVV @.
\\
 \mfT \mfa \mfn_g \times_{\mfP \mfi \mfc_g^{N(p-1)}, \mu_{N, p-1}} \mfP \mfi \mfc_{Y/T, \Gamma \neq 0}^{N} @.
\end{CD}
\end{align}
are  \'{e}tale and surjective, $ \mfT \mfa \mfn_g \times_{\mfP \mfi \mfc_g^{N(p-1)}, \mu_{N, p-1}} \mfP \mfi \mfc_{Y/T, \Gamma \neq 0}^{N}$ turns out to be of 
dimension $\geq g-2 + \frac{2g-2}{p-1}$.
This completes the proof of Theorem \ref{P0fgh7}
\end{proof}
\vspace{3mm}

\begin{rema} \label{Rrg0402}
\leavevmode\\
 \ \ \ Let us keep the notation  and assumption in Theorem \ref{P0fgh7}.
Denote by $\mfT_g$ the locus of $\mfM_g$ classifying curves   which  admits  a  Tango structure of the form $\mcL = \mcN^{\otimes (p-1)}$  for some line bundle $\mcN$ having a nontrivial global section.
According to ~\cite{Tak3}, Theorem 4.1, 
 $\mfT_g$ contains a variety of dimension $\geq g-1$.
 Moreover,   Tsuda's method (cf. ~\cite{Tsu} or  ~\cite{Tak3}, Remark 4.2) gives  a slightly better estimation: $\mr{dim} (\mfT_g) \geq 2g- \frac{(g-1) (p-1)}{p}$.
According to our result, we can prove a lower bound estimation sharper than 
  the bounds obtained previously.
 Indeed, 
 since  the stack-theoretic image of the projection $\mfG \mfT \mfa \mfn_{g, \Gamma \neq 0}^1 \migi \mfM_g$ coincides with $\mfT_g$, 
 Theorem \ref{P0fgh7} implies that
 \begin{align}
 \mr{dim} (\mfT_g) \geq g-2 + \frac{2g-2}{p-1}.
 \end{align}
\end{rema}
\vspace{3mm}

Finally, we shall conclude the paper with the following corollary (i.e., Theorem \ref{ThBf}), which follows immediately from
Theorem \ref{P0fgh7} and the construction of generalized Raynaud surfaces discussed at the beginning of \S\,\ref{wwww}.

\vspace{3mm}
\bco \label{P0ffg7} \leavevmode\\
 \ \ \ 
Let us keep the assumption in Theorem \ref{P0fgh7}.
Then,
there exists a flat family $\mfY \migi \mfT$ (i.e., $\mfP \migi \mfG \mfT \mfa \mfn^1_{g, \Gamma \neq 0}$ displayed  in (\ref{sssqw}))  of proper smooth algebraic surfaces of general type parametrized by a Deligne-Mumford stack $\mfT$  of dimension $\geq g-2 + \frac{2g-2}{p-1}$ all of whose  fibers are  pairwise non-isomorphic and have the automorphism group schemes being   non-reduced.
  \eco
\begin{proof}
The assertion follows from Theorem \ref{P0fgh7}, the discussion preceding  Theorem \ref{P0fgh7}, and the following lemma. 
\end{proof}

\vspace{3mm}
\ble \label{P0hjki55} \leavevmode\\
 \ \ \ 
 Let $\widehat{\mfX} := (X, \mcL, \mcN, \nu)$ and $\widehat{\mfX}' := (X', \mcL', \mcN', \nu')$ be generalized Tango curves of index $1$ over $k$.
 Denote by $P$ and $P'$ the generalized Raynaud surfaces associated with  $\widehat{\mfX}$ and $\widehat{\mfX}'$
 respectively.
 Suppose that $P$ is isomorphic (as a scheme  over $k$) to $P'$.
 Then, $\widehat{\mfX}$ is isomorphic to $\widehat{\mfX}'$.
  \ele
\begin{proof}
Denote by $\Psi : P \migi X^{(1)}_k$ and $\Psi' : P' \migi X'^{(1)}_k$
the natural projections of $P$ and $P'$ respectively. 
Let us fix an isomorphism $h_P : P \isom P'$.
Recall that  both $X^{(1)}_k$ and  $X'^{(1)}_k$ are of genus $g >1$, and that
the  fibers of $\Psi : P \migi X^{(1)}_k$ and $\Psi' : P' \migi X'^{(1)}_k$ are rational (cf. ~\cite{Tak3}, Theorem 3.1).
 This implies that  $h_P$ maps each   fiber of $\Psi$ to a fiber of $\Psi'$.
 Hence, $h_P$ determines  a homeomorphism $| h_{X^{(1)}_k}| : | X^{(1)}_k | \isom | X'^{(1)}_k |$  between the underlying topological spaces of $X^{(1)}_k$ and $X'^{(1)}_k$
which  is  compatible, in a natural sense, with (the underlying homeomorphism of) $h_P$ via $\Psi$ and $\Psi'$.
Since $\Psi$ and $\Psi'$ induce isomorphisms  $\mcO_{X^{(1)}_k} \isom \Psi_* (\mcO_P)$ and $\mcO_{X'^{(1)}_k} \isom  \Psi'_* (\mcO_{P'})$ respectively, $| h_{X^{(1)}_k}|$  extends to an isomorphism $h_{X^{(1)}_k} : X^{(1)}_k \isom X'^{(1)}_k$ of $k$-schemes 
which makes the following square diagram commute:
\begin{align}
\begin{CD}
P @> h_P > \sim > P'
\\
@V \Psi V \wr V @V \wr V \Psi' V
\\
X^{(1)}_k @>\sim > h_{X^{(1)}_k} > X'^{(1)}_k.
\end{CD}
\end{align}
As  $k$ is 
 algebraically closed, one may find  an isomorphism $h_X : X \isom X'$ which induces $h_{X^{(1)}_k}$, via base-change  by $F_{\mr{Spec}(k)}$.
Moreover, the isomorphism $h_X$ induces an isomorphism
\begin{align} \label{ee450}
h^*_{X^{(1)}_k} (F_{X'/k*}(\Omega_{X'/k})) \ \left(\cong F_{X/k*}(h^*_X (\Omega_{X'/k}))\right) \isom F_{X/k*}(\Omega_{X/k}).
\end{align}

Next, denote by $P^{\mr{sm}}$ and $P'^{\mr{sm}}$  the smooth loci in $P$ and $P'$ respectively  relative to $X^{(1)}_k$ and $X'^{(1)}_k$ respectively.
 $h_P$ restricts to an isomorphism 
$h_X^{\mr{sm}} : P^{\mr{sm}} \isom P'^{\mr{sm}}$.
It follows from  ~\cite{Lang}, Lemma 1, (and the fact that both $P^{\mr{sm}}$ and $P'^{\mr{sm}}$ are relative affine spaces) that  the projection $\Psi |_{P^{\mr{sm}}} : P^{\mr{sm}} \migi X^{(1)}_k$ (resp.,  $\Psi' |_{P'^{\mr{sm}}} : P'^{\mr{sm}} \migi X'^{(1)}_k$) admits a  global section, and that   the normal bundle of any global section  of this projection  is isomorphic to $\mcN$ (resp., $\mcN'$).
Hence, by passing to $h_P$ and $h_{X^{(1)}_k}$, we obtain an isomorphism 
$h_\mcN : h^*_{X^{(1)}_k} (\mcN') \isom \mcN$.
By the fact discussed in Remark \ref{RR44}, the images of the two composite injections
\begin{align}
h^*_{X^{(1)}_k} (\mcN')^{\otimes (p-1)} \xrightarrow{h_\mcN^{\otimes (p-1)}} \mcN^{\otimes (p-1)}
\xrightarrow{\nu'} \mcL \xrightarrow{\mr{incl.}} F_{X/k*}(\Omega_{X/k})
\end{align}
and 
\begin{align} \label{ee420}
h^*_{X^{(1)}_k} (\mcN')^{\otimes (p-1)}  \xrightarrow{h^*_{X^{(1)}_k} (\nu)} h^*_{X^{(1)}_k} (\mcL') \xrightarrow{\mr{incl.}} h^*_{X^{(1)}_k} (F_{X'/k*}(\Omega_{X'/k})) 
 \xrightarrow{(\ref{ee450})} F_{X/k*}(\Omega_{X/k})
\end{align}
specify the same Tango structure on $X/k$.
Thus, after  possibly replacing $h_\mcN$ with its composite with the automorphism of $\mcN$ (given by  multiplication by some element of $k^{\times}$), 
$h_\mcN$ makes the following square diagram  commute:
\begin{align}
\begin{CD}
h^*_{X^{(1)}_k} (\mcN')^{\otimes (p-1)} @> h_\mcN^{\otimes (p-1)} > \sim > \mcN^{\otimes (p-1)}
\\
@Vh^*_{X^{(1)}_k} (\nu) V \wr V @V \wr V \nu' V
\\
h^*_{X^{(1)}_k} (\mcL') @> \sim >> \mcL,
\end{CD}
\end{align} 
where the lower horizontal arrow is obtained by restricting 
(\ref{ee450}).
Consequently, the pair $(h_{X}, h_\mcN)$ specifies an isomorphism $\widehat{\mfX} \isom \widehat{\mfX}'$ of generalized Tango curves.
This completes the proof of Lemma \ref{P0hjki55}.
\end{proof}
\vspace{3mm}

\end{document}